\documentclass[12pt,oneside,reqno]{amsart}
\usepackage{graphicx}
\usepackage{mathrsfs}
\usepackage{stmaryrd}
\usepackage{amsfonts}
\usepackage{cite}
\usepackage{enumerate,amsmath,amssymb,amsthm}
\usepackage{booktabs} 
\usepackage{diagbox} 

\newcommand{\dif}{\mathrm{d}}

\newcommand{\be}{\begin{eqnarray}}
	\newcommand{\ee}{\end{eqnarray}}
\newcommand{\ce}{\begin{eqnarray*}}
	\newcommand{\de}{\end{eqnarray*}}
\newtheorem{theorem}{Theorem}[section]
\newtheorem{lemma}[theorem]{Lemma}
\newtheorem{remark}[theorem]{Remark}
\newtheorem{definition}[theorem]{Definition}
\newtheorem{proposition}[theorem]{Proposition}
\newtheorem{Examples}[theorem]{Examples}
\newtheorem{corollary}[theorem]{Corollary}
\newtheorem{condition}[theorem]{Condition}
\def\e{\varepsilon}
\def\t{\theta}
\def\a{\alpha}

\def\d{\delta}
\def\p{\partial}
\def\g{\gamma}

\def\l{\lambda}

\def\[{{\Big[}}
\def\]{{\Big]}}
\def\<{{\langle}}
\def\>{{\rangle}}
\def\({{\Big(}}
\def\){{\Big)}}

\def\no{\nonumber}
\def\bt{\begin{theorem}}
	\def\et{\end{theorem}}
\def\bl{\begin{lemma}}
	\def\el{\end{lemma}}
\def\br{\begin{remark}}
	\def\er{\end{remark}}
\def\bx{\begin{Examples}}
	\def\ex{\end{Examples}}
\def\bd{\begin{definition}}
	\def\ed{\end{definition}}
\def\bp{\begin{proposition}}
	\def\ep{\end{proposition}}
\def\bc{\begin{corollary}}
	\def\ec{\end{corollary}}
\def\bco{\begin{condition}}
	\def\eco{\end{condition}}

\def\cH{{\mathcal H}}

\def\cK{{\mathcal K}}
\def\cL{{\mathcal L}}

\def\cT{{\mathcal T}}

\def\cX{{\mathcal X}}

\def\mE{{\mathbb E}}

\def\mH{{\mathbb H}}

\def\mN{{\mathbb N}}

\def\mP{{\mathbb P}}

\def\mR{{\mathbb R}}
\def\mS{{\mathbb S}}

\def\mU{{\mathbb U}}
\def\mV{{\mathbb V}}

\def\mX{{\mathbb X}}

\def\sF{{\mathscr F}}

\def\sP{{\mathscr P}}

\def\geq{\geqslant}
\def\leq{\leqslant}

\begin{document}
	
\allowdisplaybreaks
\title{Averaging principles for nonautonomous multiscale stochastic Burgers equations with reflection}
	
\author{Huijie Qiao}

	
	
	
\date{}
	
\dedicatory{School of Mathematics,
		Southeast University\\
		Nanjing, Jiangsu 211189, China\\
		hjqiaogean@seu.edu.cn}
	
\begin{abstract}
This paper investigates averaging principles for nonautonomous multiscale stochastic Burgers equations with reflection. First, we establish a general averaging principle valid under minimal assumptions on the coefficients. Next, recognizing that the coefficients of the initially averaged equation retain explicit dependence on the small scaling parameter $\e$, we impose three distinct structural conditions--namely, asymptotic, periodic, and almost periodic dependence--on the coefficients. Under each condition, we derive a corresponding $\e$-independent averaged equation and rigorously prove the associated averaging principle in the probability sense. As a by-product, we further obtain an averaging principle in the $L^2$ norm for stochastic Burgers equations with reflection. Finally, a concrete example is provided to illustrate the applicability and validity of the theoretical results.
\end{abstract}
	
\maketitle \rm
	
\section{Introduction}

Consider the following stochastic slow-fast system on the spatial domain $[0,1]$:
\be\left\{\begin{array}{ll}
\frac{\p X_t^\e(\xi)}{\p t}=\frac{\p^2 X_t^\e(\xi)}{\p\xi^2}+f_1(X_t^\e,Y_t^\e)(\xi)+g_1(X_t^\e)\frac{\p W^1(t,\xi)}{\p t},\\
\frac{\p Y_t^\e(\xi)}{\p t}=\frac{1}{\e}\left[\frac{\p^2 Y_t^\e(\xi)}{\p\xi^2}+f_2(X_t^\e,Y_t^\e)(\xi)\right]+\frac{1}{\sqrt \e} g_2(X_t^\e,Y_t^\e)\frac{\p W^2(t,\xi)}{\p t},\\
X_t^\e(0)=X_t^\e(1)=Y_t^\e(0)=Y_t^\e(1)=0, \quad t\geq 0,\\
X_0^\e=x_0, \quad Y_0^\e=y_0,
\end{array}
\right.
\label{00eq}
\ee
where $\e>0$ represents the time-scale separation ratio between the slow component $X^\e$ and the fast component $Y^\e$. The nonlinearities $f_1$, $g_1$, $f_2$, $g_2$ satisfy regularity conditions and $W^1$, $W^2$ are mutually independent cylindrical Wiener processes defined on a complete filtered probability space $(\Omega, \sF, \{\sF_t\}_{t\geq 0}, \mP)$. These systems like (\ref{00eq}) appear to describe diffusive phenomena in reactive media, such as combustion, epidemic propagation and diffusive transport of chemical species through cells and dynamics of populations (\cite{sC1, sC2}). In \cite{sC1}, Cerrai first established an averaging principle for system \eqref{00eq} under Lipschitz continuity assumptions, proving that the slow component $X^\varepsilon$ converges in probability (as $\varepsilon \to 0$) to the solution $\tilde{X}$ of the averaged stochastic partial differential equation (SPDE):
\ce\left\{\begin{array}{ll}
\frac{\p\tilde X_t(\xi)}{\p t}=\frac{\p^2 \tilde X_t(\xi)}{\p\xi^2}+\tilde f_1(\tilde X_t)(\xi)+g_1(\tilde X_t)\frac{\p W^1(t,\xi)}{\p t},\\
\tilde X_0=x_0,
\end{array}
\right.
\de
where $\tilde f_1(x):=\int_\mH f_1(x,y)\nu^x(\dif y)$, and $\nu^x$ is the invariant measure of the following SPDE:
\be\left\{\begin{array}{ll}
\frac{\p Y_t^{x,y_0}(\xi)}{\p t}=\left[\frac{\p^2 Y_t^{x,y_0}(\xi)}{\p\xi^2}+f_2(x,Y_t^{x,y_0})(\xi)\right]+g_2(x,Y_t^{x,y_0})\frac{\p W^2(t,\xi)}{\p t},\\
Y_t^{x,y_0}(0)=Y_t^{x,y_0}(1)=0, \quad t\geq 0,\\
Y_0^{x,y_0}=y_0.
\end{array}
\right.
\label{00eqfroz}
\ee
Later, Cerrai \cite{sC2} extended the above result to more general systems with polynomially growing coefficients. Now, about averaging principles for system (\ref{00eq}) there have been many related results (cf. \cite{BYY, eB, CSS, ChLi, DSXZ, Gp1, Gp2, Gp3, GS, HLYY, HLL, LRZ, LRSX, LLS, MY, Q2, RXY, TWY, YXJ}).

Next, we consider the case where the coefficients of the fast component in system \eqref{00eq} explicitly depend on time, specifically through a rapidly oscillating temporal scale $t/\e$. This introduces significant technical challenges to the averaging analysis, as the associated frozen dynamics become nonautonomous. To illustrate this setting, we examine the following slow-fast SPDE:
\be\left\{\begin{array}{ll}
\frac{\p X_t^\e(\xi)}{\p t}=\frac{\p^2 X_t^\e(\xi)}{\p\xi^2}+f_1(X_t^\e,Y_t^\e)(\xi)+g_1(X_t^\e)\frac{\p W^1(t,\xi)}{\p t},\\
\frac{\p Y_t^\e(\xi)}{\p t}=\frac{1}{\e}\left[\frac{\p^2 Y_t^\e(\xi)}{\p\xi^2}+f_2(\frac{t}{\e},X_t^\e,Y_t^\e)(\xi)\right]+\frac{1}{\sqrt \e} g_2(\frac{t}{\e},X_t^\e,Y_t^\e)\frac{\p W^2(t,\xi)}{\p t},\\
X_t^\e(0)=X_t^\e(1)=Y_t^\e(0)=Y_t^\e(1)=0, \quad t\geq 0,\\
X_0^\e(\xi)=x_0(\xi), \quad Y_0^\e(\xi)=y_0(\xi).
\end{array}
\right.
\label{01eq}
\ee
Such systems arise in diverse applied contexts, including mathematical biology and statistical physics. A prominent example is the Hodgkin-Huxley model and its reduced variants, which describe the time-dependent activation and inactivation kinetics of voltage-gated ion channels in neuronal membranes (\cite{hT}). For system \eqref{01eq}, the standard frozen equation \eqref{00eqfroz} generalizes to the following nonautonomous SPDE:
\be\left\{\begin{array}{ll}
\frac{\p Y_t^{x,y_0}(\xi)}{\p t}=\left[\frac{\p^2 Y_t^{x,y_0}(\xi)}{\p\xi^2}+f_2(t,x,Y_t^{x,y_0})(\xi)\right]+g_2(t,x,Y_t^{x,y_0})\frac{\p W^2(t,\xi)}{\p t},\\
Y_t^{x,y_0}(0)=Y_t^{x,y_0}(1)=0, \quad t\geq 0,\\
Y_0^{x,y_0}=y_0.
\end{array}
\right.
\label{01eqfroz}
\ee
Crucially, due to the explicit time dependence, Eq.\eqref{01eqfroz} does not admit a unique invariant measure. Consequently, the classical averaging procedure, relying on ergodicity, fails. To overcome this, Cerrai and Lunardi \cite{CL} constructed an evolution family of probability measures $\{\mu_t^x\}_{t \in \mathbb{R}}$ associated with Eq.\eqref{01eqfroz} (See Subsection \ref{evolsyst}) and, under almost periodic assumptions on the coefficients, established convergence in probability: as $\varepsilon \to 0$, the slow component $X^\varepsilon$ converges to the solution $\tilde{\tilde{X}}$ of the averaged SPDE
\ce\left\{\begin{array}{ll}
\frac{\p \tilde{\tilde X}_t(\xi)}{\p t}=\frac{\p^2 \tilde{\tilde X}_t(\xi)}{\p\xi^2}+\tilde f_{1,ap}(\tilde{\tilde X}_t)(\xi)+g_1(\tilde{\tilde X}_t)\frac{\p W^1(t,\xi)}{\p t},\\
\tilde{\tilde X}_0=x_0,
\end{array}
\right.
\de
with the almost periodic averaged drift defined by
$$
\tilde f_{1,ap}(x)=\lim\limits_{T\rightarrow\infty}\frac{1}{T}\int_0^T\int_{C([0,1])}f_1(x,y)\mu^x_t(\dif y)\dif t.
$$ 
Subsequently, Cheng, Sun, and Xie \cite{CSX} extended this framework by allowing $f_1$ and $g_1$ to be time-dependent and assuming uniform almost periodicity of $f_1$, $f_2$, and $g_2$. They proved an $L^p$-averaging principle for $p \geq 2$. More recently, Li et al. \cite{LSWX} replaced the cylindrical Wiener processes $W^1$ and $W^2$ in system \eqref{01eq} with independent cylindrical $\alpha$-stable L\'evy processes ($1 < \alpha < 2$), imposed time-periodicity and asymptotic convergence conditions, assumed $f_1$ is time-dependent, and took $g_1$ and $g_2$ as identity operators, thereby establishing two distinct $L^p$-averaging principles for $1< p < \alpha$. Note that the above mentioned \cite{CL, CSX, LSWX} are for nonautonomous infinite dimensional slow-fast systems. About averaging principles for nonautonomous finite dimensional slow-fast systems there are also related results (cf. \cite{HWX, LZ, SLG, SWX1, SWX2, WXXY, XQ}).

In what follows, we address the following natural question: under what conditions on the slow component $X^\varepsilon$ does the averaging principle established in \cite{CL} remain valid, particularly when the slow dynamics are subject to pathwise constraints? In this paper, we provide an affirmative answer. Specifically, we study the following slow-fast stochastic system with reflection on the spatial domain $[0,1]$:
\be\left\{\begin{array}{ll}
\frac{\p X_t^\e(\xi)}{\p t}=\frac{\p^2 X_t^\e(\xi)}{\p\xi^2}+\frac{1}{2}\frac{\p(X_t^\e)^2(\xi)}{\p \xi}+f_1(\frac{t}{\e},X_t^\e,Y_t^\e)(\xi)+g_1(\frac{t}{\e},X_t^\e)\frac{\p W^1(t,\xi)}{\p t}\\
\qquad\qquad\qquad +\frac{\p K^\e(t,\xi)}{\p t},\\
\frac{\p Y_t^\e(\xi)}{\p t}=\frac{1}{\e}\left[\frac{\p^2 Y_t^\e(\xi)}{\p\xi^2}+f_2(\frac{t}{\e},X_t^\e,Y_t^\e)(\xi)\right]+\frac{1}{\sqrt \e} g_2(\frac{t}{\e},X_t^\e,Y_t^\e)\frac{\p W^2(t,\xi)}{\p t},\\
X_t^\e(\xi)\geq 0, \\
X_t^\e(0)=X_t^\e(1)=Y_t^\e(0)=Y_t^\e(1)=0, \quad t\geq 0,\\
X_0^\e(\xi)=x_0(\xi)\geq 0, \quad Y_0^\e(\xi)=y_0(\xi),
\end{array}
\right.
\label{0eq}
\ee
where $\frac{\p K^\e(t,\xi)}{\p t}$ denotes the formal time derivative of a reflection random measure $K^\varepsilon$, which ensures the almost-sure nonnegativity of $X^\varepsilon$ (See Subsection \ref{msbewr} for its rigorous definition). Since the slow equation is a stochastic Burgers equation with reflection, we refer to system \eqref{0eq} as a nonautonomous multiscale stochastic Burgers equation with reflection. Such systems arise in continuum mechanics, e.g., modeling interface motion under stochastic forcing, and in stochastic control problems with state constraints. They have consequently drawn sustained interest in both applied probability and SPDE theory. In this paper our primary objective is to establish rigorous averaging principles for system \eqref{0eq}. 

It is noteworthy that in \cite{Q2}, we established an averaging principle for stochastic Burgers equations with reflection. Extending this result to nonautonomous multiscale stochastic Burgers equations with reflection, however, entails substantial technical difficulties. First, moment estimates in the multiscale setting become significantly more involved due to the coupling between slow and fast dynamics. Second, as the fast component is governed by a nonautonomous SPDE, the averaged equation depends critically on the invariant measure of the corresponding frozen SPDE, i.e., the fast subsystem with $\varepsilon = 1$. Consequently, distinct properties of this invariant measure yield distinct averaged equations. Furthermore, compared with existing results on nonautonomous multiscale SPDEs without reflection \cite{CL, CSX}, the presence of the reflection term $K^\varepsilon$ introduces two principal challenges. The first concerns establishing uniform-in-$\varepsilon$ temporal modulus of continuity estimates, specifically, bounds of the form $| X_{t+r}^\varepsilon -X_t^\varepsilon|_{\mH}^2$ for arbitrary $t \geq 0$ and $r > 0$. In contrast to the penalization method adopted in \cite{Q2}, we overcome this difficulty by integrating It\^o's formula with a carefully designed stopping-time argument, thereby rendering the analysis more transparent and efficient. The second challenge lies in deriving high-order moment bounds, e.g., $\sup\limits_{t \in [0,T]} \mathbb{E}|X_t^\varepsilon|_{\mH}^p$ for $p > 2$. Here again, stopping-time truncation proves instrumental in closing the requisite estimates. Importantly, our results are not subsumed by prior work (cf. \cite{CL, CSX}), even in the absence of reflection, since the drift coefficient of the slow equation in our framework satisfies a fully local monotonicity condition, and our structural assumptions are strictly more general.

Finally, we summarize our main contributions. First, by leveraging the evolution family of measures $\{\mu_t^x\}_{t \in \mathbb{R}}$ associated with the frozen equation \eqref{01eqfroz}, we construct a well-defined averaged drift coefficient and establish a general averaging principle via a refined Khasminskii-type time discretization scheme under minimal regularity and ergodicity-type assumptions, without imposing `a priori' asymptotic behavior or periodicity. Second, since the initially constructed averaged coefficients retain explicit $\varepsilon$-dependence, we introduce three distinct structural conditions on the slow and fast coefficients: asymptotic convergence as $t \to\infty$, time-periodicity and time-almost-periodicity. Under each condition, the limiting averaged coefficients become fully autonomous (i.e., independent of both $\varepsilon$ and $t$), yielding three rigorously justified averaging principles with autonomous effective dynamics. As a by-product, we further obtain an averaging principle in the $L^2$ norm for stochastic Burgers equations with reflection (See Remark \ref{l2normaverprin}). 

The rest of this paper is organized as follows. In Section \ref{pre}, we introduce some notations and concepts. In Section \ref{main} the main results are stated. Then the proofs of four main results are placed in Section \ref{genestroapthproo}, Section \ref{asymstroapthproo}, Section \ref{peristroapthproo} and Section \ref{aperistroapthproo}, respectively. Finally, in Section \ref{exam} we give an example to explain the applicability and validity of our results.

The following convention will be used throughout the paper: $C$ with or without indices will denote different positive constants whose values may change from one place to another.

\section{Preliminaries}\label{pre}

In this section, we introduce some notations and concepts.

\subsection{Notation}

In this subsection, we introduce some notations used in the sequel.

Let $C_c([0,1])$ be the set of all continuous functions on $[0,1]$ with compact supports. Let $C^2_c([0,1])$ and  $C^\infty_c([0,1])$  be the subsets of $C_c([0,1])$ where all functions have $2$ and infinite order derivatives, respectively. 

Let $\mH=L^2([0,1], \mR)$ be the usual $L^2$-space with the norm $|\cdot|_\mH$ and inner product $\<\cdot, \cdot\>_\mH$. Denote by $\mV$ the Sobolev space of order one, i.e., $\mV$ is the completion of $C_c^{\infty}([0,1])$ under the norm $\|u\|_\mV^2=\int_0^1\left(\frac{\partial u}{\partial x}\right)^2 d x$. Remark that $\mV=\left\{u \in \mH^{1}([0,1]): u(0)=u(1)=0\right\}$, where $\mH^{1}([0,1])$ denotes the usual Sobolev space of absolutely continuous functions defined on $[0,1]$ whose derivatives belong to $\mH$. $\mV^*$ denotes the dual space of $\mV$ and the dualization between $\mV$ and $\mV^*$ is denoted by ${_{\mV^*}}\<\cdot,\cdot\>_\mV$. 

Let $\mU_i$ be the Hilbert space with norm $|\cdot|_{\mU_i}$ and inner product $\<\cdot, \cdot\>_{\mU_i}, i=1,2$. Let $\cL_2(\mU_i,\mH)$ be the collection of all Hilbert-Schmidt operators from $\mU_i$ to $\mH$ with the Hilbert-Schmidt norm $\|\cdot\|_{\cL_2(\mU_i,\mH)}, i=1,2$.

\subsection{Multiscale stochastic Burgers equations with reflection}\label{msbewr}

In this subsection, we introduce the definition of the solutions for multiscale stochastic Burgers equations with reflection.

Let $A$ be the Laplace operator on $\mathbb{H}$
$$
A x:=\frac{\partial^2  x(\xi)}{\partial \xi^2}, \quad x \in D(A)=\mH^2([0,1]) \cap \mV,
$$
where $\mH^2([0,1])$ stands for the Sobolev space of functions defined on $[0,1]$ whose derivatives up to order $2$ belong to $\mH$. Let $\left\{e_k(\xi):=\sqrt{2} \sin (k \pi \xi)\right\}_{k \geq 1}$ be an orthonormal basis of $\mathbb{H}$ consisting of the eigenvectors of $A$, i.e.
$$
A e_k=-\lambda_k e_k \quad \text { with } \lambda_k=k^2 \pi^2.
$$

Define the bilinear operator $B(x, y): \mathbb{V} \times \mathbb{V} \rightarrow \mathbb{V}^*$ by
$$
B(x, y):=x \cdot \partial_{\xi} y,
$$
and the trilinear operator $b(x,y,z): \mV\times\mV\times\mV \rightarrow\mR$ by
$$
b(x,y,z):=\int_0^1x(\xi)\p_\xi y(\xi)z(\xi)\dif \xi={_{\mV^*}}\<B(x, y),z\>_\mV.
$$
Then $B$ and $b$ have the following properties (cf. \cite{BZ}).

\bl\label{Bbprop1}
$(i)$ For all $x, y, z \in \mV$,
$$
{_{\mV^*}}\<B(x, y), z\>_\mV=b(x,y,z)=-b(x,z,y)=-{_{\mV^*}}\< B(x, z), y\>_\mV.
$$

$(ii)$ For all $x, y, z \in \mV$,
$$
|{_{\mV^*}}\<B(x,y), z\>_\mV|=|b(x, y, z)| \leq 2\|x\|_\mV^{\frac{1}{2}}|x|_\mH^{\frac{1}{2}}\|z\|_{\mV}^{\frac{1}{2}}|z|_\mH^{\frac{1}{2}}\|y\|_\mV.
$$
\el

By the above lemma, we can obtain the following corollary.

\bc\label{Bbprop2}
\ce
&&b(x, y, y)=0, \quad {_{\mV^*}}\< B(x, y),y\>_\mV=0, \quad x, y \in \mV, \\
&&\|B(x, x)\|_{V^*} \leq 2\|x\|_\mV|x|_\mH, \quad x \in \mV.
\de
\ec

\bl\label{Bbprop3}
For any $x,y\in\mV$, it holds that
\ce
\|B(x,x)-B(y,y)\|_{\mV^*}\leq 2|x-y|_\mH(\|x\|_\mV+\|y\|_\mV).
\de
\el

With the above notations, the system (\ref{0eq}) can be rewritten as 
\be\left\{\begin{array}{ll}
\dif X_t^\e=[AX_t^\e+B(X_t^\e,X_t^\e)+f_1(\frac{t}{\e},X_t^\e,Y_t^\e)]\dif t+g_1(\frac{t}{\e},X_t^\e)\dif W^1_t+\dif K^\e_t,\\
\dif Y_t^\e=\frac{1}{\e}[AY_t^\e+f_2(\frac{t}{\e},X_t^\e,Y_t^\e)]\dif t+\frac{1}{\sqrt\e} g_2(\frac{t}{\e},X_t^\e,Y_t^\e)\dif W^2_t,\\
X_t^\e(\xi)\geq 0, \quad \xi\in[0,1],\\
X_t^\e(0)=X_t^\e(1)=Y_t^\e(0)=Y_t^\e(1)=0,\quad t\geq 0,\\
X_0^\e=x_0, \quad Y_0^\e=y_0,
\end{array}
\right.
\label{1eq}
\ee
where 
$$
A: \mV\to \mV^*, \quad B: \mV \times \mathbb{V} \rightarrow \mathbb{V}^*, \quad f_1: \mR_+\times\mH\times\mH\to\mH, \quad g_1: \mR_+\times\mH\to\cL_2(\mU_1,\mH),
$$
and
$$
f_2: \mR\times\mH\times\mH\to\mH, \quad g_2: \mR\times\mH\times\mH\to\cL_2(\mU_2,\mH),
$$
are measurable mappings, $W^1$ and $W^2$ are mutually independent $\mU_1$ and $\mU_2$-valued cylindrical Wiener processes on a complete filtered probability space $(\Omega, \sF, \{\sF_t\}_{t\geq 0}, \mP)$ respectively.

In the following, we present the definition of solutions for the system (\ref{1eq}) (cf. \cite{tZ}).

\bd\label{soludefi}
$(X^\e, K^\e, Y^\e)$ is said to be a solution of the system (\ref{1eq}) if

$(i)$ $X_t^\e, Y_t^\e$ are $\mV$-valued $\sF_t$-measurable for any $t\geq 0$ and $X_t^\e(\xi)\geq 0$ a.e. for any $(t,\xi)\in\mR_+\times[0,1]$;

$(ii)$ $K^\e$ is a random measure on $\mR_+\times[0,1]$ such that 

$(a)$ $\mE\left[({\rm Var}(K^\e)([0,T]\times [0,1]))^2\right]<+\infty, \forall T \geq 0$, where ${\rm Var}(K^\e)([0,T]\times[0,1])$ denotes the total variation of $K^\e$ on $[0, T]\times[0,1]$ defined by 
\ce
{\rm Var}(K^\e)([0,T]\times[0,1]):=\sup\limits_{\pi}\sum_{i=1}^n |K^\e(E_i)|,
\de
and the supremum is taken over all partitions $\pi$ of the domain $[0,T]\times[0,1]$,

$(b)$ $K^\e$ is adapted in the sense that for any bounded measurable mapping $\varphi$:
\ce
\int_0^t\int_0^1\varphi(r,\xi)K^\e(\dif r, \dif \xi) ~\mbox{is}~\sF_t\mbox{-measurable};
\de

$(iii)$ $(X^\e, K^\e, Y^\e)$ satisfies the parabolic stochastic partial differential equations with reflection in the following sense: $\forall t \in \mathbb{R}_{+}, \psi \in C_c^2([0,1])$ with $\psi(0)=\psi(1)=0$,
\ce
&&\<X_t^\e,\psi\>_\mH=\<x_0,\psi\>_\mH+\int_0^t\[{_{\mV^*}}\<AX_r^\e+B(X_r^\e,X_r^\e),\psi\>_\mV+\<f_1(\frac r \e, X_r^\e,Y_r^\e),\psi\>_\mH\]\dif r\\
&&\qquad\qquad\qquad+\int_0^t \<\psi, g_1(\frac r \e, X_r^\e)\dif W^1_r\>_\mH+\int_0^t\int_0^1\psi(\xi) \dif K^\e(\dif r,\dif \xi), ~a.s.,\\
&&\<Y_t^\e,\psi\>_\mH=\<y_0,\psi\>_\mH+\frac{1}{\e}\int_0^t\[{_{\mV^*}}\<AY_r^\e,\psi\>_\mV+\<f_2(\frac r \e, X_r^\e,Y_r^\e),\psi\>_\mH\]\dif r\\
&&\qquad\qquad\qquad+\frac{1}{\sqrt\e}\int_0^t \<\psi,g_2(\frac r \e, X_r^\e,Y_r^\e)\dif W^2_r\>_\mH, ~a.s.;
\de

$(iv)$ for any $T>0$, $\int_0^T\int_0^1X_t^\e(\xi)\dif K^\e(\dif t, \dif \xi)=0$.
\ed

\subsection{Almost periodic functions}

In this subsection, we introduce almost periodic functions (cf. \cite{CL, ChLi, DaT, aF}).

Let $(\mX, \rho)$ be a complete separable metric space. We denote by $\sP(\mX)$ the set of all probability measures onto the Borel $\sigma$-field of $\mX$. We shall denote by $C(\mR, \mX)$ the class of all continuous functions from $\mR$ to $\mX$, and by $C_b(\mX)$ the class of all continuous and bounded functions $\phi: \mX\to\mR$ with $\|\phi\|_{\infty}:=\sup\limits_{x \in\mX}|\phi(x)|<\infty$.

For $\phi\in C_b(\mX)$ we define
\ce
&&\|\phi\|_{Lip}:=\sup\limits_{x \neq y, x,y\in\mX}\frac{|\phi(x)-\phi(y)|}{\rho(x, y)}, \\
&&\|\phi\|_{BL}:=\max \left\{\|\phi\|_{\infty},\|\phi\|_{Lip}\right\}.
\de
Set for $\mu, \nu \in \sP(\mX)$
$$
d_{BL}(\mu, \nu):=\sup\limits_{\|\phi\|_{B L} \leq 1}\left|\int_\mX \phi \dif (\mu-\nu)\right|,
$$
and $d_{BL}$ is a complete metric on $\sP(\mX)$ which generates the weak convergence topology. From now on $\sP(\mX)$ is endowed with the metric $d_{BL}$. 

Let $(\mX_1,\rho_1)$ and $(\mX_2,\rho_2)$ be two complete separable metric spaces. 

\bd\label{apuapdefi}
$(i)$ We say $\phi\in C(\mR, \mX_2)$ is almost periodic, if for any $\g>0$, there exists a constant $\iota(\g)>0$ such that $\cT(\phi,\g)\cap [a, a+\iota]\neq \emptyset$ for all $a\in\mR$, where
\ce
\cT(\phi,\g):=\left\{\tau\in\mR: \sup\limits_{t\in\mR}\rho_2(\phi(t+\tau),\phi(t))<\g\right\}.
\de

$(ii)$ We say $\phi\in C(\mR\times\mX_1,\mX_2)$ is uniformly almost periodic, if for any $\g>0$ and every compact subset $\cX\subset\mX_1$, there exists a constant $\iota(\g, \cX)>0$ such that $\cT(\phi,\g,\cX)\cap [a, a+\iota]\neq \emptyset$ for all $a\in\mR$, where
\ce
\cT(\phi,\g,\cX):=\left\{\tau\in\mR: \sup\limits_{(t,x)\in\mR\times\cX}\rho_2(\phi(t+\tau,x),\phi(t,x))<\g\right\}.
\de
\ed

Next, we give a characterization of uniformly almost periodic functions.

\bt(\cite[Theorem 2.3]{DaT})\label{uapchartheo}
For $\phi\in C(\mR\times\mX_1,\mX_2)$, the following statements are equivalent:

$(i)$ $\phi$ is uniformly almost periodic.

$(ii)$ For any sequence $\zeta^{\prime}=\{\zeta_n^{\prime}\}_{n\in\mN}\subset\mR$, there exists a subsequence $\zeta\subset\zeta^\prime$ and a continuous function $\phi^\prime: \mR\times\mX_1\to \mX_2$ such that 
$$
\lim\limits_{n\rightarrow\infty}\phi(t+\zeta_n,x)=\phi^\prime(t,x)
$$
uniformly on $\mR\times\cX$, for any compact subset $\cX\subset\mX_1$.

$(iii)$ For every pair of sequences $\zeta^\prime=\{\zeta_n^{\prime}\}_{n\in\mN}\subset\mR$ and $\zeta^{\prime\prime}=\{\zeta^{\prime\prime}_m\}_{m\in\mN}\subset\mR$, there exist common subsequences $\zeta\subset\zeta^\prime$ and $\check\zeta\subset\zeta^{\prime\prime}$ 
such that 
$$
\lim\limits_{n\rightarrow\infty}\phi(t+\zeta_n+\check\zeta_n,x)=\lim\limits_{n\rightarrow\infty}\lim\limits_{m\rightarrow\infty}\phi(t+\zeta_n+\check\zeta_m,x)
$$
uniformly on $\mR\times\cX$, for any compact subset $\cX\subset\mX_1$.
\et

Finally, about the mean values of uniformly almost periodic functions, we have the following result. 

\bt(\cite[Theorem 3.4]{CL} or \cite[Theorem 3.1 and Corollary 3.2]{aF})\label{uapmeantheo}
If $\phi\in C(\mR\times\mX_1,\mX_2)$ is uniformly almost periodic, then
\ce
\exists \lim\limits_{T\rightarrow\infty}\frac{1}{T}\int_t^{t+T}\phi(s,x)\dif s=\lim\limits_{T\rightarrow\infty}\frac{1}{T}\int_0^{T}\phi(s,x)\dif s
\de
uniformly on $\mR\times\cX$, for any compact subset $\cX\subset\mX_1$.
\et

\section{Main results}\label{main}

In this section, we formulate the main results in this paper.

\subsection{Nonautonomous multiscale stochastic Burgers equations with reflection in the general case}\label{genecase}

In this subsection, we will present an averaging principle for nonautonomous multiscale stochastic Burgers equations with reflection in the general setting.

We recall system (\ref{1eq}), i.e.
\ce\left\{\begin{array}{ll}
\dif X_t^\e=[AX_t^\e+B(X_t^\e,X_t^\e)+f_1(\frac{t}{\e},X_t^\e,Y_t^\e)]\dif t+g_1(\frac{t}{\e},X_t^\e)\dif W^1_t+\dif K^\e_t,\\
\dif Y_t^\e=\frac{1}{\e}[AY_t^\e+f_2(\frac{t}{\e},X_t^\e,Y_t^\e)]\dif t+\frac{1}{\sqrt\e} g_2(\frac{t}{\e},X_t^\e,Y_t^\e)\dif W^2_t,\\
X_t^\e(\xi)\geq 0, \quad \xi\in[0,1],\\
X_t^\e(0)=X_t^\e(1)=Y_t^\e(0)=Y_t^\e(1)=0,\quad t\geq 0,\\
X_0^\e=x_0, \quad Y_0^\e=y_0.
\end{array}
\right.
\de
Assume:
\begin{enumerate}[$({\bf H}_{f_1, g_1})$]
\item There exist two constants $L_1, L_2\geq 0$ such that for any $t\in\mR_+, x_1, x_2, y_1, y_2\in\mH$,
\ce
&&|f_1(t,x_1,y_1)-f_1(t,x_2,y_2)|_\mH\leq L_1(|x_1-x_2|_\mH+|y_1-y_2|_\mH),\\
&&|f_1(t,0,0)|_\mH\leq L_1,\\
&&\|g_1(t,x_1)-g_1(t,x_2)\|_{\cL_2(\mU_1,\mH)}\leq L_2|x_1-x_2|_\mH,\\
&&\|g_1(t,0)\|_{\cL_2(\mU_1,\mH)}\leq L_2.
\de
\end{enumerate}
\begin{enumerate}[$({\bf H}^1_{f_2, g_2})$]
\item There exist four constants $L_3, L_4, L_5, L_6\geq 0$ such that for any $t\in\mR, x_1, x_2, y_1, y_2\in\mH$,
\ce
&&|f_2(t,x_1,y_1)-f_2(t,x_2,y_2)|_\mH\leq L_3|x_1-x_2|_\mH+L_4|y_1-y_2|_\mH,\\
&&|f_2(t,0,0)|_\mH\leq L_3,\\
&&\|g_2(t,x_1,y_1)-g_2(t,x_2,y_2)\|_{\cL_2(\mU_2,\mH)}\leq L_5|x_1-x_2|_\mH+L_6|y_1-y_2|_\mH,\\
&&\|g_2(t,0,0)\|_{\cL_2(\mU_2,\mH)}\leq L_5.
\de
\end{enumerate}
\begin{enumerate}[$({\bf H}^2_{f_2, g_2})$]
\item $L_4, L_6$ satisfy
\ce
\l_1-\frac{3}{2}L_4-L^2_6>0,
\de
where $\l_1>0$ is the smallest eigenvalue of $-A$.
\end{enumerate}

\br
$(i)$ By $({\bf H}_{f_1, g_1})$, it is easy to see that for any $t\in\mR_+, x, y\in\mH$,
\be
|f_1(t,x,y)|_\mH\leq L_1(1+|x|_\mH+|y|_\mH), \quad \|g_1(t,x)\|_{\cL_2(\mU_1,\mH)}\leq L_2(1+|x|_\mH).
\label{f1g1linegrow}
\ee

$(ii)$ $({\bf H}^1_{f_2, g_2})$ yields that for any $t\in\mR, x, y\in\mH$,
\be
|f_2(t,x,y)|_\mH\leq L_3(1+|x|_\mH)+L_4|y|_\mH, \quad \|g_2(t,x,y)\|_{\cL_2(\mU_2,\mH)}\leq L_5(1+|x|_\mH)+L_6|y|_\mH.
\label{f2g2linegrow}
\ee

$(iii)$ The Poincar\'e inequality, the Young inequality, $({\bf H}^1_{f_2, g_2})$ and $({\bf H}^2_{f_2, g_2})$ imply that for $t\in\mR, x_1, x_2, x\in\mH, y_1, y_2, y\in\mV$
\be
&&2{_{\mV^*}}\<Ay_1-Ay_2, y_1-y_2\>_\mV+2\<f_2(t,x_1,y_1)-f_2(t,x_2,y_2), y_1-y_2\>_\mH\no\\
&&+\|g_2(t,x_1,y_1)-g_2(t,x_2,y_2)\|^2_{\cL_2(\mU_2,\mH)}\no\\
&\leq&-2\l_1|y_1-y_2|^2_\mH+2L_3|x_1-x_2|_\mH|y_1-y_2|_\mH+2L_4|y_1-y_2|^2_\mH\no\\
&&+2L^2_5|x_1-x_2|^2_\mH+2L^2_6|y_1-y_2|^2_\mH\no\\
&\leq&-2\l_1|y_1-y_2|^2_\mH+3L_4|y_1-y_2|^2_\mH+2L^2_6|y_1-y_2|^2_\mH+C|x_1-x_2|^2_\mH\no\\
&=&-2\a|y_1-y_2|^2_\mH+C|x_1-x_2|^2_\mH,
\label{discond1}
\ee
and similarly
\be
&&2{_{\mV^*}}\<Ay, y\>_\mV+2\<f_2(t,x,y), y\>_\mH+\|g_2(t,x,y)\|^2_{\cL_2(\mU_2,\mH)}\no\\
&\leq& -2\a|y|^2_\mH+C(1+|x|^2_\mH),
\label{discond2}
\ee
where $\a:=\l_1-\frac{3}{2}L_4-L^2_6>0$.
\er

Based on Theorem 2.3 in \cite{LRSX}, following the idea of Theorem 2.1 in \cite{MY}, we can obtain the following well-posedness result, whose proof is omitted.

\bt\label{wellpose}
Assume that $({\bf H}_{f_1, g_1})$, $({\bf H}^1_{f_2, g_2})$ and $({\bf H}^2_{f_2, g_2})$ hold. Then for any $x_0\in\mV, x_0(\xi)\geq 0, y_0\in\mH$, system (\ref{1eq}) has a unique solution $(X^\e, K^\e, Y^\e)$ such that 
$$
\mE\left[\sup\limits_{t\in[0,T]}|X^\e_t|_\mH^2+\int_0^T\|X^\e_t\|^2_\mV\dif t\right]<\infty, \quad \mE\left[\sup\limits_{t\in[0,T]}|Y^\e_t|_\mH^2+\int_0^T\|Y^\e_t\|^2_\mV\dif t\right]<\infty.
$$
\et

Next, we define 
\ce
\bar{f}_1(t,x):=\int_\mH f_1(t,x,y)\mu^x_t(\dif y), \quad t\geq 0,
\de
where $\mu^x_t$ is from Lemma \ref{mux}, and construct the following averaged equation
\be\left\{\begin{array}{ll}
\dif \bar X_t^\e=[A\bar X_t^\e+B(\bar X_t^\e,\bar X_t^\e)+\bar f_1(\frac{t}{\e},\bar X_t^\e)]\dif t+g_1(\frac{t}{\e},\bar X_t^\e)\dif W^1_t+\dif \bar K^\e_t,\\
\bar X_t^\e(\xi)\geq 0, \quad \xi\in[0,1],\\
\bar X_t^\e(0)=\bar X_t^\e(1)=0,\quad t\geq 0,\\
\bar X_0^\e=x_0.
\end{array}
\right.
\label{eaverequa}
\ee
The following theorem describes the averaging principle for system (\ref{1eq}) in the general case.

\bt\label{genestroapth}
Suppose that $({\bf H}_{f_1, g_1})$, $({\bf H}^1_{f_2, g_2})$ and $({\bf H}^2_{f_2, g_2})$ hold. Then for any $T>0$ and $\eta>0$, 
\ce
\lim\limits_{\e\rightarrow0}\mP\left(\sup\limits_{t\in[0,T]}|X^\e_t-\bar X^\e_t|_\mH\geq \eta\right)=0,
\de
where $\bar X^\e$ is the solution of Eq.(\ref{eaverequa}).
\et

\br
It is worth noting that the coefficients in Eq.(\ref{eaverequa}) still retain explicit dependence on the small parameter $\e$. Therefore, in oder to derive averaged equations whose coefficients are independent of $\e$, we introduce asymptotic, periodic and almost periodic conditions separately in the following three subsections.
\er

The proof of the above theorem is postponed to Section \ref{genestroapthproo}.

\subsection{Nonautonomous multiscale stochastic Burgers equations with reflection in the asymptotic case}

In this subsection, we present an averaging principle for nonautonomous multiscale stochastic Burgers equations with reflection in the asymptotic case.

Following Subsection \ref{genecase}, we continue to assume:
\begin{enumerate}[$({\bf H}^a_{f_1})$]
\item There exists $\check f_1: \mH\times\mH\rightarrow\mH$ such that for any $T>0$, $x,y\in\mH$
\ce
&&\sup\limits_{t\geq 0}\left|\frac{1}{T}\int_t^{t+T}f_1(s,x,y)\dif s-\check f_1(x,y)\right|_\mH\leq\kappa_1(T)(1+|x|_\mH+|y|_\mH),\\
\de
where $\kappa_1: [0,\infty)\rightarrow(0,\infty)$ is locally bounded and satisfies $\lim\limits_{T\rightarrow\infty}\kappa_1(T)=0$.
\end{enumerate} 
\begin{enumerate}[$({\bf H}^a_{g_1})$]
\item There exists $\bar g_1: \mH\to \cL_2(\mU_1,\mH)$ such that for any $T>0$ and $x\in\mH$,
\ce
\sup\limits_{t\geq 0}\frac{1}{T}\int_t^{t+T}\left\|g_1(s,x)-\bar g_1(x)\right\|^2_{\cL_2(\mU_1,\mH)}\dif s\leq\kappa_2(T)(1+|x|^2_\mH),
\de
where $\kappa_2: [0,\infty)\rightarrow(0,\infty)$ is locally bounded and satisfies $\lim\limits_{T\rightarrow\infty}\kappa_2(T)=0$.
\end{enumerate}
\begin{enumerate}[$({\bf H}^a_{f_2, g_2})$]
\item There exist $\bar f_{2}: \mH\rightarrow\mH$ and $\bar g_{2}: \mH\rightarrow\cL_2(\mU_2,\mH)$ such that for any $t\geq 0$, $x,y\in\mH$
\ce
|f_2(t,x,y)-\bar f_{2}(x,y)|_\mH+\|g_2(t,x,y)-\bar g_{2}(x,y)\|_{\cL_2(\mU_2,\mH)}\leq \kappa_3(t)(1+|x|_\mH+|y|_\mH),
\de
where $\kappa_3: [0,\infty)\rightarrow(0,\infty)$ is locally bounded and satisfies $\lim\limits_{t\rightarrow\infty}\kappa_3(t)=0$ and 
$$
\lim\limits_{T\to\infty}\sup\limits_{t\geq 0}\frac{1}{T}\int_t^{t+T}\kappa^2_3(r)\dif r=0.
$$
\end{enumerate}

\br
By $({\bf H}^a_{f_2, g_2})$, we derive the corresponding autonomous limit frozen equation. Subsequently, by leveraging $({\bf H}^a_{f_1})$ and $({\bf H}^a_{g_1})$ we construct the associated averaged equation with coefficients that are entirely independent of the small parameter $\e$.
\er

Next, we set 
\ce
\check f_{1,a}(x):=\int_\mH \check f_1(x,y)\mu^x(\dif y),
\de
where $\mu^x$ is from Lemma \ref{mutmu}, and consider the following averaged equation:
\be\left\{\begin{array}{ll}
\dif \check X_t=[A\check X_t+B(\check X_t,\check X_t)+\check f_{1,a}(\check X_t)]\dif t+\bar g_1(\check X_t)\dif W^1_t+\dif \check K_t,\\
\check X_t(\xi)\geq 0, \quad \xi\in[0,1],\\
\check X_t(0)=\check X_t(1)=0,\quad t\geq 0,\\
\check X_0=x_0.
\end{array}
\right.
\label{averequas}
\ee
The following theorem is the second main result in this paper.

\bt\label{asymstroapth}
Suppose that $({\bf H}_{f_1, g_1})$, $({\bf H}^1_{f_2, g_2})$, $({\bf H}^2_{f_2, g_2})$, $({\bf H}^a_{f_1})$, $({\bf H}^a_{g_1})$ and $({\bf H}^a_{f_2, g_2})$ hold. Then for any $T>0$ and $\eta>0$, 
\ce
\lim\limits_{\e\rightarrow 0}\mP\left(\sup\limits_{t\in[0,T]}|X_t^\e-\check X_t|_\mH\geq\eta\right)=0,
\de
where $\check X$ is the solution of Eq.(\ref{averequas}).
\et

The proof of the above theorem is placed in Section \ref{asymstroapthproo}.

\subsection{Nonautonomous multiscale stochastic Burgers equations with reflection in the periodic case}

In this subsection, we state an averaging principle for nonautonomous multiscale stochastic Burgers equations with reflection in the periodic setting.

In order to obtain the averaged equation whose coefficients don't depend on $\e$, we assume more:
\begin{enumerate}[$({\bf H}^p_{f_1})$]
\item For any $x, y\in\mH$, $f_1(t,x,y)$ is $\varrho_1$-periodic with respect to $t$, i.e.
\ce
f_1(t+\varrho_1,x,y)=f_1(t,x,y), \quad t\geq 0,
\de
where $\varrho_1>0$ is a constant.
\end{enumerate}
\begin{enumerate}[$({\bf H}^p_{f_2, g_2})$]
\item For any $x, y\in\mH$, $f_2(t,x,y), g_2(t,x,y)$ are $\varrho_2$-periodic with respect to $t$, where $\varrho_2>0$ is a constant. Moreover, $\frac{\varrho_1}{\varrho_2}$ is a rational number.
\end{enumerate}

\br
The hypothesis $({\bf H}^p_{f_2, g_2})$ assures the periodicity of the evolution system $\{\mu^x_t\}_{t\in\mR}$ defined in Lemma \ref{mux}. Based on this periodicity property, we further utilize the assumptions $({\bf H}^p_{f_1})$ and $({\bf H}^a_{g_1})$ to construct the corresponding averaged equation, whose coefficients are completely independent of the small parameter $\e$.
\er

By $({\bf H}^p_{f_1})$ and $({\bf H}^p_{f_2, g_2})$, we know that there exists a constant $\varrho>0$ satisfying $\varrho=n_2\varrho_1=n_1\varrho_2$ for some $n_1, n_2\in \mN_+=\{1,2,\cdots\}$. Set
\ce
\bar f_{1,p}(x):=\frac{1}{\varrho}\int_0^\varrho \bar f_1(t,x)\dif t=\frac{1}{\varrho}\int_0^\varrho\int_\mH f_1(t,x,y)\mu^x_t(\dif y)\dif t,
\de
and we construct the corresponding averaged equation as follows:
\be\left\{\begin{array}{ll}
\dif \bar X_t=[A\bar X_t+B(\bar X_t,\bar X_t)+\bar f_{1,p}(\bar X_t)]\dif t+\bar g_1(\bar X_t)\dif W^1_t+\dif \bar K_t,\\
\bar X_t(\xi)\geq 0, \quad \xi\in[0,1],\\
\bar X_t(0)=\bar X_t(1)=0,\quad t\geq 0,\\
\bar X_0=x_0.
\end{array}
\right.
\label{averequap}
\ee
The following theorem is the third main result in this paper.

\bt\label{peristroapth}
Suppose that $({\bf H}_{f_1, g_1})$, $({\bf H}^1_{f_2, g_2})$, $({\bf H}^2_{f_2, g_2})$, $({\bf H}^p_{f_1})$, $({\bf H}^a_{g_1})$ and $({\bf H}^p_{f_2, g_2})$ hold. Then for any $T>0$ and $\eta>0$, 
\ce
\lim\limits_{\e\rightarrow 0}\mP\left(\sup\limits_{t\in[0,T]}|X_t^\e-\bar X_t|_\mH\geq\eta\right)=0,
\de
where $\bar X$ is the solution of Eq.(\ref{averequap}).
\et

The proof of the above theorem is placed in Section \ref{peristroapthproo}.

\subsection{Nonautonomous multiscale stochastic Burgers equations with reflection in the almost periodic case}

In this subsection, we state an averaging principle for nonautonomous multiscale stochastic Burgers equations with reflection in the almost periodic setting.

To derive an averaged equation with coefficients independent of $\e$, we assume more:
\begin{enumerate}[$({\bf H}^{ap}_{f_1})$]
\item $f_1$ is continuous with respect to $(t,x,y)$ and uniformly almost periodic. 
\end{enumerate}
\begin{enumerate}[$({\bf H}^{ex}_{f_2, g_2})$]
\item There exists a closed subset $\mS \subset \mH$ equipped with the norm $\|\cdot\|_{\mS}$ such that $\mV \subset \mS$ is continuous and $\mS \subset \mH$ is compact. Let $\{T_n\}_{n\in\mN}$ be a sequence of positive definite self-adjoint operators on $\mH$ such that for each $n\geq 1$,
$$
\langle y_1, y_2\rangle_n:=\left\langle y_1, T_n y_2\right\rangle_\mH, \quad y_1, y_2 \in \mH,
$$
defines a new inner product on $\mH$. And the norms $\|\cdot\|_n$ generated by $\langle \cdot, \cdot\rangle_n$ are all equivalent to $|\cdot|_\mH$ and for all $y \in \mS$ we have
$$
\|y\|_n \uparrow\|y\|_\mS \quad \text { as } n \rightarrow \infty.
$$
Moreover, $T_n: \mV \rightarrow \mV$ is continuous and there exist constants $L_7>0, L_8>0$ such that for any $x\in\mH, y\in\mV, t\in\mR$
$$
2_{\mV^*}\left\langle Ay, T_n y\right\rangle_\mV+2\<f_2(t,x,y), T_n y\>_\mH+ \|g_2(t,x,y)\|_{\cL_2(\mU_2, \mH_n)}^2 \leq-L_7\|y\|_n^2+L_8(1+|x|^2_\mH),
$$
where $\mH_n:=(\mH, \langle \cdot, \cdot\rangle_n)$ and $\|y\|_n^2:=\langle y, T_n y\rangle_\mH$.
\end{enumerate}
\begin{enumerate}[$({\bf H}^{ap}_{f_2, g_2})$]
\item $f_2, g_2$ are continuous with respect to $(t,x,y)$ and uniformly almost periodic.
\end{enumerate}

\br
The conditions $({\bf H}^{ex}_{f_2, g_2})$ and $({\bf H}^{ap}_{f_2, g_2})$ guarantee the uniform almost periodicity of the mapping $\mR\times\mH\ni (t,x)\mapsto \mu_t^x\in\sP(\mH)$ from Lemma \ref{mux}. Combined with $({\bf H}^{ap}_{f_1})$, this implies that the averaged function $\bar{f}_1(t,x)=\int_{\mathbb{H}} f_1(t,x,y)\,\mu_t^x(\mathrm{d}y)$ is likewise uniformly almost periodic. Leveraging the mean values of $\bar{f}_1(t,x)$ and $\bar{g}_1$ specified in $({\bf H}^a_{g_1})$, we formulate the corresponding averaged equation, whose coefficients are entirely independent of the small parameter $\e$.
\er

Under $({\bf H}_{f_1, g_1})$, $({\bf H}^1_{f_2, g_2})$, $({\bf H}^2_{f_2, g_2})$, $({\bf H}^{ap}_{f_1})$, $({\bf H}^{ex}_{f_2, g_2})$ and $({\bf H}^{ap}_{f_2, g_2})$, by Lemma \ref{barf1unifalmoperi} and Theorem \ref{uapmeantheo} we conclude that there exists 
$\hat f_{1,ap}: \mH\rightarrow\mH$ such that 
\ce
\lim\limits_{T\rightarrow\infty}\frac{1}{T}\int_t^{t+T}\bar f_1(s,x)\dif s=\hat f_{1,ap}(x)
\de
uniformly on $\mR_+\times\cH$, for any compact subset $\cH\subset\mH$. Thus, we construct the corresponding averaged equation as follows:
\be\left\{\begin{array}{ll}
\dif \hat X_t=[A\hat X_t+B(\hat X_t,\hat X_t)+\hat f_{1,ap}(\hat X_t)]\dif t+\bar g_{1}(\hat X_t)\dif W^1_t+\dif \hat K_t,\\
\hat X_t(\xi)\geq 0, \quad \xi\in[0,1],\\
\hat X_t(0)=\hat X_t(1)=0,\quad t\geq 0,\\
\hat X_0=x_0.
\end{array}
\right.
\label{averequaap}
\ee
The following theorem is the fourth main result in this paper.

\bt\label{aperistroapth}
Suppose that $({\bf H}_{f_1, g_1})$, $({\bf H}^1_{f_2, g_2})$, $({\bf H}^2_{f_2, g_2})$, $({\bf H}^{ap}_{f_1})$, $({\bf H}^a_{g_1})$, $({\bf H}^{ex}_{f_2, g_2})$ and $({\bf H}^{ap}_{f_2, g_2})$ hold. Then for any $T>0$ and $\eta>0$, 
\ce
\lim\limits_{\e\rightarrow 0}\mP\left(\sup\limits_{t\in[0,T]}|X_t^\e-\hat X_t|_\mH\geq\eta\right)=0,
\de
where $\hat X$ is the solution of Eq.(\ref{averequaap}).
\et

The proof of the above theorem is placed in Section \ref{aperistroapthproo}.

\section{Proof of Theorem \ref{genestroapth}}\label{genestroapthproo}

In this section, we prove Theorem \ref{genestroapth}. In Subsection \ref{someestixeye}, we give some estimates about $X^\e$ and $Y^\e$. Then an evolution system of measures is constructed in Subsection \ref{evolsyst}. In Subsection \ref{estiaverequagene}, we derive some estimates for the averaged equation (\ref{eaverequa}). Finally, we prove Theorem \ref{genestroapth}.

\subsection{Some estimates about $X^\e$ and $Y^\e$}\label{someestixeye}

In this subsection, we derive some estimates about $X^\e$ and $Y^\e$ used in the sequel.

\bl
Suppose that $({\bf H}_{f_1, g_1})$, $({\bf H}^1_{f_2, g_2})$ and $({\bf H}^2_{f_2, g_2})$ hold. Then for $T>0$, there is $C>0$ such that 
\be
\sup\limits_{\e}\mE\left[\sup\limits_{t\in[0,T]}|X_t^\e|_\mH^2+\int_0^T\|X_t^\e\|_\mV^2\dif t\right]\leq C(1+|x_0|^2_\mH+|y_0|^2_\mH),
\label{xehvesti}
\ee
and
\be
\sup\limits_{\e}\sup\limits_{t\in[0,T]}\mE|Y_t^\e|_\mH^2\leq C(1+|x_0|^2_\mH+|y_0|^2_\mH).
\label{yeesti}
\ee
\el
\begin{proof}
First of all, applying the It\^o formula to $|X_t^\e|_\mH^2$ for any $t\in[0,T]$, we have that 
\ce
|X_t^\e|_\mH^2&=&|x_0|_\mH^2+2\int_0^t {_{\mV^*}}\<AX_r^\e,X_r^\e\>_\mV \dif r+2\int_0^t{_{\mV^*}}\<B(X_r^\e,X_r^\e),X_r^\e\>_\mV\dif r\\
&& +2\int_0^t\<f_1(\frac{r}{\e},X_r^\e,Y_r^\e), X_r^\e\>_\mH\dif r+2\int_0^t\<X_r^\e, g_1(\frac{r}{\e},X_r^\e)\dif W^1_r\>_\mH\\
&& +2\int_0^t\int_0^1X_r^\e(\xi)K^\e(\dif r,\dif \xi)+\int_0^t\|g_1(\frac{r}{\e},X_r^\e)\|_{\cL_2(\mU_1,\mH)}^2\dif r. 
\de
Then the definition of $A$, Corollary \ref{Bbprop2} and Definition \ref{soludefi} $(iv)$ imply that
\ce
&&{_{\mV^*}}\<AX_r^\e,X_r^\e\>_\mV=-\|X_r^\e\|^2_\mV,\\
&&{_{\mV^*}}\<B(X_r^\e,X_r^\e),X_r^\e\>_\mV=0,\\
&&\int_0^t\int_0^1X_r^\e(\xi)K^\e(\dif r,\dif \xi)=0.
\de
 So, by (\ref{f1g1linegrow}) and the Young inequality we infer that
\ce
\left|X_t^\e\right|_\mH^2+2\int_0^t\|X_r^\e\|_\mV^2\dif r&\leq& |x_0|_\mH^2+C+C\int_0^t|X_r^\e|_\mH^2\dif r+C\int_0^t|Y_r^\e|_\mH^2\dif r\\
 &&+2\left|\int_0^t\<X_r^\e, g_1(\frac{r}{\e},X_r^\e)\dif W^1_r\>_\mH\right|. 
\de
The Burkholder-Davis-Gundy inequality and  the Young inequality imply that
\ce
&&\mE\sup\limits_{r\in[0,t]}\left|X_r^\e\right|_\mH^2+2\mE\int_0^t\|X_r^\e\|_\mV^2\dif r\\
&\leq& |x_0|_\mH^2+C+C\mE\int_0^t|X_r^\e|_\mH^2\dif r+C\mE\int_0^t|Y_r^\e|_\mH^2\dif r\no\\
&&+C \mE\left(\int_0^t|X_r^\e|_\mH^2\|g_1(\frac{r}{\e},X_r^\e)\|_{\cL_2(\mU_1,\mH)}^2\dif r\right)^{1/2}\no\\
&\leq& |x_0|_\mH^2+C+C\mE\int_0^t|X_r^\e|_\mH^2\dif r+C\mE\int_0^t|Y_r^\e|_\mH^2\dif r\no\\
&&+\frac{1}{2}\mE\sup\limits_{r\in[0,t]}|X_r^\e|_\mH^2+C\mE\int_0^t\|g_1(\frac{r}{\e},X_r^\e)\|_{\cL_2(\mU_1,\mH)}^2\dif r.
\de
and furthermore by (\ref{f1g1linegrow})
\be
\mE\sup\limits_{r\in[0,t]}\left|X_r^\e\right|_\mH^2+2\mE\int_0^t\|X_r^\e\|_\mV^2\dif r\leq C(|x_0|_\mH^2+1)+C\mE\int_0^t|X_r^\e|_\mH^2\dif r+C\mE\int_0^t|Y_r^\e|_\mH^2\dif r.
\label{xegueyegue1}
\ee

Next, we are devoted to estimating $\mE|Y_r^\e|_\mH^2$. From the It\^o formula, it follows that for any $l_1>0$
\ce
e^{l_1 t}|Y_t^\e|_\mH^2&=&|y_0|_\mH^2+l_1\int_0^te^{l_1 r}|Y_r^\e|_\mH^2\dif r+\frac{2}{\e}\int_0^te^{l_1 r}{_{\mV^*}}\<AY_r^\e,Y_r^\e\>_\mV \dif r\\
&&+\frac{2}{\e}\int_0^te^{l_1 r}\<f_2(\frac{r}{\e},X_r^\e,Y_r^\e), Y_r^\e\>_\mH\dif r+\frac{2}{\sqrt{\e}}\int_0^te^{l_1 r}\<Y_r^\e,g_2(\frac{r}{\e},X_r^\e,Y_r^\e)\dif W^2_r\>_\mH\\
&&+\frac{1}{\e}\int_0^t e^{l_1 r}\|g_2(\frac{r}{\e},X_r^\e,Y_r^\e)\|_{\cL_2(\mU_2,\mH)}^2\dif r.
\de
Besides, (\ref{discond2}) yields that
\ce
&&\frac{2}{\e}{_{\mV^*}}\<AY_r^\e,Y_r^\e\>_\mV+\frac{2}{\e}\<f_2(\frac{r}{\e},X_r^\e,Y_r^\e), Y_r^\e\>_\mH+\frac{1}{\e}\|g_2(\frac{r}{\e},X_r^\e,Y_r^\e)\|_{\cL_2(\mU_2,\mH)}^2\\
&\leq&-\frac{2\a}{\e}|Y_r^\e|_\mH^2+\frac{C}{\e}(1+|X_r^\e|_\mH^2).
\de
So, we infer that
\ce
\mE e^{l_1 t}|Y_t^\e|_\mH^2\leq |y_0|_\mH^2+(l_1-\frac{2\a}{\e})\mE\int_0^te^{l_1 r}|Y_r^\e|_\mH^2\dif r+\frac{C}{\e}\mE\int_0^te^{l_1 r}(1+|X_r^\e|_\mH^2)\dif r.
\de
By taking $l_1=\frac{2\a}{\e}$, it holds that
\be
\mE|Y_t^\e|_\mH^2\leq |y_0|_\mH^2+\frac{C}{2\a}(1+\mE\sup\limits_{r\in[0,t]}|X_r^\e|_\mH^2).
\label{yegue}
\ee

Finally, by inserting (\ref{yegue}) into (\ref{xegueyegue1}), the Gronwall inequality implies (\ref{xehvesti}). Then by (\ref{xehvesti}) and (\ref{yegue}), we obtain (\ref{yeesti}). The proof is complete.
\end{proof}

Next, in order to estimate the increments of $X^\e$, for any $R>0$ we define the stopping time
\ce
\tau_1:=\inf\{t>0: |X_t^\e|_\mH>R\},
\de
and have the following result.

\bl
Under $({\bf H}_{f_1, g_1})$, $({\bf H}^1_{f_2,g_2})$ and $({\bf H}^2_{f_2,g_2})$ it holds that for any $\d>0$ small enough,
\be
\mE\left[\int_0^{T\wedge\tau_1}|X_t^\e-X_{t(\d)}^\e|_\mH^2\dif t\right]\leq C_R\d^{1/2}\left(1+|x_0|_\mH^2+|y_0|_\mH^2\right)^{3/2},
\label{xeettd}
\ee
where $t(\d):=[\frac{t}{\d}]\d$, $[\frac{t}{\d}]$ denotes the largest integer which is less than $\frac{t}{\d}$ and $C_R>0$ is a constant independent of $\e$ and dependent on $R$.
\el
\begin{proof} 
Since by (\ref{xehvesti})
\ce
&&\mE\left[\int_0^{T\wedge\tau_1}|X_t^\e-X_{t(\d)}^\e|_\mH^2\dif t\right]=\mE\left[\int_0^{T}|X_t^\e-X_{t(\d)}^\e|_\mH^2 I_{t\leq \tau_1}\dif t\right]\\
&=& \mE\left[\int_0^\d|X_t^\e-X_{t(\d)}^\e|_\mH^2 I_{t\leq \tau_1}\dif t\right]+\mE\left[\int_\d^{T}|X_t^\e-X_{t(\d)}^\e|_\mH^2I_{t\leq \tau_1}\dif t\right]\\
&\leq&\mE\left[\int_0^\d|X_t^\e-x_0|_\mH^2I_{t\leq \tau_1}\dif t\right]+2\mE\left[\int_\d^{T}|X_t^\e-X_{t-\d}^\e|_\mH^2I_{t\leq \tau_1}\dif t\right]\\
&&+2\mE\left[\int_\d^{T}|X_{t-\d}^\e-X_{t(\d)}^\e|_\mH^2I_{t\leq \tau_1}\dif t\right]\\
&\leq&C\left(1+|x_0|_\mH^2+|y_0|_\mH^2\right)\d+2\mE\left[\int_\d^{T}|X_t^\e-X_{t-\d}^\e|_\mH^2I_{t\leq \tau_1}\dif t\right]\\
&&+2\mE\left[\int_\d^{T}|X_{t-\d}^\e-X_{t(\d)}^\e|_\mH^2I_{t\leq \tau_1}\dif t\right],
\de
and the estimates for the second term and the third term in the right hand of the above inequality are similar, we only deal with the second term $2\mE\left[\int_\d^{T}|X_t^\e-X_{t-\d}^\e|_\mH^2I_{t\leq \tau_1}\dif t\right]$.

Next, applying the It\^o formula to $|X_r^\e-X_{t-\d}^\e|_\mH^2$ for $r\in[t-\d,t]$, we have that
\ce
&&|X_t^\e-X_{t-\d}^\e|_\mH^2\\
&=&2\int_{t-\d}^t{_{\mV^*}}\<A X_r^\e, X_r^\e-X_{t-\d}^\e\>_\mV\dif r+2\int_{t-\d}^t{_{\mV^*}}\<B(X_r^\e,X_r^\e), X_r^\e-X_{t-\d}^\e\>_\mV\dif r\\
&& +2 \int_{t-\d}^t\<f_1(\frac{r}{\e},X_r^\e,Y_r^\e), X_r^\e-X_{t-\d}^\e\>_\mH\dif r +2 \int_{t-\d}^t\<X_r^\e-X_{t-\d}^\e, g_1(\frac{r}{\e},X_r^\e)\dif W^1_r\>_\mH \\
&&+2 \int_{t-\d}^t\int_0^1(X_r^\e(\xi)-X_{t-\d}^\e(\xi))K^\e(\dif r,\dif \xi)+ \int_{t-\d}^t\|g_1(\frac{r}{\e},X_r^\e)\|_{\cL_2(\mU_1,\mH)}^2\dif r.
\de
By the definition of $A$, it holds that
\ce
2\int_{t-\d}^t{_{\mV^*}}\<A X_r^\e, X_r^\e-X_{t-\d}^\e\>_\mV\dif r
&\leq&-2\int_{t-\d}^t\|X_r^\e\|_\mV^2\dif r+2\int_{t-\d}^t\|X_r^\e\|_\mV\|X_{t-\d}^\e\|_\mV\dif r\\
&\leq&2\int_{t-\d}^t\|X_r^\e\|_\mV\cdot\|X_{t-\d}^\e\|_\mV\dif r.
\de
And Corollary \ref{Bbprop2} implies that
\ce
2\int_{t-\d}^t{_{\mV^*}}\<B(X_r^\e,X_r^\e), X_r^\e-X_{t-\d}^\e\>_\mV\dif r&=&-2\int_{t-\d}^t{_{\mV^*}}\<B(X_r^\e,X_r^\e), X_{t-\d}^\e\>_\mV\dif r\\
&\leq&4\int_{t-\d}^t\|X_r^\e\|_\mV\cdot|X_r^\e|_\mH\cdot\|X_{t-\d}^\e\|_\mV\dif r.
\de
From (\ref{f1g1linegrow}), it follows that
\ce
&&2 \int_{t-\d}^t\<f_1(\frac{r}{\e},X_r^\e,Y_r^\e), X_r^\e-X_{t-\d}^\e\>_\mH\dif r\\
&\leq&2L_1\int_{t-\d}^t(1+|X_r^\e|_\mH+|Y_r^\e|_\mH)(|X_r^\e|_\mH+|X_{t-\d}^\e|_\mH)\dif r\\
&\leq&3L_1\d+7L_1\d\sup\limits_{r\in[0,T]}|X_r^\e|_\mH^2+3L_1\int_{t-\d}^t|Y_r^\e|_\mH^2\dif r,
\de
and
\ce
\int_{t-\d}^t\|g_1(\frac{r}{\e},X_r^\e)\|_{\cL_2(\mU_1,\mH)}^2\dif r\leq 2L^2_2\int_{t-\d}^t(1+|X_r^\e|_\mH^2)\dif r\leq 2L^2_2\d+2L^2_2\d\sup\limits_{r\in[0,T]}|X_r^\e|_\mH^2.
\de
Definition \ref{soludefi} yields that
\ce
2 \int_{t-\d}^t\int_0^1(X_r^\e(\xi)-X_{t-\d}^\e(\xi))K^\e(\dif r,\dif \xi)\leq 0.
\de

So, we infer that
\ce
|X_t^\e-X_{t-\d}^\e|_\mH^2&\leq& 2\int_{t-\d}^t\|X_r^\e\|_\mV(1+2|X_r^\e|_\mH)\|X_{t-\d}^\e\|_\mV\dif r\\
&&+\left(C\d+C\d\sup\limits_{r\in[0,T]}|X_r^\e|_\mH^2+3L_1\int_{t-\d}^t|Y_r^\e|_\mH^2\dif r\right)\\
&& +2\left|\int_{t-\d}^t\<X_r^\e-X_{t-\d}^\e, g_1(\frac{r}{\e},X_r^\e)\dif W^1_r\>_\mH\right| \\
&=:& J_1(t)+J_2(t)+J_3(t).
\de

For $J_1(t)$, by the H\"older inequality, the Fubini theorem and (\ref{xehvesti}), it holds that
\ce
&&\mE\left[\int_\d^{T}J_1(t)I_{t\leq \tau_1}\dif t\right]\\
&\leq&2(1+2R)\left[\mE\int_\d^{T}\int_{t-\d}^t\|X_r^\e\|_\mV^2\dif r\dif t\right]^{1/2}\left[\mE\int_\d^{T}\int_{t-\d}^t\|X_{t-\d}^\e\|_\mV^2\dif r\dif t\right]^{1/2}\\
&\leq&2(1+2R)\d \left[\mE\int_0^{T}\|X_r^\e\|_\mV^2\dif r\right]\\
&\leq&2(1+2R)\d C\left(1+|x_0|_\mH^2+|y_0|_\mH^2\right).
\de

For $J_2(t)$, by the Fubini theorem, (\ref{xehvesti}) and (\ref{yeesti}), we have that
\ce
\mE\left[\int_\d^{T}J_2(t)I_{t\leq \tau_1}\dif t\right]&\leq&CT\d+CT\d\mE\sup\limits_{r\in[0,T]}|X_r^\e|_\mH^2+3L_1\mE\int_\d^{T}\int_{t-\d}^t|Y_r^\e|_\mH^2\dif r\dif t\\
&\leq&CT\d+CT\d\mE\sup\limits_{r\in[0,T]}|X_r^\e|_\mH^2+3L_1\d\mE\int_0^{T}|Y_r^\e|_\mH^2\dif r\\
&\leq&C\d\left(1+|x_0|_\mH^2+|y_0|_\mH^2\right).
\de

For $J_3(t)$, by the Burkholder-Davis-Gundy inequality, the H\"older inequality and $(\ref{f1g1linegrow})$, it holds that
\ce
&&\mE\left[\int_\d^{T}J_3(t)I_{t\leq \tau_1}\dif t\right]\\
&\leq&2C\int_\d^{T}\mE\left(\int_{t-\d}^t|X_r^\e-X_{t-\d}^\e|_\mH^2\|g_1(\frac{r}{\e},X_r^\e)\|_{\cL_2(\mU_1,\mH)}^2\dif r\right)^{1/2}I_{t\leq \tau_1}\dif t\\
&\leq&2 CT^{1/2}\left[\mE\int_\d^{T}\int_{t-\d}^t(|X_r^\e|_\mH^2+|X_{t-\d}^\e|_\mH^2)(1+|X_r^\e|_\mH^2)I_{t\leq \tau_1}\dif r\dif t\right]^{1/2}\\
&\leq&2CT^{1/2}(1+R^2)^{1/2}\left[\mE\int_\d^{T}\int_{t-\d}^t(|X_r^\e|_\mH^2+|X_{t-\d}^\e|_\mH^2)\dif r\dif t\right]^{1/2}\\
&\leq&2CT(1+R^2)^{1/2}\d^{1/2}\left[\mE\sup\limits_{r\in[0,T]}|X_r^\e|_\mH^2\right]^{1/2}\\
&\leq&2CT(1+R^2)^{1/2}\d^{1/2}\left(1+|x_0|_\mH^2+|y_0|_\mH^2\right)^{1/2}.
\de

Finally, combining the above deduction, we obtain that
\ce
2\mE\left[\int_\d^{T}|X_{t-\d}^\e-X_{t(\d)}^\e|_\mH^2I_{t\leq \tau_1}\dif t\right]\leq C_R\d^{1/2}\left(1+|x_0|_\mH^2+|y_0|_\mH^2\right),
\de
which completes the proof.
\end{proof}

Finally we introduce the following auxiliary SPDE: 
\be\left\{\begin{array}{l}
\dif\hat{Y}_{t}^\e=\frac{1}{\e}\left[A\hat{Y}_{t}^\e+f_{2}(\frac{t}{\e},X_{t(\d)}^\e,\hat{Y}_{t}^\e)\right]\dif t+\frac{1}{\sqrt{\e}}g_{2}(\frac{t}{\e},X_{t(\d)}^\e,\hat{Y}_{t}^\e)\dif W^2_{t},\\
\hat{Y}_{t}^\e(0)=\hat{Y}_{t}^\e(1)=0,\quad t\geq 0,\\
\hat{Y}_{0}^\e=y_0.
\end{array}
\right.
\label{hatzu}
\ee
Under $({\bf H}^1_{f_2,g_2})$, the above SPDE has a unique solution $\hat{Y}_{\cdot}^\e$ in $L^2(\Omega, \sF, \mP; C([0,T],\mH))\cap L^2(\Omega, \sF, \mP; L^2([0,T],\mV))$ in the following sense (\cite{DaZ}): $\forall t \in \mathbb{R}_{+}, \psi \in C_c^2([0,1])$ with $\psi(0)=\psi(1)=0$,
\ce
\<\hat{Y}_{t}^\e,\psi\>_\mH&=&\<y_0,\psi\>_\mH+\frac{1}{\e}\int_0^t\left[{_{\mV^*}}\<A\hat{Y}_{r}^\e,\psi\>_\mV+\<f_{2}(\frac{r}{\e},X_{r(\d)}^\e,\hat{Y}_{r}^\e),\psi\>_\mH\right]\dif r\\
&&+\frac{1}{\sqrt{\e}}\int_0^t\<\psi,g_{2}(\frac{r}{\e},X_{r(\d)}^\e,\hat{Y}_{r}^\e)\dif W^2_{r}\>_\mH.
\de

\bl
Under $({\bf H}_{f_1, g_1})$, $({\bf H}^1_{f_2,g_2})$ and $({\bf H}^2_{f_2,g_2})$, it holds that 
 \be
 &&\sup\limits_{t\in[0,T]}\mE|\hat{Y}_{t}^\e|_\mH^2\leq C(1+|x_0|_\mH^{2}+|y_0|_\mH^{2}), \label{hatzub}\\
&&\mE\int_0^{T\wedge\tau_1}|Y_{t}^\e-\hat{Y}_{t}^\e|_\mH^2\dif t\leq C_R\d^{1/2}\left(1+|x_0|_\mH^2+|y_0|_\mH^2\right).
\label{unztu}
\ee
\el
\begin{proof}
Since the proof of (\ref{hatzub}) is similar to that for (\ref{yeesti}), we only prove (\ref{unztu}).

First of all, applying the It\^{o} formula to $|Y_{t}^\e-\hat{Y}_{t}^\e|_\mH^2 e^{l_2 t}$ for any $l_2>0$, one could obtain that
\ce
|Y_{t}^\e-\hat{Y}_{t}^\e|_\mH^2 e^{l_2 t}&=&l_2\int_{0}^{t}|Y_{r}^\e-\hat{Y}_{r}^\e|_\mH^2 e^{l_2 r}\dif r\\
&&+\frac{1}{\e}\int_0^t2e^{l_2 r}{_{\mV^*}}\<AY_{r}^\e-A\hat{Y}_{r}^\e, Y_{r}^\e-\hat{Y}_{r}^\e\>_\mV\dif r\\
&&+\frac{1}{\e}\int_{0}^{t}2e^{l_2 r}\<f_{2}(\frac{r}{\e},X_{r}^\e,Y_{r}^\e)-f_{2}(\frac{r}{\e},X_{r(\d)}^\e,\hat{Y}_{r}^\e), Y_{r}^\e-\hat{Y}_{r}^\e\>_\mH\dif r\\
&&+\frac{1}{\sqrt{\e}}\int_{0}^{t}2e^{l_2 r}\<Y_{r}^\e-\hat{Y}_{r}^\e, \(g_{2}(\frac{r}{\e},X_{r}^\e,Y_{r}^\e)-g_{2}(\frac{r}{\e},X_{r(\d)}^\e,\hat{Y}_{r}^\e)\)\dif W^2_{r}\>_\mH\\
&&+\frac{1}{\e}\int_{0}^{t}e^{l_2 r}\|g_{2}(\frac{r}{\e},X_{r}^\e,Y_{r}^\e)-g_{2}(\frac{r}{\e},X_{r(\d)}^\e,\hat{Y}_{r}^\e)\|^2_{\cL_2(\mU_2,\mH)}\dif r.
\de

By the definition of $A$ and the Poincar\'e inequality, we know that
\ce
2{_{\mV^*}}\<AY_{r}^\e-A\hat{Y}_{r}^\e, Y_{r}^\e-\hat{Y}_{r}^\e\>_\mV\leq -2\l_1|Y_{r}^\e-\hat{Y}_{r}^\e|_\mH^2.
\de
Then $({\bf H}^1_{f_2, g_2})$ and the Young inequality imply  that
\ce
&&2\<f_{2}(\frac{r}{\e},X_{r}^\e,Y_{r}^\e)-f_{2}(\frac{r}{\e},X_{r(\d)}^\e,\hat{Y}_{r}^\e), Y_{r}^\e-\hat{Y}_{r}^\e\>_\mH\\
&\leq&2\<f_{2}(\frac{r}{\e},X_{r}^\e,Y_{r}^\e)-f_{2}(\frac{r}{\e},X_{r}^\e,\hat{Y}_{r}^\e), Y_{r}^\e-\hat{Y}_{r}^\e\>_\mH\\
&&+2\<f_{2}(\frac{r}{\e},X_{r}^\e,\hat Y_{r}^\e)-f_{2}(\frac{r}{\e},X_{r(\d)}^\e,\hat{Y}_{r}^\e), Y_{r}^\e-\hat{Y}_{r}^\e\>_\mH\\
&\leq&3L_4|Y_{r}^\e-\hat{Y}_{r}^\e|_\mH^2 +C|X_{r}^\e-X_{r(\d)}^\e|_\mH^2,
\de
and
\ce
&&\|g_{2}(\frac{r}{\e},X_{r}^\e,Y_{r}^\e)-g_{2}(\frac{r}{\e},X_{r(\d)}^\e,\hat{Y}_{r}^\e)\|^2_{\cL_2(\mU_2,\mH)}\\
&\leq&2L_6^2|Y_{r}^\e-\hat{Y}_{r}^\e|_\mH^2 +2L_5^2|X_{r}^\e-X_{r(\d)}^\e|_\mH^2.
\de

Combining the above deduction, we obtain that
\ce
|Y_{t}^\e-\hat{Y}_{t}^\e|_\mH^2 e^{l_2 t}
&\leq&(l_2-\frac{2\a}{\e})\int_{0}^{t}e^{l_2 r}|Y_{r}^\e-\hat{Y}_{r}^\e|_\mH^2 \dif r+\frac{C}{\e}\int_{0}^{t}e^{l_2 r}|X_{r}^\e-X_{r(\d)}^\e|_\mH^2\dif r\\
&&+\frac{1}{\sqrt{\e}}\int_{0}^{t}2e^{l_2 r}\<Y_{r}^\e-\hat{Y}_{r}^\e, \(g_{2}(\frac{r}{\e},X_{r}^\e,Y_{r}^\e)-g_{2}(\frac{r}{\e},X_{r(\d)}^\e,\hat{Y}_{r}^\e)\)\dif W^2_{r}\>_\mH.
\de
which together with $l_2=\frac{2\a}{\e}$ yields that
\ce
|Y_{t}^\e-\hat{Y}_{t}^\e|_\mH^2
&\leq&\frac{1}{\sqrt{\e}}\int_{0}^{t}2e^{-\frac{2\a}{\e}(t-r)}\<Y_{r}^\e-\hat{Y}_{r}^\e,\(g_{2}(\frac{r}{\e},X_{r}^\e,Y_{r}^\e)-g_{2}(\frac{r}{\e},X_{r(\d)}^\e,\hat{Y}_{r}^\e)\)\dif W^2_{r}\>_\mH\\
&&+\frac{C}{\e}\int_{0}^{t}e^{-\frac{2\a}{\e}(t-r)}|X_{r}^\e-X_{r(\d)}^\e|_\mH^2\dif r.
\de
Integrating from $0$ to $T\wedge\tau_1$ and taking the expectation on two sides, by (\ref{xeettd}) we conclude that
\ce
\mE\int_0^{T\wedge\tau_1}|Y_{t}^\e-\hat{Y}_{t}^\e|_\mH^2 \dif t
&\leq&\frac{1}{\sqrt{\e}}\mE\int_0^{T\wedge\tau_1}\int_{0}^{t}2e^{-\frac{2\a}{\e}(t-r)}\<Y_{r}^\e-\hat{Y}_{r}^\e, \\
&&\qquad\qquad\qquad\qquad\(g_{2}(\frac{r}{\e},X_{r}^\e,Y_{r}^\e)-g_{2}(\frac{r}{\e},X_{r(\d)}^\e,\hat{Y}_{r}^\e)\)\dif W^2_{r}\>_\mH\dif t\\
&&+\frac{C}{\e}\mE\int_0^{T\wedge\tau_1}\int_{0}^{t}e^{-\frac{2\a}{\e}(t-r)}|X_{r}^\e-X_{r(\d)}^\e|_\mH^2\dif r\dif t\\
&\leq&C_R\d^{1/2}\left(1+|x_0|_\mH^2+|y_0|_\mH^2\right).
\de
The proof is complete.
\end{proof}

\subsection{An evolution system of measures}\label{evolsyst}

In this subsection, we construct an evolution system of measures which is important to define the averaged coefficient in the general case.

First of all, we take the independent version $W^{2,v}$ of the cylindrical Wiener process $W^2$ and define 
\ce
\tilde W^2_t=\left\{\begin{array}{ll}
W^2_t, ~~~\quad t\in[0,\infty),\\
W^{2,v}_{-t}, \quad t\in (-\infty,0).
\end{array}
\right.
\de
Then we fix $x\in\mH$ and consider the following SPDE: for any $s\in\mR$
\be\left\{\begin{array}{ll}
\dif Y_t^{x,y_0}=[AY_t^{x,y_0}+f_2(t,x,Y_t^{x,y_0})]\dif t+g_2(t,x,Y_t^{x,y_0})\dif \tilde W^2_t,\\
Y_t^{x,y_0}(0)=Y_t^{x,y_0}(1)=0, \quad t\geq s,\\
Y_s^{x,y_0}=y_0.
\end{array}
\right.
\label{frozequa}
\ee
Under $({\bf H}^1_{f_2,g_2})$, Eq.(\ref{frozequa}) has a unique solution denoted by $Y^{s,x,y_0}$ such that $\mE[\sup\limits_{t\in[s,T]}|Y^{s,x,y_0}_t|_\mH^2+\int_s^T\|Y^{s,x,y_0}_t\|^2_\mV\dif t]<\infty$ for any $T\geq s$ (\cite[Theorem 7.5]{DaZ}).

\bl\label{yxymomeincres}
Under $({\bf H}^1_{f_2,g_2})$ and $({\bf H}^2_{f_2,g_2})$ it holds that for any $x, x_1, x_2, y_1, y_2\in\mH$ and $t\geq s$
\be
&&\mE|Y_t^{s,x,y_0}|_\mH^2\leq C(1+|x|_\mH^2)+|y_0|_\mH^2 e^{-2\a(t-s)},\label{frozmomeesti}\\
&&\mE|Y_t^{s,x_1,y_1}-Y_t^{s,x_2,y_2}|_\mH^2\leq C|x_1-x_2|_\mH^2+|y_1-y_2|_\mH^2 e^{-2\a(t-s)}. \label{frozdiffesti}
\ee
\el

Since the proof of the above lemma is similar to that for Lemma 4.1 in \cite{Q1}, we omit it.

\bl\label{yxeseq}
Suppose that $(\mathbf{H}^1_{f_{2}, g_{2}})$ and $(\mathbf{H}^2_{f_{2},g_{2}})$ hold. Then for any $t\in\mR$ and $x\in \mH$, there exists $Y_t^x$ such that for all $t\geq s$ and $x\in \mH$,
\be
&&\mE|Y_t^{s,x,y_0}-Y_t^x|^2_\mH\leq C(1+|x|_\mH^2+|y_0|_\mH^2) e^{-2\a(t-s)},\label{ysxyyx}\\
&&\sup_{t\in\mR}\mE|Y_t^x|^2_\mH\leq C\left(1+|x|^2_\mH\right),\label{yxboun}\\
&&\sup_{t\in\mR}\mE|Y_t^{x_1}-Y_t^{x_2}|^2_\mH\leq C|x_1-x_2|_\mH^2, \label{yxdiffesti}
\ee
and $\{Y_t^x\}_{t\in\mR}$ solves
\be\left\{\begin{array}{ll}
\dif Y_t^{x}=[AY_t^{x}+f_2(t,x,Y_t^{x})]\dif t+g_2(t,x,Y_t^{x})\dif \tilde W^2_t,\\
Y_t^{x}(0)=Y_t^{x}(1)=0.
\end{array}
\right.
\label{yxeq}
\ee
\el
\begin{proof}
For any $\eta>0$, $t\geq s$ and $x\in \mH$, $Y^{s-\eta,x,y_0}$ is the solution of the following SPDE
\be\left\{\begin{array}{ll}
\dif Y_t^{s-\eta,x,y_0}=[AY_t^{s-\eta,x,y_0}+f_2(t,x,Y_t^{s-\eta,x,y_0})]\dif t+g_2(t,x,Y_t^{s-\eta,x,y_0})\dif \tilde W^2_t,\\
Y_t^{s-\eta,x,y_0}(0)=Y_t^{s-\eta,x,y_0}(1)=0, \quad t\geq s-\eta,\\
Y_{s-\eta}^{s-\eta,x,y_0}=y_0.
\end{array}
\right.
\label{ysexy}
\ee
Consequently, $Y^{s-\eta,x,y_0}-Y^{s,x,y_0}$ solves the following SPDE
\ce\left\{\begin{array}{ll}
\dif (Y_t^{s-\eta,x,y_0}-Y_t^{s,x,y_0})=\[A(Y_t^{s-\eta,x,y_0}-Y_t^{s,x,y_0})+f_2(t,x,Y_t^{s-\eta,x,y_0})-f_2(t,x,Y_t^{s,x,y_0})\]\dif t\\
\qquad\qquad\qquad\qquad\qquad+\(g_2(t,x,Y_t^{s-\eta,x,y_0})-g_2(t,x,Y_t^{s,x,y_0})\)\dif \tilde W^2_t,\\
Y_s^{s-\eta,x,y_0}-Y_s^{s,x,y_0}=Y_s^{s-\eta,x,y_0}-y_0.
\end{array}
\right.
\de
Applying the It\^{o} formula to $e^{l_3(t-s)}|Y_t^{s-\eta,x,y_0}-Y_t^{s,x,y_0}|^2_\mH$ for any $l_3>0$, we obtain that
\ce
&&e^{l_3(t-s)}|Y_t^{s-\eta,x,y_0}-Y_t^{s,x,y_0}|^2_\mH\\
&=&|Y_s^{s-\eta,x,y_0}-y_0|^2_\mH+l_3\int_s^te^{l_3(r-s)}|Y_r^{s-\eta,x,y_0}-Y_r^{s,x,y_0}|^2_\mH\dif r\\
&&+2\int_s^te^{l_3(r-s)}{_{\mV^*}}\<AY_r^{s-\eta,x,y_0}-AY_r^{s,x,y_0},Y_r^{s-\eta,x,y_0}-Y_r^{s,x,y_0}\>_\mV\dif r\\
&&+2\int_s^te^{l_3(r-s)}\<f_2(r,x,Y_r^{s-\eta,x,y_0})-f_2(r,x,Y_r^{s,x,y_0}),Y_r^{s-\eta,x,y_0}-Y_r^{s,x,y_0}\>_\mH\dif r\\
&&+2\int_s^te^{l_3(r-s)}\<Y_r^{s-\eta,x,y_0}-Y_r^{s,x,y_0}, (g_2(r,x,Y_r^{s-\eta,x,y_0})-g_2(r,x,Y_r^{s,x,y_0}))\dif \tilde W^2_r\>_\mH\\
&&+\int_s^te^{l_3(r-s)}\|g_2(r,x,Y_r^{s-\eta,x,y_0})-g_2(r,x,Y_r^{s,x,y_0})\|^2_{\cL_2(\mU_2,\mH)}\dif r.
\de
(\ref{discond1}) implies that
\ce
&&2\<AY_r^{s-\eta,x,y_0}-AY_r^{s,x,y_0},Y_r^{s-\eta,x,y_0}-Y_r^{s,x,y_0}\>_\mV\\
&&+2\<f_2(r,x,Y_r^{s-\eta,x,y_0})-f_2(r,x,Y_r^{s,x,y_0}),Y_r^{s-\eta,x,y_0}-Y_r^{s,x,y_0}\>_\mH\\
&&+\|g_2(r,x,Y_r^{s-\eta,x,y_0})-g_2(r,x,Y_r^{s,x,y_0})\|^2_{\cL_2(\mU_2,\mH)}\\
&\leq&-2\a|Y_r^{s-\eta,x,y_0}-Y_r^{s,x,y_0}|^2_\mH.
\de
So, taking $l_3=2\a$, by (\ref{frozmomeesti}) we conclude that for any $\eta>0$, $t\geq s$ and $x\in \mH$,
\ce
\mE|Y_t^{s-\eta,x,y_0}-Y_t^{s,x,y_0}|^2_\mH&\leq& \mE|Y_s^{s-\eta,x,y_0}-y_0|^2_\mH e^{-2\a(t-s)}\leq C(1+|x|_\mH^2+|y_0|_\mH^2) e^{-2\a(t-s)}.
\de
As $s\rightarrow-\infty$, it follows that
\ce
\lim_{s\rightarrow-\infty}\mE|Y_t^{s-\eta,x,y_0}-Y_t^{s,x,y_0}|^2_\mH=0.
\de
Then for any $t\in\mR$ and $x\in \mH$, $\{Y_t^{\cdot,x,y_0}\}_{t\geq \cdot}$ is a Cauchy sequence in $L^2(\Omega,\sF_t,\mP;\mH)$. Thus there exists a random variable $Y_t^{x,y_0}\in L^2(\Omega,\sF_t,\mP;\mH)$ such that for all $t\geq s$ and $x\in \mH$, 
\be
\mE|Y_t^{x,y_0}-Y_t^{s,x,y_0}|^2_\mH\leq C(1+|x|_\mH^2+|y_0|_\mH^2) e^{-2\a(t-s)}.
\label{yvsyl2}
\ee

Next, we show that $Y_t^{x,y_0}$ is independent of $y_0$. Indeed, (\ref{frozdiffesti}) implies that for any $y_1, y_2 \in \mH$,
\ce
\lim_{s\rightarrow-\infty}\mE|Y_t^{s,x,y_1}-Y_t^{s,x,y_2}|^2_\mH=0.
\de
Note that
\ce
\mE|Y_t^{x,y_1}-Y_t^{x,y_2}|^2_\mH&\leq& 3\mE|Y_t^{x,y_1}-Y_t^{s,x,y_1}|^2_\mH+3\mE|Y_t^{s,x,y_1}-Y_t^{s,x,y_2}|^2_\mH\\
&&+3\mE|Y_t^{s,x,y_2}-Y_t^{x,y_2}|^2_\mH.
\de
As $s\rightarrow-\infty$, by (\ref{frozdiffesti}) and (\ref{yvsyl2}) we get $\mE|Y_t^{x,y_1}-Y_t^{x,y_2}|^2_\mH=0$. So, $Y_t^{x,y_1}=Y_t^{x,y_2}$ in $L^2(\Omega,\sF_t,\mP;\mH)$. We henceforth write $Y_t^{x,y_0}$ by $Y_t^x$ and then (\ref{yvsyl2}) is just (\ref{ysxyyx}).

In order to prove (\ref{yxboun}), from (\ref{yvsyl2}) and (\ref{frozmomeesti}) it follows
\ce
\mE|Y_t^x|^2_\mH&\leq& 2\mE|Y_t^x-Y_t^{s,x,y_0}|^2_\mH+2\mE|Y_t^{s,x,y_0}|^2_\mH\\
&\leq&C(1+|x|_\mH^2+|y_0|_\mH^2) e^{-2\a(t-s)}+C(1+|x|_\mH^2)+2|y_0|_\mH^2 e^{-2\a(t-s)}.
\de
Letting $s\rightarrow-\infty$, we obtain (\ref{yxboun}). Then the similar deduction implies (\ref{yxdiffesti}).

Finally, letting $\eta\rightarrow+\infty$ in (\ref{ysexy}), by (\ref{yvsyl2}) we conclude that for any $t\in\mR$ and $x\in \mH$,
\ce\left\{\begin{array}{ll}
\dif Y_t^{x}=[AY_t^{x}+f_2(t,x,Y_t^{x})]\dif t+g_2(t,x,Y_t^{x})\dif \tilde W^2_t,\\
Y_t^{x}(0)=Y_t^{x}(1)=0,
\end{array}
\right.
\de
which is exactly (\ref{yxeq}). The proof is complete.
\end{proof}

Next, for each $t\in\mR$ and $x\in\mH$, we denote the distribution of $Y_t^x$ by $\mu_t^x$. Besides, let $\{P_{s,t}^x\}_{t\geq s}$ denote the semigroup associated with the process $\{Y_t^{s,x,y_0}\}_{t\geq s}$, i.e., for any bounded 
measurable function $h: \mH\rightarrow\mR$,
\ce
P_{s,t}^x h(y_0)=\mE h(Y_t^{s,x,y_0}).
\de
Then we say that a family $\{\mu_t^x\}_{t\in\mR}$ defines an evolution system of measures for the semigroup $\{P_{s,t}^x\}_{t\geq s}$, if for any $t\geq s$, $x\in\mH$ and $h\in C_b(\mH)$, where $C_b(\mH)$ denotes all continuous and bounded functions on $\mH$,
\be
\int_{\mH}P_{s,t}^x h(y)\mu_s^x(\dif y)=\int_{\mH}h(y)\mu_t^x(\dif y).
\label{evo}
\ee

\bl\label{mux}
Suppose that $(\mathbf{H}^1_{f_{2}, g_{2}})$ and $(\mathbf{H}^2_{f_{2}, g_{2}})$ hold. Then $\{\mu_t^x\}_{t\in\mR}$ defines an evolution system of measures for the semigroup $\{P_{s,t}^x\}_{t\geq s}$ and satisfies 
\be
\sup_{t\in\mR}\int_{\mH}|y|_\mH^2\mu_t^x(\dif y)\leq C\left(1+|x|^2_\mH\right).
\label{muxboun}
\ee
Moreover, there exists a constant $C>0$ such that for any $t\geq s$, $x\in\mH$, and Lipschitz continuous function $h$ on $\mH$
\be
\left|P_{s,t}^x h(y_0)-\int_\mH h(y)\mu_t^x(\dif y)\right|\leq C\|h\|_{Lip}(1+|x|_\mH+|y_0|_\mH) e^{-\a(t-s)}.
\label{evoesti}
\ee
Besides, evolution systems $\{\mu_t^x\}_{t\in\mR}$ of measures for the semigroup $\{P_{s,t}^x\}_{t\geq s}$ with (\ref{muxboun}) are unique. 
\el

Since the proof of the above lemma is similar to that for Proposition 5.3 in \cite{CL}, we omit it.

\subsection{Some estimates for the averaged equation (\ref{eaverequa})}\label{estiaverequagene}

In this subsection, we present some estimates for the averaged equation (\ref{eaverequa}) to obtain the averaging principle in the general setting.

\bl
Assume that $({\bf H}_{f_1, g_1})$, $({\bf H}^1_{f_2,g_2})$ and $({\bf H}^2_{f_2, g_2})$ hold. Then Eq.(\ref{eaverequa}) has a unique solution $(\bar{X}^\e, \bar{K}^\e)$ such that 
\be
\sup\limits_\e\mE\left[\sup _{t \in[0, T]}\left|\bar X_t^\e\right|_\mH^2+\int_0^T\left\|\bar X_t^\e\right\|_\mV^2\dif t\right]\leq C\left(1+|x_0|_\mH^2\right).
\label{barxee}
\ee
Moreover, set the stopping time $\tau_2:=\inf\{t>0: |\bar X_t^\e|_\mH>R\}$, and for any $\d>0$ small enough,
\be
\mE\int_0^{T\wedge\tau_2}|\bar X_t^\e-\bar X_{t(\d)}^\e|_\mH^2\dif t\leq C\d^{1/2}\left(1+|x_0|_\mH^2\right).
\label{barxeeincr}
\ee
\el
\begin{proof} 
First, we notice that 
\ce
\bar{f}_1(t,x)=\int_\mH f_1(t,x,y)\mu^x_t(\dif y), \quad t\geq 0.
\de
By $({\bf H}_{f_1, g_1})$, (\ref{frozdiffesti}) and (\ref{evoesti}), the same deduction to that in \cite[Lemma 4.3]{Q1} yields that for $t\in\mR_+, x_1,x_2,x\in\mH$
\be
&&|\bar{f}_1(t,x_1)-\bar{f}_1(t,x_2)|_\mH\leq C|x_1-x_2|_\mH,\label{barf1txlip}\\
&&|\bar{f}_1(t,x)|_\mH\leq C(1+|x|_\mH).\label{barf1txlinegrow}
\ee
So, by Theorem 3.2 in \cite{Q2}, we obtain that Eq.(\ref{eaverequa}) has a unique solution $(\bar{X}^\e, \bar{K}^\e)$. Then since the proofs of (\ref{barxee}) and (\ref{barxeeincr}) are similar to that for (\ref{xehvesti}) and (\ref{xeettd}), respectively, we omit them.
\end{proof}

Next, in order to estimate the difference between $X^\e$ and $\bar{X}^\e$, we define the stopping time 
\ce
\tau_3:=\inf\left\{t>0, |X^\e_t|_\mH+|\bar{X}^\e_t|_\mH+\int_0^t\(\|X^\e_{r}\|_\mV^2+\|\bar{X}^\e_{r}\|^2_\mV\)\dif r>R\right\}.
\de
So, $\tau_3\leq \tau_1\wedge\tau_2$.

\bp\label{xebarxees}
Assume that $({\bf H}_{f_1, g_1})$, $({\bf H}^1_{f_2,g_2})$ and $({\bf H}^2_{f_2, g_2})$ hold. Then for any $T>0$, there exists a constant $C_R>0$ such that
\be
\mE\sup\limits_{t\in[0,T\wedge \tau_3]}|X^\e_t-\bar{X}^\e_t|_\mH^{2}\leq C_{R}\left(\d^{1/2}+\d^{1/4}+(\frac{\e}{\d})^{1/2}\right)\left(1+|x_0|_\mH^2+|y_0|_\mH^2\right).
\label{xebarxeesti}
\ee
\ep
\begin{proof}
The proof is divided into two steps. In the first step, we prove (\ref{xebarxeesti}). Then a key inequality in the first step is shown in the second step.

{\bf Step 1.} We prove (\ref{xebarxeesti}).

Set $Z^\e_t:=X^\e_{t}-\bar{X}^\e_{t}$ and by It\^o's formula we get that
\ce
|Z^\e_t|_\mH^2&=&2\int_{0}^t{_{\mV^*}}\<AX^\e_{r}-A\bar{X}^\e_{r},Z^\e_r\>_\mV\dif r+2\int_{0}^t{_{\mV^*}}\<B(X^\e_{r},X^\e_{r})-B(\bar{X}^\e_{r},\bar{X}^\e_{r}),Z^\e_r\>_\mV\dif r\no\\
&&+2\int_{0}^t\< Z^\e_r, f_1(\frac{r}{\e},X^\e_{r},Y^\e_{r})-\bar{f}_1(\frac{r}{\e},\bar{X}^\e_{r})\>_\mH\dif r   \no\\
&&+2 \int_{0}^t\<Z^\e_r, (g_1(\frac{r}{\e},X^\e_{r})-g_1(\frac{r}{\e},\bar X_r^\e))\dif W^1_r\>_\mH \\
&&+2\int_{0}^t\int_0^1Z^\e_r(\xi)K^\e(\dif r,\dif \xi)-2\int_{0}^{t}\int_0^1Z^\e_r(\xi)\bar{K}^\e(\dif r,\dif\xi)\no\\
&&+\int_{0}^{t}\|g_1(\frac{r}{\e},X^\e_{r})-g_1(\frac{r}{\e},\bar X_r^\e)\|_{\cL_2(\mU_1,\mH)}^{2} \dif r.
\de

By the definition of $A$ and Lemma \ref{Bbprop3}, it holds that
\ce
2{_{\mV^*}}\<AX^\e_{r}-A\bar{X}^\e_{r},Z^\e_r\>_\mV=-2\|Z^\e_r\|_\mV^2.
\de
and
\ce
2{_{\mV^*}}\<B(X^\e_{r},X^\e_{r})-B(\bar{X}^\e_{r},\bar{X}^\e_{r}),Z^\e_r\>_\mV&\leq& 4|Z^\e_r|_\mH(\|X^\e_{r}\|_\mV+\|\bar{X}^\e_{r}\|_\mV)\|Z^\e_r\|_\mV\\
&\leq&8|Z^\e_r|_\mH^2(\|X^\e_{r}\|^2_\mV+\|\bar{X}^\e_{r}\|^2_\mV)+\|Z^\e_r\|_\mV^2.
\de
Next, note that
\ce
&&2\< Z^\e_r, f_1(\frac{r}{\e},X^\e_{r},Y^\e_{r})-\bar{f}_1(\frac{r}{\e},\bar{X}^\e_{r})\>_\mH\\
&=&\bigg[2\<Z^\e_r, f_1(\frac{r}{\e},X^\e_{r},Y^\e_{r})-f_1(\frac{r}{\e},X^\e_{r(\d)},\hat{Y}^\e_{r}) \>_\mH+2\<Z^\e_r,-\bar{f}_1(\frac{r}{\e},X^\e_{r})+\bar{f}_1(\frac{r}{\e},X^\e_{r(\d)})\>_\mH\\
&&+2\<Z^\e_r,\bar{f}_1(\frac{r}{\e},X^\e_{r})-\bar{f}_1(\frac{r}{\e},\bar{X}^\e_{r}) \>_\mH\bigg]\\
&&+2\<Z^\e_r-Z^\e_{r(\d)},f_1(\frac{r}{\e},X^\e_{r(\d)},\hat{Y}^\e_{r})-\bar{f}_1(\frac{r}{\e},X^\e_{r(\d)})\>_\mH\\
&&+2\<Z^\e_{r(\d)},f_1(\frac{r}{\e},X^\e_{r(\d)},\hat{Y}^\e_{r})-\bar{f}_1(\frac{r}{\e},X^\e_{r(\d)})\>_\mH\\
&=:&\cK_1(r)+\cK_2(r)+\cK_3(r).
\de
So, by the Lipschitz continuity of $f_1, \bar{f}_1$ and the Young inequality, we get that
\ce
\int_{0}^t\cK_1(r)\dif r&\leq& C\int_{0}^t|Z^\e_r|_\mH^2\dif r+C\int_{0}^t|X^\e_{r}-X^\e_{r(\d)}|_\mH^2\dif r\\
&&+C\int_{0}^t|Y^\e_{r}-\hat{Y}^\e_{r}|_\mH^2\dif r,
\de
and
\ce
\int_{0}^t\cK_2(r)\dif r\leq C\int_{0}^t(|X^\e_r-X^\e_{r(\d)}|_\mH+|\bar{X}^\e_r-\bar{X}^\e_{r(\d)}|_\mH)(1+|X^\e_{r(\d)}|_\mH+|\hat{Y}^\e_{r}|_\mH)\dif r.
\de
And Definition \ref{soludefi} yields that $\int_{0}^t\int_0^1X^\e_{r}(\xi)K^\e(\dif r,\dif \xi)=0$ and $\bar{X}^\e_{r}(\xi)\geq 0$ for $r\in [0,t]$ and $\xi\in[0,1]$. So,
\ce
2\int_{0}^t\int_0^1Z^\e_r(\xi)K^\e(\dif r,\dif \xi)=2\int_{0}^t\int_0^1X^\e_{r}(\xi)K^\e(\dif r,\dif \xi)-2\int_{0}^t\int_0^1\bar X^\e_{r}(\xi)K^\e(\dif r,\dif \xi)\leq 0.
\de
By the same deduction to the above inequality, we know that
\ce
-2\int_{0}^t\int_0^1Z^\e_r(\xi)\bar{K}^\e(\dif r,\dif\xi)\leq 0.
\de
By $({\bf H}_{f_1, g_1})$, it holds that
\ce
\int_{0}^{t}\|g_1(\frac{r}{\e},X^\e_{r})-g_1(\frac{r}{\e},\bar X_t^\e)\|_{\cL_2(\mU_1,\mH)}^{2} \dif r\leq L_2^2\int_{0}^{t}|Z^\e_r|_\mH^2\dif r.
\de

Combining the above deduction, we obtain that
\ce
&&|Z^\e_t|^{2}_\mH+\int_{0}^t\|Z^\e_r\|_\mV^2\dif r\\
&\leq& C\int_{0}^t|Z^\e_r|_\mH^2(1+\|X^\e_{r}\|^2_\mV+\|\bar{X}^\e_{r}\|^2_\mV)\dif r\\
&&+C\int_{0}^t|X^\e_{r}-X^\e_{r(\d)}|_\mH^2\dif r+C\int_{0}^t|Y^\e_{r}-\hat{Y}^\e_{r}|_\mH^2\dif r\\
&&+C\int_{0}^t(|X^\e_r-X^\e_{r(\d)}|_\mH+|\bar{X}^\e_r-\bar{X}^\e_{r(\d)}|_\mH)(1+|X^\e_{r(\d)}|_\mH+|\hat{Y}^\e_{r}|_\mH)\dif r\\
&&+2\left|\int_{0}^t\<Z^\e_{r(\d)},f_1(\frac{r}{\e},X^\e_{r(\d)},\hat{Y}^\e_{r})-\bar{f}_1(\frac{r}{\e},X^\e_{r(\d)})\>_\mH\dif r\right|\\
&&+2\left|\int_{0}^t\<Z^\e_r,  (g_1(\frac{r}{\e},X^\e_{r})-g_1(\frac{r}{\e},\bar X_r^\e))\dif W^1_r\>_\mH\right|,
\de
which together with the Gronwall inequality yields that
\be
&&\mE\sup\limits_{t\in[0,T\wedge \tau_3]}|Z^\e_t|_\mH^2+\mE\int_{0}^{T\wedge \tau_3}\|Z^\e_r\|_\mV^2\dif r\no\\
&\leq& C_R\mE\int_{0}^{T\wedge \tau_3}|X^\e_{r}-X^\e_{r(\d)}|_\mH^2\dif r+C_R\mE\int_{0}^{T\wedge \tau_3}|Y^\e_{r}-\hat{Y}^\e_{r}|_\mH^2\dif r\no\\
&&+C_R\mE\int_0^{T\wedge \tau_3}(|X^\e_r-X^\e_{r(\d)}|_\mH+|\bar{X}^\e_r-\bar{X}^\e_{r(\d)}|_\mH)(1+|X^\e_{r(\d)}|_\mH+|\hat{Y}^\e_{r}|_\mH)\dif r\no\\
&&+2C_R\mE\sup\limits_{t\in[0,T\wedge \tau_3]}\left|\int_{0}^t\<Z^\e_{r(\d)},f_1(\frac{r}{\e},X^\e_{r(\d)},\hat{Y}^\e_{r})-\bar{f}_1(\frac{r}{\e},X^\e_{r(\d)})\>_\mH\dif r\right|\no\\
&&+2C_R\mE\sup\limits_{t\in[0,T\wedge \tau_3]}\left|\int_{0}^t\<Z^\e_r,  (g_1(\frac{r}{\e},X^\e_{r})-g_1(\frac{r}{\e},\bar X_r^\e))\dif W^1_r\>_\mH\right|\no\\
&=:&I_1+I_2+I_3+I_4+I_5.
\label{i12345}
\ee

In the following, (\ref{xeettd}) and (\ref{unztu}) imply that
\be
I_1+I_2\leq C_R\d^{1/2}\left(1+|x_0|_\mH^2+|y_0|_\mH^2\right).
\label{i12}
\ee
Then by (\ref{xeettd}), (\ref{barxeeincr}), (\ref{xehvesti}), (\ref{hatzub}) and the H\"older inequality, it holds that
\be
I_3&\leq& C_R\left(\mE\int_0^{T\wedge \tau_3}(|X^\e_r-X^\e_{r(\d)}|_\mH^2+|\bar{X}^\e_r-\bar{X}^\e_{r(\d)}|_\mH^2)\dif r\right)^{1/2}\no\\
&&\times\left(\mE\int_0^{T\wedge \tau_3}(1+|X^\e_{r(\d)}|_\mH^2+|\hat{Y}^\e_{r}|_\mH^2)\dif r\right)^{1/2}\no\\
&\leq&C_{R}\d^{1/4}\left(1+|x_0|_\mH^2+|y_0|_\mH^2\right).
\label{i3}
\ee

By the deduction in {\bf Step 2}, we infer that
\be
I_4\leq C_R\((\frac{\e}{\d})^{1/2}+\d^{1/2}\)\left(1+|x_0|_\mH^2+|y_0|_\mH^2\right).
\label{i4}
\ee

For $I_5$, from the Burkholder-Davis-Gundy inequality, the Young inequality and $({\bf H}_{f_1, g_1})$, it follows that
\be
I_5&\leq& 2C\mE\left(\int_{0}^{T\wedge \tau_3}|Z^\e_r|_\mH^2\|g_1(\frac{r}{\e},X^\e_{r})-g_1(\frac{r}{\e},\bar X_r^\e)\|_{\cL_2(\mU_1,\mH)}^2 \dif r\right)^{1/2}\no\\
&\leq&\frac{1}{2}\mE\sup\limits_{t\in[0,T\wedge \tau_3]}|Z^\e_t|_\mH^{2}+C\mE\int_{0}^{T\wedge \tau_3}\|g_1(\frac{r}{\e},X^\e_{r})-g_1(\frac{r}{\e},\bar X_r^\e)\|_{\cL_2(\mU_1,\mH)}^2 \dif r\no\\
&\leq&\frac{1}{2}\mE\sup\limits_{t\in[0,T\wedge \tau_3]}|Z^\e_t|_\mH^{2}+C\int_{0}^{T}\mE\sup\limits_{t\in[0,r\wedge \tau_3]}|Z^\e_t|_\mH^{2}\dif r.
\label{i5}
\ee

Combining (\ref{i12})-(\ref{i5}) with (\ref{i12345}), we can get
\ce
&&\mE\sup\limits_{t\in[0,T\wedge \tau_3]}|Z^\e_t|_\mH^{2}+\mE\int_{0}^{T\wedge \tau_3}\|Z^\e_r\|_\mV^2\dif r\no\\
&\leq&C_{R}\left(\d^{1/2}+\d^{1/4}+(\frac{\e}{\d})^{1/2}\right)\left(1+|x_0|_\mH^2+|y_0|_\mH^2\right)+C\int_{0}^{T}\mE\sup\limits_{t\in[0,r\wedge \tau_3]}|Z^\e_t|_\mH^{2}\dif r,
\de
which together with the Gronwall inequality implies (\ref{xebarxeesti}).

{\bf Step 2.} We prove (\ref{i4}).

For $I_4$, it holds that
\ce
I_4&\leq&2C_R\mE\sup\limits_{t\in[0,T\wedge \tau_3]}\left|\int_{0}^{[\frac{t}{\d}]\d}\<Z^\e_{r(\d)},f_1(\frac{r}{\e},X^\e_{r(\d)},\hat{Y}^\e_{r})-\bar{f}_1(\frac{r}{\e},X^\e_{r(\d)})\>_\mH\dif r\right|\no\\
&&+2C_R\mE\sup\limits_{t\in[0,T\wedge \tau_3]}\left|\int_{[\frac{t}{\d}]\d}^t\<Z^\e_{r(\d)},f_1(\frac{r}{\e},X^\e_{r(\d)},\hat{Y}^\e_{r})-\bar{f}_1(\frac{r}{\e},X^\e_{r(\d)})\>_\mH\dif r\right|\no\\
 &=:&I_{41}+I_{42}.
\de

Next, we are devoted to estimating $I_{41}$. Note that
\be
I_{41}
&=&2C_R\mE\sup\limits_{t\in[0,T\wedge \tau_3]}
\Big|\sum\limits_{k=0}^{[\frac{t}{\d}]-1}\int_{k\d}^{(k+1)\d}\<Z^\e_{k\d},f_1(\frac{r}{\e},X^\e_{k\d},\hat{Y}^\e_{r})-\bar{f}_1(\frac{r}{\e},X^\e_{k\d})\>_\mH\dif r\Big|\no\\
&\leq&2C_R\mE\sup\limits_{t\in[0,T\wedge \tau_3]}\sum\limits_{k=0}^{[\frac{t}{\d}]-1}\left|\int_{k\d}^{(k+1)\d}\<Z^\e_{k\d},f_1(\frac{r}{\e},X^\e_{k\d},\hat{Y}^\e_{r})-\bar{f}_1(\frac{r}{\e},X^\e_{k\d})\>_\mH\dif r\right|\no\\
&\leq&2C_R\sum\limits_{k=0}^{[\frac{T}{\d}]-1}\mE\Big|\int_{k\d}^{(k+1)\d}\<Z^\e_{k\d},f_1(\frac{r}{\e},X^\e_{k\d},\hat{Y}^\e_{r})-\bar{f}_1(\frac{r}{\e},X^\e_{k\d})\>_\mH\dif r\Big|\no\\
&\leq&2C_R[\frac{T}{\d}]\sup_{0\leq k\leq [\frac{T}{\d}]-1}\mE\Big|\int_{k\d}^{(k+1)\d}\<Z^\e_{k\d},f_1(\frac{r}{\e},X^\e_{k\d},\hat{Y}^\e_{r})-\bar{f}_1(\frac{r}{\e},X^\e_{k\d})\>_\mH\dif r\Big|\no\\
&\leq&2C_R(\frac{T}{\d})\sup_{0\leq k\leq [\frac{T}{\d}]-1}\mE\Big|\<Z^\e_{k\d},\int_{k\d}^{(k+1)\d}(f_1(\frac{r}{\e},X^\e_{k\d},\hat{Y}^\e_{r})-\bar{f}_1(\frac{r}{\e},X^\e_{k\d}))\dif r\>_\mH\Big|\no\\
&\leq&2C_R(\frac{T}{\d})\sup_{0\leq k\leq [\frac{T}{\d}]-1}\mE|Z^\e_{k\d}|_\mH\left|\int_{k\d}^{(k+1)\d}(f_1(\frac{r}{\e},X^\e_{k\d},\hat{Y}^\e_{r})-\bar{f}_1(\frac{r}{\e},X^\e_{k\d}))\dif r\right|_\mH\no\\
 &\leq&2C_R(\frac{T}{\d})\sup_{0\leq k\leq [\frac{T}{\d}]-1}(\mE|Z^\e_{k\d}|_\mH^2)^{1/2}\no\\
 &&\times\left(\mE\left|\int_{k\d}^{(k+1)\d}(f_1(\frac{r}{\e},X^\e_{k\d},\hat{Y}^\e_{r})-\bar{f}_1(\frac{r}{\e},X^\e_{k\d}))\dif r\right|_\mH^2\right)^{1/2}\no\\
  &\leq&2C_R(\frac{T}{\d})\left(1+|x_0|_\mH^2+|y_0|_\mH^2\right)^{1/2}\sup_{0\leq k\leq [\frac{T}{\d}]-1}\left(2\int_{k\d}^{(k+1)\d}\int_{r}^{(k+1)\d}\Psi(u,r)\dif u\dif r\right)^{1/2},\no\\
\label{b41c}
\ee
where for $k\d\leq r\leq u\leq(k+1)\d$
$$
\Psi(u,r):=\mE\<f_1(\frac{u}{\e},X^\e_{k\d},\hat{Y}^\e_{u})-\bar{f}_1(\frac{u}{\e},X^\e_{k\d}),f_1(\frac{r}{\e},X^\e_{k\d},\hat{Y}^\e_{r})-\bar{f}_1(\frac{r}{\e},X^\e_{k\d})\>_\mH.
$$

In the following, for any random variables $X, Y\in L^2(\Omega,\mathscr{F}_{k\d},\mP; \mH)$, we construct the following SPDE
\ce\left\{\begin{array}{l}
\dif \check{Y}_t^{\e, X,Y}=\frac{1}{\e}[A\check{Y}_t^{\e,X,Y}+f_2(\frac{t}{\e},X, \check{Y}_t^{\e,X,Y})]\dif t+\frac{1}{\sqrt{\e}}g_2(\frac{t}{\e},X, \check{Y}_t^{\e,X,Y})\dif W^2_t,\\
\check{Y}_{t}^{\e,X,Y}(0)=\check{Y}_{t}^{\e,X,Y}(1)=0, \quad t\geq k\d,\\
\check{Y}_{k\d}^{\e,X,Y}=Y.
\end{array}
\right.
\de
Under $({\bf H}^1_{f_2, g_2})$, the above SPDE has a unique solution $\check{Y}^{\e,X,Y}$. Moreover, it holds that for $t \in[k \delta,(k+1) \delta]$
\ce
\hat{Y}_t^\e=\check{Y}_t^{\e,X_{k \delta}^\e, \hat{Y}_{k \delta}^\e},
\de
and
\ce
\Psi(u,r)=\mE\<f_1(\frac{u}{\e},X^\e_{k\d},\check{Y}^{\e,X_{k \delta}^\e, \hat{Y}_{k \delta}^\e}_{u})-\bar{f}_1(\frac{u}{\e},X^\e_{k\d}),f_1(\frac{r}{\e},X^\e_{k\d},\check{Y}^{\e,X_{k \delta}^\e, \hat{Y}_{k \delta}^\e}_{r})-\bar{f}_1(\frac{r}{\e},X^\e_{k\d})\>_\mH.
\de
Since $X_{k \delta}^\e, \hat{Y}_{k \delta}^\e$ are $\sF_{k\d}$-measurable, and for any $x,y\in\mH$, $\check{Y}_t^{\e,x,y}$ is independent of $\sF_{k\d}$, we obtain that
\ce
&&\Psi(u,r)\\
 &=&\mE\Bigg[\mE\Bigg[\<f_1(\frac{u}{\e},X^\e_{k\d},\check{Y}^{\e,X_{k \delta}^\e, \hat{Y}_{k \delta}^\e}_{u})-\bar{f}_1(\frac{u}{\e},X^\e_{k\d}),f_1(\frac{r}{\e},X^\e_{k\d},\check{Y}^{\e,X_{k \delta}^\e, \hat{Y}_{k \delta}^\e}_{r})-\bar{f}_1(\frac{r}{\e},X^\e_{k\d})\>_\mH\Bigg{|}\sF_{k\d}\Bigg]\Bigg]\\
&=&\mE\Bigg[\mE\Bigg[\<f_{1}(\frac{u}{\e},x,\check{Y}^{\e,x, y}_{u})-\bar{f}_{1}(\frac{u}{\e},x), f_{1}(\frac{r}{\e},x,\check{Y}^{\e,x, y}_{r})-\bar{f}_{1}(\frac{r}{\e},x)\>_\mH\Bigg]\Bigg{|}_{(x,y)=(X_{k\d}^\e,\hat{Y}_{k \delta}^\e)}\Bigg].
\de

Now, we investigate $\check{Y}^{\e,x, y}_{u}$. On the one hand, it holds that
\ce
\check{Y}^{\e, x, y}_{u}&=&y+\frac{1}{\e} \int_{k\d}^{u}[A\check{Y}^{\e,x, y}_{\iota}+f_2(\frac{\iota}{\e}, x,\check{Y}^{\e,x, y}_{\iota})]\dif \iota+\frac{1}{\sqrt{\e}} \int_{k\d}^{u} g_2(\frac{\iota}{\e}, x, \check{Y}^{\e,x, y}_{\iota})\dif W^2_\iota\\
&=&y+\int_{\frac{k\d}{\e}}^{\frac{u}{\e}}\left[A\check{Y}^{\e,x, y}_{\e v}+f_2(v,x,\check{Y}^{\e,x, y}_{\e v})\right]\dif v+\int_{\frac{k\d}{\e}}^{\frac{u}{\e}} g_2(v, x, \check{Y}^{\e,x, y}_{\e v})\dif \check{W}^2_v,
\de
where $\check{W}^2_\cdot:=\frac{1}{\sqrt{\e}}W^2_{\e \cdot}$ is a $\mU_2$-valued cylindrical Wiener process. On the other hand, the frozen equation (\ref{frozequa}) is written as
\ce
Y_{\frac{u}{\e}}^{\frac{k\d}{\e},x, y}=y+\int_{\frac{k\d}{\e}}^{\frac{u}{\e}}\left[AY_{v}^{\frac{k\d}{\e},x, y}+f_{2}(v, x, Y_{v}^{\frac{k\d}{\e},x, y})\right]\dif v+\int_{\frac{k\d}{\e}}^{\frac{u}{\e}} g_2(v, x,Y_{v}^{\frac{k\d}{\e},x, y})\dif W^2_{v}.
\de
Thus, for $v\in[\frac{k\d}{\e},\frac{(k+1)\d}{\e}]$, $\check{Y}^{\e,x, y}_{\e v}$ and $Y_{v}^{\frac{k\d}{\e},x, y}$ have the same distribution, which implies that
\ce
&&\mE\Bigg[\<f_{1}(\frac{u}{\e},x,\check{Y}^{\e,x, y}_{u})-\bar{f}_{1}(\frac{u}{\e},x), f_{1}(\frac{r}{\e},x,\check{Y}^{\e,x, y}_{r})-\bar{f}_{1}(\frac{r}{\e},x)\>_\mH\Bigg]\\
&=&\mE\<f_{1}(\frac{u}{\e},x,Y_{\frac{u}{\e}}^{\frac{k\d}{\e},x, y})-\bar{f}_{1}(\frac{u}{\e},x), f_{1}(\frac{r}{\e},x,Y_{\frac{r}{\e}}^{\frac{k\d}{\e},x, y})-\bar{f}_{1}(\frac{r}{\e},x)\>_\mH\\
&=&\mE\left[\mE\left[\<f_{1}(\frac{u}{\e},x,Y_{\frac{u}{\e}}^{\frac{k\d}{\e},x, y})-\bar{f}_{1}(\frac{u}{\e},x), f_{1}(\frac{r}{\e},x,Y_{\frac{r}{\e}}^{\frac{k\d}{\e},x, y})-\bar{f}_{1}(\frac{r}{\e},x)\>_\mH\Bigg{|}\sF_{\frac{r}{\e}}^{Y^{\frac{k\d}{\e},x, y}}\right]\right]\\
&=&\mE\left[\Bigg<\mE\left[f_{1}(\frac{u}{\e},x,Y_{\frac{u}{\e}}^{\frac{k\d}{\e},x, y})\Bigg{|}\sF_{\frac{r}{\e}}^{Y^{\frac{k\d}{\e},x, y}}\right]-\bar{f}_{1}(\frac{u}{\e},x),f_{1}(\frac{r}{\e},x,Y_{\frac{r}{\e}}^{\frac{k\d}{\e},x, y})-\bar{f}_{1}(\frac{r}{\e},x)\Bigg>_\mH\right]\\
&\leq&\left(\mE\left|\mE\left[f_{1}(\frac{u}{\e},x,Y_{\frac{u}{\e}}^{\frac{r}{\e},x,Y_{\frac{r}{\e}}^{\frac{k\d}{\e},x,y}})\right]-\bar{f}_{1}(\frac{u}{\e},x)\right|_\mH^2\right)^{1/2}\\
&&\times\left(\mE|f_{1}(\frac{r}{\e},x,Y_{\frac{r}{\e}}^{\frac{k\d}{\e},x, y})-\bar{f}_{1}(\frac{r}{\e},x)|_\mH^2\right)^{1/2},
\de
where $\mathscr{F}_{\frac{r}{\e}}^{Y^{\frac{k\d}{\e},x, y}}=\sigma\{Y_t^{\frac{k\d}{\e},x, y}, \frac{k\d}{\e}\leq t\leq \frac{r}{\e}\}$. 

Besides, by (\ref{evoesti}) and (\ref{frozmomeesti}) it holds that
\ce
&&\left(\mE\left|\mE\left[f_{1}(\frac{u}{\e},x,Y_{\frac{u}{\e}}^{\frac{r}{\e},x,Y_{\frac{r}{\e}}^{\frac{k\d}{\e},x,y}})\right]-\bar{f}_{1}(\frac{u}{\e},x)\right|_\mH^2\right)^{1/2}\\
&\leq&\(Ce^{-2\a\frac{u-r}{\e}}(1+|x|_\mH^2+\mE|Y_{\frac{r}{\e}}^{\frac{k\d}{\e},x,y}|_\mH^2)\)^{\frac{1}{2}}\\
&\leq&\(Ce^{-2\a\frac{u-r}{\e}}(1+|x|_\mH^2+|y|_\mH^2)\)^{\frac{1}{2}}\\
&\leq& Ce^{-\a\frac{u-r}{\e}}(1+|x|_\mH+|y|_\mH),
\de
and by (\ref{frozmomeesti}) and (\ref{muxboun})
\ce
&&\left(\mE|f_{1}(\frac{r}{\e},x,Y_{\frac{r}{\e}}^{\frac{k\d}{\e},x, y})-\bar{f}_{1}(\frac{r}{\e},x)|_\mH^2\right)^{1/2}\\
&=&\left(\mE|f_{1}(\frac{r}{\e},x,Y_{\frac{r}{\e}}^{\frac{k\d}{\e},x, y})-\int_{\mH}f_{1}(\frac{r}{\e},x,z)\mu^{x}_{\frac{r}{\e}}(\dif z)|_\mH^2\right)^{1/2}\\
&\leq&\left(\mE\int_{\mH}|f_{1}(\frac{r}{\e},x,Y_{\frac{r}{\e}}^{\frac{k\d}{\e},x, y})-f_{1}(\frac{r}{\e},x,z)|_\mH^2\mu^{x}_{\frac{r}{\e}}(\dif z)\right)^{1/2}\\
&\leq&L_1\left(\int_{\mH}\mE|Y_{\frac{r}{\e}}^{\frac{k\d}{\e},x, y}-z|_\mH^2\mu^{x}_{\frac{r}{\e}}(\dif z)\right)^{1/2}\\
&\leq&C\left(|y|_\mH^2 e^{-2\a\frac{r-k\d}{\e}}+C(1+|x|_\mH^2)\right)^{1/2}\\
&\leq&C(1+|x|_\mH+|y|_\mH).
\de

The above deduction with (\ref{xehvesti}) and (\ref{hatzub}) yields that
\ce
\Psi(u,r)\leq Ce^{-\a\frac{u-r}{\e}}\left(1+|x_0|_\mH^2+|y_0|_\mH^2\right).
\de
Inserting the above inequality in (\ref{b41c}), we get that
\be
I_{41}&\leq& 2C_R(\frac{T}{\d})\left(1+|x_0|_\mH^2+|y_0|_\mH^2\right)\sup_{0\leq k\leq [\frac{T}{\d}]-1}\left(C\int_{k\d}^{(k+1)\d}\int_{r}^{(k+1)\d}e^{-\a\frac{u-r}{\e}}\dif u\dif r\right)^{1/2}\no\\
&\leq& C_R(\frac{\e}{\d})^{1/2}\left(1+|x_0|_\mH^2+|y_0|_\mH^2\right).
\label{b4de1}
\ee

Next, we deal with $I_{42}$. By (\ref{f1g1linegrow}), (\ref{xehvesti}), (\ref{barxee}), (\ref{hatzub}) and the H\"older inequality, one could get that
\be
I_{42}&\leq&2C_R\mE\sup\limits_{t\in[0,T\wedge \tau_3]}\int_{[\frac{t}{\d}]\d}^t|Z^\e_{r(\d)}|_\mH|f_1(\frac{r}{\e},X^\e_{r(\d)},\hat{Y}^\e_{r})-\bar{f}_1(\frac{r}{\e},X^\e_{r(\d)})|_\mH\dif r\no\\
 &\leq&2C_R\d^{1/2}\left(\mE\sup_{t\in[0,T]}\int_{[\frac{t}{\d}]\d}^{t}|f_1(\frac{r}{\e},X^\e_{r(\d)},\hat{Y}^\e_{r})-\bar{f}_1(\frac{r}{\e},X^\e_{r(\d)})|_\mH^2\dif r\right)^{1/2}\no\\
 &&\times\left(\mE\sup\limits_{r\in[0,T]}|Z^\e_r|_\mH^2\right)^{1/2}\no\\ 
 &\leq&2C_R\d^{1/2}\left(\mE\int_{0}^{T}|f_1(\frac{r}{\e},X^\e_{r(\d)},\hat{Y}^\e_{r})-\bar{f}_1(\frac{r}{\e},X^\e_{r(\d)})|_\mH^2\dif r\right)^{1/2}\no\\
 &&\times\left(\mE\sup\limits_{r\in[0,T]}|Z^\e_r|_\mH^2\right)^{1/2}\no\\ 
 &\leq&C_R\d^{1/2}\left(\int_{0}^{T}(1+\mE|X_{r(\d)}^\e|_\mH^2+\mE|\hat{Y}_{r}^\e|_\mH^2)\dif r\right)^{1/2}\no\\
 &&\times\left(\mE\sup_{r\in[0,T]}|X^\e_{r}|_\mH^2+\mE\sup_{r\in[0,T]}|\bar{X}^\e_{r}|_\mH^2\right)^{1/2}\no\\
&\leq& C_R\d^{1/2}\left(1+|x_0|_\mH^2+|y_0|_\mH^2\right).
\label{i42}
\ee

Finally, combining (\ref{b4de1}) and (\ref{i42}) we conclude (\ref{i4}).
\end{proof}

Now, we are ready to prove Theorem \ref{genestroapth}.

{\bf Proof of Theorem \ref{genestroapth}.} For any $\eta>0$, by the Chebychev inequality and (\ref{xebarxeesti}) it holds that
\ce
&&\mP(\sup\limits_{t\in[0,T]}|X^\e_t-\bar X^\e_t|_\mH\geq \eta)\\
&=&\mP(\sup\limits_{t\in[0,T]}|X^\e_t-\bar X^\e_t|_\mH\geq \eta, T\leq \tau_3)+\mP(\sup\limits_{t\in[0,T]}|X^\e_t-\bar X^\e_t|_\mH\geq \eta, T>\tau_3)\\
&\leq&\mP(\sup\limits_{t\in[0,T\wedge\tau_3]}|X^\e_t-\bar X^\e_t|_\mH\geq \eta)+\mP(T>\tau_3)\\
&\leq&\frac{1}{\eta^2}\mE\sup\limits_{t\in[0,T\wedge\tau_3]}|X^\e_t-\bar X^\e_t|^2_\mH+\frac{1}{R}\mE\left[|X^\e_T|_\mH+|\bar{X}^\e_T|_\mH+\int_0^T\(\|X^\e_{r}\|_\mV^2+\|\bar{X}^\e_{r}\|^2_\mV\)\dif r\right]\\
&\leq&\frac{C_{R}}{\eta^2}\left(\d^{1/2}+\d^{1/4}+(\frac{\e}{\d})^{1/2}\right)\left(1+|x_0|_\mH^2+|y_0|_\mH^2\right)+\frac{C}{R}\left(1+|x_0|_\mH^2+|y_0|_\mH^2\right).
\de
Taking $\d=\e^\t$ for any $\t\in (0,1)$ and letting $\e\rightarrow0$ first and then $R\rightarrow\infty$, we obtain that
\ce
\lim\limits_{\e\rightarrow0}\mP(\sup\limits_{t\in[0,T]}|X^\e_t-\bar X^\e_t|_\mH\geq \eta)=0,
\de
which completes the proof.

\section{Proof of Theorem \ref{asymstroapth}}\label{asymstroapthproo}

In this section, we prove Theorem \ref{asymstroapth}. First, we study the asymptotic behavior for the evolution system $\{\mu_t^{x}\}_{t\geq0}$ of measures by means of an auxiliary SPDE. Then some estimates for the averaged equation (\ref{averequas}) are presented. Finally, we prove Theorem \ref{asymstroapth}.

\subsection{Asymptotic behavior for evolution system of measures}\label{asyevo}

In this subsection, we study the asymptotic behavior for the evolution system of measures.

In terms of $(\mathbf{H}^1_{f_{2}, g_{2}})$ and $(\mathbf{H}^a_{f_{2}, g_{2}})$, it holds that for any $t\geq 0$ and $x_1, x_2, y_1, y_2\in\mH$,
\ce
|\bar{f}_2(x_1,y_1)-\bar{f}_2(x_2,y_2)|_\mH&\leq& |\bar{f}_2(x_1,y_1)-f_2(t, x_1,y_1)|_\mH+|f_2(t, x_1,y_1)-f_2(t, x_2,y_2)|_\mH\\
&&+|f_2(t,x_2,y_2)-\bar{f}_2(x_2,y_2)|_\mH\\
&\leq&\kappa_3(t)(2+|x_1|_\mH+|y_1|_\mH+|x_2|_\mH+|y_2|_\mH)+L_3|x_1-x_2|_\mH\\
&&+L_4|y_1-y_2|_\mH.
\de
Letting $t\rightarrow\infty$ and noticing $\lim\limits_{t\rightarrow\infty}\kappa_3(t)=0$, we obtain that
\be
|\bar{f}_2(x_1,y_1)-\bar{f}_2(x_2,y_2)|_\mH\leq L_3|x_1-x_2|_\mH+L_4|y_1-y_2|_\mH.
\label{barf2lip}
\ee
By the same deduction to that for $\bar{f}_2$, we infer that
\be
\|\bar{g}_2(x_1,y_1)-\bar{g}_2(x_2,y_2)\|_{\cL_2(\mU_2,\mH)}\leq L_5|x_1-x_2|_\mH+L_6|y_1-y_2|_\mH.
\label{barg2lip}
\ee
Moreover, it is easy to see that (\ref{discond1}) and (\ref{discond2}) hold for $\bar{f}_2, \bar{g}_2$.

Consequently, the following SPDE
\be\left\{\begin{array}{ll}
\dif \bar{Y}_t=\left[A\bar{Y}_t+\bar{f}_2(x,\bar{Y}_t)\right]\dif t+\bar{g}_2(x,\bar{Y}_t)\dif W^2_t,\\
\bar{Y}_t(0)=\bar{Y}_t(1)=0,  \quad t\geq 0,\\
\bar{Y}_0=y_0,
\end{array}
\right.
\label{frozeqlimi}
\ee
admits a unique solution $\{\bar{Y}_t^{x,y_0}\}_{t\geq0}$ such that $\mE[\sup\limits_{t\in[0,T]}|\bar Y^{x,y_0}_t|_\mH^2+\int_0^T\|\bar Y^{x,y_0}_t\|^2_\mV\dif t]<\infty$ for any $T>0$ (\cite[Theorem 7.5]{DaZ}). Moreover, by the argument used in the proofs of (\ref{frozmomeesti}) and (\ref{frozdiffesti}), one can obtain that for any $t\geq0$ and $x_1, x_2, y_1, y_2\in\mH$,
\be
&&\mE|\bar{Y}_t^{x,y_0}|_\mH^2\leq C(1+|x|_\mH^2)+|y_0|^2_\mH e^{-2\a t},\label{baryes}\\
&&\mE|\bar{Y}_t^{x_1,y_1}-\bar{Y}_t^{x_2,y_2}|_\mH^2\leq C|x_1-x_2|^2_\mH+|y_1-y_2|^2_\mH e^{-2\alpha t}. \label{baryxyes}
\ee

Let $\{\bar{P}_t^{x}\}_{t\geq0}$ denote the transition semigroup associated with $\{\bar{Y}_t^{x,y_0}\}_{t\geq0}$, i.e., for any bounded measurable function $h:\mH\rightarrow\mR$,
\ce
\bar{P}_t^xh(y_0)=\mE h(\bar{Y}_t^{x,y_0}).
\de
Then it admits a unique invariant measure $\mu^x$ satisfying
\be
\int_{\mH}|y|^2\mu^x(\dif y)\leq C(1+|x|^2_\mH),
\label{musees}
\ee
which with (\ref{baryxyes}) yields that for any Lipschitz function $h$ on $\mH$,
\be
\left|\mE h(\bar{Y}_t^{x,y_0})-\int_\mH h(y)\mu^x(\dif y)\right|\leq C\|h\|_{Lip}(1+|x|_\mH+|y_0|_\mH)e^{-\a t}.
\label{conbarYnu}
\ee

\bl\label{mutmu}
Suppose that $(\mathbf{H}^1_{f_{2},g_{2}})$, $(\mathbf{H}^2_{f_{2},g_{2}})$ and $(\mathbf{H}^a_{f_{2}, g_{2}})$ hold. Then for any Lipschitz function $h$,
\be
&&\left|\int_\mH h(y)\mu_t^x(\dif y)-\int_\mH h(y)\mu^x(\dif y)\right|\no\\
&\leq& C\|h\|_{Lip}\left[(1+|x|_\mH+|y_0|_\mH)e^{-\a t}+(1+|x|_\mH)\left(\int_0^te^{-2\a(t-r)}\kappa_3^2(r)\dif r\right)^{1/2}\right],
\label{nutnuesti}
\ee
where $\{\mu_t^{x}\}_{t\geq0}$ is the evolution system of measures of the semigroup $\{P_{s,t}^{x}\}_{t\geq s}$ given in Lemma \ref{mux}.
\el
\begin{proof}
First of all, we estimate $\mE|Y_t^x-\bar Y_t^{x,y_0}|_\mH^2$. By the It\^o formula, it holds that for any $t\geq 0$ and $l_4>0$
\ce
e^{l_4 t}|Y_t^x-\bar Y_t^{x,y_0}|_\mH^2&=&|Y_0^x-y_0|_\mH^2+l_4\int_0^te^{l_4 r}|Y_r^x-\bar Y_r^{x,y_0}|_\mH^2\dif r\\
&&+ 2\int_0^te^{l_4 r}{_{\mV^*}}\<AY_{r}^x-A\bar{Y}_{r}^{x,y_0}, Y_r^x-\bar Y_r^{x,y_0}\>_\mV\dif r\\
&&+2\int_0^te^{l_4 r}\<f_2(r,x,Y_{r}^x)-\bar f_2(x,\bar Y_r^{x,y_0}), Y_r^x-\bar Y_r^{x,y_0}\>_\mH\dif r\\
&&+2\int_0^te^{l_4 r}\<Y_r^x-\bar Y_r^{x,y_0}, \(g_2(r,x,Y_{r}^x)-\bar g_2(x,\bar Y_r^{x,y_0})\)\dif W^2_r\>_\mH\\
&&+\int_0^te^{l_4 r}\|g_2(r,x,Y_{r}^x)-\bar g_2(x,\bar Y_r^{x,y_0})\|^2_{\cL_2(\mU_2,\mH)}\dif r.
\de
By the definition of $A$ and the Poincar\'e inequality, we know that
\ce
2\<AY_{r}^x-A\bar{Y}_{r}^{x,y_0}, Y_r^x-\bar Y_r^{x,y_0}\>_\mV\leq -2\l_1|Y_r^x-\bar Y_r^{x,y_0}|^2_\mH.
\de
Then $(\mathbf{H}^a_{f_{2}, g_{2}})$, (\ref{barf2lip}) and the Young inequality imply that
\ce
&&2\<f_2(r,x,Y_{r}^x)-\bar f_2(x,\bar Y_r^{x,y_0}), Y_r^x-\bar Y_r^{x,y_0}\>_\mH\\
&=&2\<f_2(r,x,Y_{r}^x)-\bar f_2(x,Y_{r}^x), Y_r^x-\bar Y_r^{x,y_0}\>_\mH+2\<\bar f_2(x,Y_{r}^x)-\bar f_2(x,\bar Y_r^{x,y_0}), Y_r^x-\bar Y_r^{x,y_0}\>_\mH\\
&\leq&C|f_2(r,x,Y_{r}^x)-\bar f_2(x,Y_{r}^x)|_\mH^2+L_4|Y_r^x-\bar Y_r^{x,y_0}|_\mH^2+2L_4|Y_r^x-\bar Y_r^{x,y_0}|_\mH^2\\
&\leq&C\kappa^2_3(r)(1+|x|_\mH^2+|Y_{r}^x|_\mH^2)+3L_4|Y_r^x-\bar Y_r^{x,y_0}|_\mH^2.
\de
In terms of $(\mathbf{H}^a_{f_{2}, g_{2}})$ and (\ref{barg2lip}), the similar argument yields that
\ce
\|g_2(r,x,Y_{r}^x)-\bar g_2(x,\bar Y_r^{x,y_0})\|^2_{\cL_2(\mU_2,\mH)}\leq 6\kappa^2_3(r)(1+|x|_\mH^2+|Y_{r}^x|_\mH^2)+2L_6^2|Y_r^x-\bar Y_r^{x,y_0}|_\mH^2.
\de

Collecting the above deduction, we conclude that
\ce
 e^{l_4 t}\mE|Y_t^x-\bar Y_t^{x,y_0}|_\mH^2&\leq&\mE|Y_0^x-y_0|_\mH^2+C\int_0^te^{l_4 r}\kappa^2_3(r)(1+|x|_\mH^2+\mE|Y_{r}^x|_\mH^2)\dif r\\
&&+(l_4-2\a)\int_0^te^{l_4 r}\mE|Y_r^x-\bar Y_r^{x,y_0}|_\mH^2\dif r.
\de
Taking $l_4=2\a$, by (\ref{yxboun}) one can obtain that
\ce
\mE|Y_t^x-\bar Y_t^{x,y_0}|_\mH^2\leq C(1+|x|^2_\mH+|y_0|^2_\mH)e^{-2\a t}+C(1+|x|^2_\mH)\int_0^te^{-2\a (t-r)}\kappa^2_3(r)\dif r.
\de

Finally, by (\ref{conbarYnu}) and the above inequality, it holds that for any Lipschitz function $h$,
\ce
&&\left|\int_\mH h(y)\mu_t^x(\dif y)-\int_\mH h(y)\mu^x(\dif y)\right|\\
&\leq&\left|\mE h(Y_t^x)-\mE h(\bar{Y}_t^{x,y_0})\right|+\left|\mE h(\bar{Y}_t^{x,y_0})-\int_\mH h(y)\mu^x(\dif y)\right|\\
&\leq&\|h\|_{Lip}\mE|Y_t^x-\bar{Y}_t^{x,y_0}|_\mH+C\|h\|_{Lip}(1+|x|+|y_0|_\mH)e^{-\a t}\\
&\leq&C\|h\|_{Lip}(1+|x|+|y_0|_\mH)e^{-\a t}+C\|h\|_{Lip}(1+|x|_\mH)\left[\int_0^te^{-2\a(t-r)}\kappa_3^2(r)\dif r\right]^{1/2}.
\de
The proof is complete.
\end{proof}

\br
As shown in \cite[Remark 4.2]{SWX1},  
\ce
\lim\limits_{t\to\infty}\int_0^t e^{-2\alpha(t-r)}\kappa_3^2(r)\dif r = 0.
\de  
This vanishing integral implies that the time-dependent measure $\mu^x_t$ converges weakly to the measure $\mu^x$ as $t\to\infty$. That is, the ``invariant measure" $\mu^x_t$ for the nonautonomous SPDE (\ref{frozequa}) converges weakly to the invariant measure $\mu^x$ for the autonomous SPDE (\ref{frozeqlimi})  as $t\to\infty$. This is reasonable since $f_2(t,x,y)\to \bar f_2(x,y)$ and $g_2(t,x,y)\to \bar g_2(x,y)$.
\er

Finally, by Lemma \ref{mutmu}, we compute $\left|\frac{1}{T}\int_t^{t+T}\bar f_1(r,x)\dif r-\check f_{1,a}(x)\right|_\mH$ which is used in the next subsection. Note that for any $T>0, t\geq 0$ and $x\in\mH$
\ce
\frac{1}{T}\int_t^{t+T}\bar f_1(r,x)\dif r=\frac{1}{T}\int_t^{t+T}\int_\mH f_1(r,x,y)\mu^x_r(\dif y)\dif r, \quad \check f_{1,a}(x)=\int_\mH \check f_1(x,y)\mu^x(\dif y).
\de
Thus, by (\ref{nutnuesti}) and $(\mathbf{H}^a_{f_{1}})$ it holds that
\ce
&&\left|\frac{1}{T}\int_t^{t+T}\bar f_1(r,x)\dif r-\check f_{1,a}(x)\right|_\mH\\
&\leq&\left|\frac{1}{T}\int_t^{t+T}\int_\mH f_1(r,x,y)\mu^x_r(\dif y)\dif r-\frac{1}{T}\int_t^{t+T}\int_\mH f_1(r,x,y)\mu^x(\dif y)\dif r\right|_\mH\\
&&+\left|\int_\mH\left(\frac{1}{T}\int_t^{t+T} f_1(r,x,y)\dif r\right)\mu^x(\dif y)-\int_\mH \check f_1(x,y)\mu^x(\dif y)\right|_\mH\\
&\leq&CL_1(1+|x|+|y_0|_\mH)\frac{1}{T}\int_t^{t+T}e^{-\a r}\dif r\\
&&+CL_1(1+|x|_\mH)\frac{1}{T}\int_t^{t+T}\left[\int_0^re^{-2\a(r-v)}\kappa_3^2(v)\dif v\right]^{1/2}\dif r\\
&&+\int_\mH \kappa_1(T)(1+|x|_\mH+|y|_\mH)\mu^x(\dif y)\\
&\leq&C(1+|x|_\mH+|y_0|_\mH)\left[\frac{1}{T}\int_t^{t+T}e^{-\a r}\dif r+\frac{1}{T}\int_t^{t+T}\left[\int_0^re^{-2\a(r-v)}\kappa_3^2(v)\dif v\right]^{1/2}\dif r+\kappa_1(T)\right],
\de
where we use (\ref{musees}) in the last inequality. 

Set 
$$
\kappa^*(T):=\left[\frac{1}{T}\int_t^{t+T}e^{-\a r}\dif r+\frac{1}{T}\int_t^{t+T}\left[\int_0^re^{-2\a(r-v)}\kappa_3^2(v)\dif v\right]^{1/2}\dif r+\kappa_1(T)\right],
$$
and the above argument implies that
\be
\left|\frac{1}{T}\int_t^{t+T}\bar f_1(r,x)\dif r-\check f_{1,a}(x)\right|_\mH\leq C(1+|x|_\mH+|y_0|_\mH)\kappa^*(T).
\label{barf1checkf1a}
\ee

For $\kappa^*$, on the one hand, since $\kappa_1, \kappa_3$ are locally bounded, by some simple calculation it is easy to see that it  is locally bounded. On the other hand, the H\"older inequality and the Fubini theorem imply that
\ce
\frac{1}{T}\int_t^{t+T}\left[\int_0^re^{-2\a(r-v)}\kappa_3^2(v)\dif v\right]^{1/2}\dif r&\leq& \left[\frac{1}{T}\int_t^{t+T}\int_0^re^{-2\a(r-v)}\kappa_3^2(v)\dif v\dif r\right]^{1/2}\\
&\leq& \left[\frac{1}{T}\int_t^{t+T} \kappa_3^2(v)\dif v\right]^{1/2}.
\de
In terms of $(\mathbf{H}^a_{f_{2}, g_{2}})$, we conclude that $\kappa^*$ satisfies $\lim\limits_{T\rightarrow\infty}\kappa^*(T)=0$. 

\subsection{Some estimates for the averaged equation (\ref{averequas})}

In this subsection, we present some estimates for the averaged equation (\ref{averequas}).

\bl\label{checkxwellpose}
Assume that $(\mathbf{H}_{f_{1}, g_{1}})$, $(\mathbf{H}^1_{f_{2}, g_{2}})$, $(\mathbf{H}^2_{f_{2}, g_{2}})$, $(\mathbf{H}^a_{f_{1}})$, $(\mathbf{H}^a_{g_{1}})$ and $(\mathbf{H}^a_{f_{2}, g_{2}})$ hold. Then Eq.(\ref{averequas}) has a unique solution $(\check{X}, \check{K})$ such that 
\be
\mE\left[\sup _{t \in[0, T]}\left|\check X_t\right|_\mH^2+\int_0^T\left\|\check X_t\right\|_\mV^2\dif t\right]\leq C\left(1+|x_0|_\mH^2\right).
\label{checkxesti}
\ee
Moreover, set the stopping time $\check\tau_1:=\inf\{t>0, |\check X_t|_\mH>R\}$, and for any $\d>0$ small enough,
\be
\mE\int_0^{T\wedge\check\tau_1}|\check X_t-\check X_{t(\d)}|_\mH^2\dif t\leq C\d^{1/2}\left(1+|x_0|_\mH^2\right).
\label{checkxincr}
\ee
\el
\begin{proof}
First of all, we prove that $\check{f}_{1}$ is Lipschitz continuous. Indeed, by $(\mathbf{H}_{f_{1}, g_{1}})$ and $(\mathbf{H}^a_{f_{1}})$ it holds that for any $T>0$, $t\geq 0$ and $x_1, x_2, y_1, y_2\in\mH$, 
\be
|\check f_1(x_1,y_1)-\check f_1(x_2,y_2)|_\mH&\leq& \left|\check f_1(x_1,y_1)-\frac{1}{T}\int_t^{t+T} f_1(s,x_1,y_1)\dif s\right|_\mH\no\\
&&+\left|\frac{1}{T}\int_t^{t+T} f_1(s,x_1,y_1)\dif s-\frac{1}{T}\int_t^{t+T} f_1(s,x_2,y_2)\dif s\right|_\mH\no\\
&&+\left|\frac{1}{T}\int_t^{t+T} f_1(s,x_2,y_2)\dif s-\check f_1(x_2,y_2)\right|_\mH\no\\
&\leq&\kappa_1(T)(1+|x_1|_\mH+|y_1|_\mH+|x_2|_\mH+|y_2|_\mH)\no\\
&&+L_1(|x_1-x_2|_\mH+|y_1-y_2|_\mH).
\label{checkf1plip}
\ee
Letting $T\rightarrow\infty$, we obtain that
\ce
|\check f_1(x_1,y_1)-\check f_1(x_2,y_2)|_\mH\leq L_1(|x_1-x_2|_\mH+|y_1-y_2|_\mH).
\de
Since $\check f_{1,a}(x)=\int_\mH \check f_1(x,y)\mu^x(\dif y)$, in terms of (\ref{conbarYnu}) and (\ref{baryxyes}) the same deduction to that in \cite[Lemma 4.3]{Q1} yields that
\ce
|\check f_{1,a}(x_1)-\check f_{1,a}(x_2)|_\mH\leq C|x_1-x_2|_\mH.
\de
By $(\mathbf{H}_{f_{1}, g_{1}})$ and $(\mathbf{H}^a_{g_{1}})$, the same deduction to that for $\check{f}_{1}$ yields that
\be
|\bar{g}_1(x_1)-\bar{g}_1(x_2)|_\mH\leq C|x_1-x_2|_\mH.
\label{barg1plip}
\ee
Thus, by Theorem 3.2 in \cite{Q2} or Theorem 3.1 in \cite{tZ}, we obtain that Eq.(\ref{averequas}) has a unique solution $(\check{X}, \check{K})$. Then since the proofs of (\ref{checkxesti}) and (\ref{checkxincr}) are similar to that for (\ref{xehvesti}) and (\ref{xeettd}), respectively, we omit them.
\end{proof}

Next, we investigate the difference between $\bar X^\e$ and $\check X$. In order to do this, we define the stopping time
\ce
\check\tau_2:=\inf\left\{t>0, |\bar X^\e_t|_\mH+|\check X_t|_\mH+\int_0^t\(\|\bar X^\e_{r}\|_\mV^2+\|\check{X}_{r}\|^2_\mV\)\dif r>R\right\}.
\de
Then $\check\tau_2\leq\tau_2\wedge\check\tau_1$.

\bp\label{barxecheckxespr}
Suppose that $(\mathbf{H}_{f_{1}, g_{1}})$, $(\mathbf{H}^1_{f_{2}, g_{2}})$, $(\mathbf{H}^2_{f_{2}, g_{2}})$, $(\mathbf{H}^a_{f_{1}})$, $(\mathbf{H}^a_{g_{1}})$ and $(\mathbf{H}^a_{f_{2}, g_{2}})$ hold. Then for any $T>0$
\be
\mE\sup\limits_{t\in[0,T\wedge\check\tau_2]}|\bar X^\e_t-\check X_t|^2_\mH\leq C_R\left(\d+\d^{1/2}+\d^{1/4}+\kappa^*(\frac{\d}{\e})+\kappa_2(\frac{\d}{\e})\right)(1+|x_0|^2_\mH+|y_0|^2_\mH).
\label{barxecheckxesti}
\ee
\ep
\begin{proof}
We divide the proof into two steps. In the first step, we prove the required estimate (\ref{barxecheckxesti}). Two inequalities in the first step are proved in the second step.

{\bf Step 1.} We prove the required estimate (\ref{barxecheckxesti}).

Set $\bar Z^\e_t:=\bar X^\e_{t}-\check{X}_{t}$ and by the It\^o formula we get that
\ce
|\bar Z^\e_t|_\mH^2&=&2\int_{0}^t{_{\mV^*}}\<A\bar X^\e_{r}-A\check{X}_{r},\bar Z^\e_r\>_\mV\dif r+2\int_{0}^t{_{\mV^*}}\<B(\bar X^\e_{r},\bar X^\e_{r})-B(\check{X}_{r},\check{X}_{r}),\bar Z^\e_r\>_\mV\dif r\no\\
&&+2\int_{0}^t\<\bar Z^\e_r, \bar f_1(\frac{r}{\e},\bar X^\e_{r})-\check{f}_{1,a}(\check{X}_{r})\>_\mH\dif r+2 \int_{0}^t\<\bar Z^\e_r, (g_1(\frac{r}{\e},\bar X^\e_{r})-\check g_1(\check X_r))\dif W^1_r\>_\mH \\
&&+2\int_{0}^t\int_0^1\bar Z^\e_r(\xi)\bar K^\e(\dif r,\dif \xi)-2\int_{0}^{t}\int_0^1\bar Z^\e_r(\xi)\check{K}(\dif r,\dif\xi)\no\\
&&+\int_{0}^{t}\|g_1(\frac{r}{\e},\bar X^\e_{r})-\check g_1(\check X_r)\|_{\cL_2(\mU_1,\mH)}^{2} \dif r.
\de
By the similar deduction to that in Proposition \ref{xebarxees}, it holds that
\ce
&&2{_{\mV^*}}\<A\bar X^\e_{r}-A\check{X}_{r},\bar Z^\e_r\>_\mV=-2\|\bar Z^\e_r\|^2_\mV,\\
&&2{_{\mV^*}}\<B(\bar X^\e_{r},\bar X^\e_{r})-B(\check{X}_{r},\check{X}_{r}),\bar Z^\e_r\>_\mV\leq 8|\bar Z^\e_r|^2_\mH(\|\bar X^\e_{r}\|^2_\mV+\|\check X_{r}\|^2_\mV)+\|\bar Z^\e_r\|^2_\mV,\\
&&2\int_{0}^t\int_0^1\bar Z^\e_r(\xi)\bar K^\e(\dif r,\dif \xi)-2\int_{0}^{t}\int_0^1\bar Z^\e_r(\xi)\check{K}(\dif r,\dif\xi)\leq 0.
\de
Then the Lipschitz continuity of $\bar f_1, \check{f}_{1,a}, g_1, \bar g_1$ imply that
\ce
&&2\<\bar Z^\e_r, \bar f_1(\frac{r}{\e},\bar X^\e_{r})-\check{f}_{1,a}(\check{X}_{r})\>_\mH\\
&=&2\<\bar Z^\e_r, \bar f_1(\frac{r}{\e},\bar X^\e_{r})-\bar f_1(\frac{r}{\e},\check X_{r})\>_\mH+2\<\bar Z^\e_r, \bar f_1(\frac{r}{\e},\check X_{r})-\bar f_1(\frac{r}{\e},\check X_{r(\d)})\>_\mH\\
&&+2\<\bar Z^\e_r, \check{f}_{1,a}(\check{X}_{r(\d)})-\check{f}_{1,a}(\check{X}_{r})\>_\mH+2\<\bar Z^\e_r-\bar Z^\e_{r(\d)},\bar f_1(\frac{r}{\e},\check X_{r(\d)})-\check{f}_{1,a}(\check{X}_{r(\d)})\>_\mH\\
&&+2\<\bar Z^\e_{r(\d)},\bar f_1(\frac{r}{\e},\check X_{r(\d)})-\check{f}_{1,a}(\check{X}_{r(\d)})\>_\mH\\
&\leq&2C|\bar Z^\e_r|^2_\mH+C|\check X_{r}-\check X_{r(\d)}|^2_\mH+C(|\bar X^\e_{r}-\bar X^\e_{r(\d)}|_\mH+|\check X_{r}-\check X_{r(\d)}|_\mH)(1+|\check{X}_{r(\d)}|)\\
&&+2\<\bar Z^\e_{r(\d)},\bar f_1(\frac{r}{\e},\check X_{r(\d)})-\check{f}_{1,a}(\check{X}_{r(\d)})\>_\mH,
\de
and
\be
&&\|g_1(\frac{r}{\e},\bar X^\e_{r})-\bar g_1(\check X_r)\|_{\cL_2(\mU_1,\mH)}^{2}\no\\
&\leq& 4\|g_1(\frac{r}{\e},\bar X^\e_{r})-g_1(\frac{r}{\e},\check X_{r})\|_{\cL_2(\mU_1,\mH)}^{2}+ 4\|g_1(\frac{r}{\e},\check X_{r})-g_1(\frac{r}{\e},\check X_{r(\d)})\|_{\cL_2(\mU_1,\mH)}^{2}\no\\
&&+4\|g_1(\frac{r}{\e},\check X_{r(\d)})-\bar g_1(\check X_{r(\d)})\|_{\cL_2(\mU_1,\mH)}^{2}+4\|\bar g_1(\check X_{r(\d)})-\bar g_1(\check X_r)\|_{\cL_2(\mU_1,\mH)}^{2}\no\\
&\leq& 4L^2_2|\bar Z^\e_r|^2_\mH+C|\check X_{r}-\check X_{r(\d)}|^2_\mH+4\|g_1(\frac{r}{\e},\check X_{r(\d)})-\bar g_1(\check X_{r(\d)})\|_{\cL_2(\mU_1,\mH)}^{2}.
\label{g1deco}
\ee
Thus, collecting the above deduction, we infer that
\ce
&&|\bar Z^\e_t|^{2}_\mH+\int_{0}^t\|\bar Z^\e_r\|_\mV^2\dif r\\
&\leq& C\int_{0}^t|\bar Z^\e_r|_\mH^2(1+\|\bar X^\e_{r}\|^2_\mV+\|\check{X}_{r}\|^2_\mV)\dif r+C\int_{0}^t|\check X_{r}-\check X_{r(\d)}|^2_\mH\dif r\\
&&+C\int_{0}^t(|\bar X^\e_{r}-\bar X^\e_{r(\d)}|_\mH+|\check X_{r}-\check X_{r(\d)}|_\mH)(1+|\check{X}_{r(\d)}|)\dif r\\
&&+2\left|\int_{0}^t\<\bar Z^\e_{r(\d)},\bar f_1(\frac{r}{\e},\check X_{r(\d)})-\check{f}_{1,a}(\check{X}_{r(\d)})\>_\mH\dif r\right|\\
&&+2 \left|\int_{0}^t\<\bar Z^\e_r, (g_1(\frac{r}{\e},\bar X^\e_{r})-\bar g_1(\check X_r))\dif W^1_r\>_\mH\right|\\
&&+4\int_{0}^{t}\|g_1(\frac{r}{\e},\check X_{r(\d)})-\bar g_1(\check X_{r(\d)})\|_{\cL_2(\mU_1,\mH)}^{2}\dif r,
\de
which together with the definition of $\bar\tau_2$ and the Gronwall inequality implies that
\ce
&&\mE\sup\limits_{t\in[0,T\wedge\bar\tau_2]}|\bar Z^\e_t|^{2}_\mH+\mE\int_{0}^{T\wedge\bar\tau_2}\|\bar Z^\e_r\|_\mV^2\dif r\\
&\leq&C_R\mE\int_{0}^{T\wedge\bar\tau_2}|\check X_{r}-\check X_{r(\d)}|^2_\mH\dif r\\
&&+C_R\mE\int_{0}^{T\wedge\bar\tau_2}(|\bar X^\e_{r}-\bar X^\e_{r(\d)}|_\mH+|\check X_{r}-\check X_{r(\d)}|_\mH)(1+|\check{X}_{r(\d)}|)\dif r\\
&&+C_R\mE\sup\limits_{t\in[0,T\wedge\bar\tau_2]}\left|\int_{0}^t\<\bar Z^\e_{r(\d)},\bar f_1(\frac{r}{\e},\check X_{r(\d)})-\check{f}_{1,a}(\check{X}_{r(\d)})\>_\mH\dif r\right|\\
&&+C_R\mE\sup\limits_{t\in[0,T\wedge\bar\tau_2]}\left|\int_{0}^t\<\bar Z^\e_r, (g_1(\frac{r}{\e},\bar X^\e_{r})-\bar g_1(\check X_r))\dif W^1_r\>_\mH\right|\\
&&+C_R\mE\int_{0}^{T\wedge\bar\tau_2}\|g_1(\frac{r}{\e},\check X_{r(\d)})-\bar g_1(\check X_{r(\d)})\|_{\cL_2(\mU_1,\mH)}^{2}\dif r.
\de
Besides, the Burkholder-Davis-Gundy inequality and the Young inequality admit us to obtain that
\ce
&&C_R\mE\sup\limits_{t\in[0,T\wedge\bar\tau_2]}\left|\int_{0}^t\<\bar Z^\e_r, (g_1(\frac{r}{\e},\bar X^\e_{r})-\bar g_1(\check X_r))\dif W^1_r\>_\mH\right|\\
&\leq&C_R\mE\left(\int_{0}^{T\wedge\bar\tau_2}|\bar Z^\e_r|^2_\mH\|g_1(\frac{r}{\e},\bar X^\e_{r})-\bar g_1(\check X_r)\|_{\cL_2(\mU_1,\mH)}^{2} \dif r\right)^{1/2}\\
&\leq&\frac{1}{2}\mE\sup\limits_{t\in[0,T\wedge\bar\tau_2]}|\bar Z^\e_t|^{2}_\mH+C_R\mE\int_{0}^{T\wedge\bar\tau_2}\|g_1(\frac{r}{\e},\bar X^\e_{r})-\bar g_1(\check X_r)\|_{\cL_2(\mU_1,\mH)}^{2} \dif r\\
&\leq&\frac{1}{2}\mE\sup\limits_{t\in[0,T\wedge\bar\tau_2]}|\bar Z^\e_t|^{2}_\mH+C_R\mE\int_{0}^{T\wedge\bar\tau_2}|\bar Z^\e_r|^{2}_\mH\dif r+C_R\mE\int_{0}^{T\wedge\bar\tau_2}|\check X_{r}-\check X_{r(\d)}|^2_\mH\dif r\\
&&+C_R\mE\int_{0}^{T\wedge\bar\tau_2}\|g_1(\frac{r}{\e},\check X_{r(\d)})-\bar g_1(\check X_{r(\d)})\|_{\cL_2(\mU_1,\mH)}^{2}\dif r,
\de
where the last inequality is from (\ref{g1deco}). So, according to the above deduction we conclude that
\be
&&\mE\sup\limits_{t\in[0,T\wedge\bar\tau_2]}|\bar Z^\e_t|^{2}_\mH+\mE\int_{0}^{T\wedge\bar\tau_2}\|\bar Z^\e_r\|_\mV^2\dif r\no\\
&\leq&C_R\mE\int_{0}^{T\wedge\bar\tau_2}|\check X_{r}-\check X_{r(\d)}|^2_\mH\dif r\no\\
&&+C_R\mE\int_{0}^{T\wedge\bar\tau_2}(|\bar X^\e_{r}-\bar X^\e_{r(\d)}|_\mH+|\check X_{r}-\check X_{r(\d)}|_\mH)(1+|\check{X}_{r(\d)}|)\dif r\no\\
&&+C_R\mE\sup\limits_{t\in[0,T\wedge\bar\tau_2]}\left|\int_{0}^t\<\bar Z^\e_{r(\d)},\bar f_1(\frac{r}{\e},\check X_{r(\d)})-\check{f}_{1,a}(\check{X}_{r(\d)})\>_\mH\dif r\right|\no\\
&&+C_R\mE\int_{0}^{T\wedge\bar\tau_2}\|g_1(\frac{r}{\e},\check X_{r(\d)})-\bar g_1(\check X_{r(\d)})\|_{\cL_2(\mU_1,\mH)}^{2}\dif r\no\\
&&+C_R\int_{0}^{T}\mE\sup\limits_{t\in[0,r\wedge\bar\tau_2]}|\bar Z^\e_t|^{2}_\mH\dif r\no\\
&=:&J_1+J_2+J_3+J_4+C_R\int_{0}^{T}\mE\sup\limits_{t\in[0,r\wedge\bar\tau_2]}|\bar Z^\e_t|^{2}_\mH\dif r.
\label{j1234}
\ee

Next, (\ref{checkxincr}) implies that
\be
J_1\leq C_R\d^{1/2}(1+|x_0|^2_\mH).
\label{j1}
\ee
For $J_2$, in terms of (\ref{barxeeincr}), (\ref{checkxincr}), (\ref{checkxesti}) and the H\"older inequality, it holds that
\be
J_2&\leq& C_R\left(\mE\int_{0}^{T\wedge\bar\tau_2}(|\bar X^\e_{r}-\bar X^\e_{r(\d)}|^2_\mH+|\check X_{r}-\check X_{r(\d)}|^2_\mH)\dif r\right)^{1/2}\no\\
&&\times\left(\mE\int_{0}^{T\wedge\bar\tau_2}(1+|\check{X}_{r(\d)}|^2)\dif r\right)^{1/2}\no\\
&\leq& C_R\d^{1/4}(1+|x_0|^2_\mH).
\label{j2}
\ee
For $J_3$ and $J_4$, the result in {\bf Step 2} yields that
\be
&&J_3\leq C_R\left(\kappa^*(\frac{\d}{\e})+\d\right)(1+|x_0|^2_\mH+|y_0|^2_\mH),\label{j3}\\
&&J_4\leq C_R\left(\kappa_2(\frac{\d}{\e})+\d\right)(1+|x_0|^2_\mH).\label{j4}
\ee
Finally, combining (\ref{j1})-(\ref{j4}) with (\ref{j1234}) and using the Gronwall inequality, we establish (\ref{barxecheckxesti}).

{\bf Step 2.} We prove (\ref{j3}) and (\ref{j4}).

In order to prove (\ref{j3}), we observe that
\ce
J_3&\leq&C_R\mE\sup\limits_{t\in[0,T\wedge\bar\tau_2 ]}\left|\int_{0}^{[\frac{t}{\d}]\d}\<\bar Z^\e_{r(\d)},\bar f_1(\frac{r}{\e},\check X_{r(\d)})-\check{f}_{1,a}(\check X_{r(\d)})\>_\mH\dif r\right|\\
&&+C_R\mE\sup\limits_{t\in[0,T\wedge\bar\tau_2]}\left|\int_{[\frac{t}{\d}]\d}^{t}\<\bar Z^\e_{r(\d)},\bar f_1(\frac{r}{\e},\check X_{r(\d)})-\check{f}_{1,a}(\check X_{r(\d)})\>_\mH\dif r\right|\\
&=:&J_{31}+J_{32}.
\de
For $J_{31}$, the H\"older inequality implies that
\ce
J_{31}&=&C_R\mE\sup\limits_{t\in[0,T\wedge\bar\tau_2]}\left|\sum\limits_{k=0}^{[\frac{t}{\d}]-1}\int_{k\d}^{(k+1)\d}\<\bar Z^\e_{r(\d)},\bar f_1(\frac{r}{\e},\check X_{r(\d)})-\check{f}_{1,a}(\check X_{r(\d)})\>_\mH\dif r\right|\\
&\leq&C_R\mE\sup\limits_{t\in[0,T\wedge\bar\tau_2]}\sum\limits_{k=0}^{[\frac{t}{\d}]-1}\left|\bar Z^\e_{k\d}\right|_\mH\left|\int_{k\d}^{(k+1)\d}(\bar f_1(\frac{r}{\e},\check X_{k\d})-\check{f}_{1,a}(\check X_{k\d}))\dif r\right|_\mH\\
&\leq&C_R\sum\limits_{k=0}^{[\frac{T}{\d}]-1}\mE\left|\bar Z^\e_{k\d}\right|_\mH\left|\int_{k\d}^{(k+1)\d}(\bar f_1(\frac{r}{\e},\check X_{k\d})-\check{f}_{1,a}(\check X_{k\d}))\dif r\right|_\mH \\
&\leq&C_R\frac{T}{\d}\sup\limits_{0\leq k\leq [\frac{T}{\d}]-1}\left(\mE\left|\bar Z^\e_{k\d}\right|^2_\mH\right)^{1/2}\left(\mE\left|\int_{k\d}^{(k+1)\d}(\bar f_1(\frac{r}{\e},\check X_{k\d})-\check{f}_{1,a}(\check X_{k\d}))\dif r\right|^2_\mH \right)^{1/2}.
\de
On the one hand, by (\ref{barxee}) and (\ref{checkxesti}) we have that
\ce
\left(\mE\left|\bar Z^\e_{k\d}\right|^2_\mH\right)^{1/2}\leq C(1+|x_0|^2_\mH)^{1/2}.
\de
On the other hand, (\ref{barf1checkf1a}) yields that
\ce
&&\left(\mE\left|\int_{k\d}^{(k+1)\d}(\bar f_1(\frac{r}{\e},\check X_{k\d})-\check{f}_{1,a}(\check X_{k\d}))\dif r\right|^2_\mH\right)^{1/2}\\
&=&\d\left(\mE\left|\frac{\e}{\d}\int_{\frac{k\d}{\e}}^{\frac{(k+1)\d}{\e}}\bar f_1(v,\check X_{k\d})\dif v-\check{f}_{1,a}(\check X_{k\d})\right|^2_\mH\right)^{1/2}\\
&\leq&C\d(1+|x_0|^2_\mH+|y_0|^2_\mH)^{1/2}\kappa^*(\frac{\d}{\e}).
\de
So, we infer that
\ce
J_{31}\leq C_R\kappa^*(\frac{\d}{\e})(1+|x_0|^2_\mH+|y_0|^2_\mH).
\de
For $J_{32}$, by the H\"older inequality it holds that
\ce
J_{32}&\leq& C_R\mE\sup\limits_{t\in[0,T\wedge\bar\tau_2]}\int_{[\frac{t}{\d}]\d}^{t}\left|\bar Z^\e_{r(\d)}\right|_\mH\left|\bar f_1(\frac{r}{\e},\check X_{r(\d)})-\check{f}_{1,a}(\check X_{r(\d)})\right|_\mH\dif r\\
&\leq& C_R\mE\sup\limits_{t\in[0,T\wedge\bar\tau_2]}\left(\int_{[\frac{t}{\d}]\d}^{t}\left|\bar Z^\e_{r(\d)}\right|_\mH^2\dif r\right)^{1/2}\left(\int_{[\frac{t}{\d}]\d}^{t}\left|\bar f_1(\frac{r}{\e},\check X_{r(\d)})-\check{f}_{1,a}(\check X_{r(\d)})\right|^2_\mH\dif r\right)^{1/2}\\
&\leq& C_R\d\mE\left(\sup\limits_{t\in[0,T]}|\bar Z^\e_t|^2_\mH\right)^{1/2}\left(1+\sup\limits_{t\in[0,T]}|\check X_t|^2_\mH\right)^{1/2}\\
&\leq& C_R\d(1+|x_0|^2_\mH).
\de
Then we establish (\ref{j3}). By the similar deduction to that for (\ref{j3}), we can obtain (\ref{j4}).
\end{proof}

\br\label{l2normaverprin}
Comparing Proposition \ref{barxecheckxespr} with Theorem 3.4 in \cite{Q2}, we observe that the convergence of $\bar X^\e$ to $\check X$ established here is in the mean square sense, whereas the corresponding result in \cite{Q2} only proves convergence in probability. This demonstrates that our result is strictly stronger.
\er

At present, we are ready to prove Theorem \ref{asymstroapth}.

{\bf Proof of Theorem \ref{asymstroapth}.}  For any $\eta>0$, by (\ref{barxecheckxesti}) and the Chebyshev inequality, it holds that
\ce
&&\mP\left\{\sup\limits_{t\in[0,T]}|\bar X^\e_t-\check X_t|_\mH\geq\frac{\eta}{2}\right\}\\
&=&\mP\left\{\sup\limits_{t\in[0,T]}|\bar X^\e_t-\check X_t|_\mH\geq\frac{\eta}{2}, T\leq\check\tau_2\right\}+\mP\left\{\sup\limits_{t\in[0,T]}|\bar X^\e_t-\check X_t|_\mH\geq\frac{\eta}{2}, T>\check\tau_2\right\}\\
&\leq&\mP\left\{\sup\limits_{t\in[0,T\wedge\check\tau_2]}|\bar X^\e_t-\check X_t|_\mH\geq\frac{\eta}{2}\right\}+\mP\left\{T>\check\tau_2\right\}\\
&\leq&\frac{4}{\eta^2}\mE\sup\limits_{t\in[0,T\wedge\check\tau_2]}|\bar X^\e_t-\check X_t|^2_\mH+\frac{1}{R}\mE\left[|\bar X^\e_T|_\mH+|\check X_T|_\mH+\int_0^T\(\|\bar X^\e_{r}\|_\mV^2+\|\check{X}_{r}\|^2_\mV\)\dif r\right]\\
&\leq&\frac{4C_R}{\eta^2}\left(\d+\d^{1/2}+\d^{1/4}+\kappa^*(\frac{\d}{\e})+\kappa_2(\frac{\d}{\e})\right)(1+|x_0|^2_\mH+|y_0|^2_\mH)+\frac{C}{R}\left(1+|x_0|_\mH^2\right).
\de
Letting $\e\rightarrow0$ first, $\d\rightarrow0$ and $R\rightarrow\infty$, we obtain that
\ce
\lim\limits_{\e\rightarrow0}\mP\left\{\sup\limits_{t\in[0,T]}|\bar X^\e_t-\check X_t|_\mH\geq\frac{\eta}{2}\right\}=0.
\de
Besides, by Theorem \ref{genestroapth}, it holds that
\ce
\lim\limits_{\e\rightarrow0}\mP\left\{\sup\limits_{t\in[0,T]}|X^\e_t-\bar X^\e_t|_\mH\geq\frac{\eta}{2}\right\}=0.
\de
Note that
\ce
\mP\left\{\sup\limits_{t\in[0,T]}|X^\e_t-\check X_t|_\mH\geq\eta\right\}\leq \mP\left\{\sup\limits_{t\in[0,T]}|X^\e_t-\bar X^\e_t|_\mH\geq\frac{\eta}{2}\right\}+\mP\left\{\sup\limits_{t\in[0,T]}|\bar X^\e_t-\check X_t|_\mH\geq\frac{\eta}{2}\right\}.
\de
Thus, we conclude that
\ce
\lim\limits_{\e\rightarrow0}\mP\left\{\sup\limits_{t\in[0,T]}|X^\e_t-\check X_t|_\mH\geq\eta\right\}=0.
\de

\section{Proof of Theorem \ref{peristroapth}}\label{peristroapthproo}

In this section, we prove Theorem \ref{peristroapth}. First, we investigate that the evolution system $\{\mu_t^x\}_{t\in\mR}$ of measures is $\varrho_2$-periodic under $(\mathbf{H}^1_{f_{2}, g_{2}})$, $(\mathbf{H}^2_{f_{2}, g_{2}})$ and $(\mathbf{H}^p_{f_{2}, g_{2}})$. Then based on this, $(\mathbf{H}^p_{f_{1}})$ and $(\mathbf{H}^a_{g_{1}})$, some important estimates are presented. Finally, we prove Theorem \ref{peristroapth}.

First of all, we state the following result, which is from \cite[Theorem 4.1]{DaT}.

\bl\label{muxperi}
Assume that $(\mathbf{H}^1_{f_{2}, g_{2}})$, $(\mathbf{H}^2_{f_{2}, g_{2}})$ and $(\mathbf{H}^p_{f_{2}, g_{2}})$ hold. Then the evolution system $\{\mu_t^x\}_{t\in\mR}$ of measures in Lemma \ref{mux} is $\varrho_2$-periodic, 
i.e. for any $t\in\mR$, $\mu_{t+\varrho_2}^x=\mu_t^x$.
\el

By the above lemma and $(\mathbf{H}^p_{f_{1}})$, picking a constant $\varrho>0$ satisfying $\varrho=n_2\varrho_1=n_1\varrho_2$ for some $n_1, n_2\in \mN_+=\{1,2,\cdots\}$, we have that for any $t\geq 0$ and $x\in\mH$
\ce
\bar f_1(t+\varrho, x)=\int_\mH f_1(t+\varrho, x,y)\mu_{t+\varrho}^x(\dif y)=\int_\mH f_1(t, x,y)\mu_t^x(\dif y)=\bar f_1(t,x),
\de
that is, $\bar f_1$ is $\varrho$-periodic. Moreover, it holds that for any $T>0$ and $t\geq0$
\ce
&&\left|\frac{1}{T}\int_t^{t+T}\bar f_1(r,x)\dif r-\frac{1}{\varrho}\int_0^\varrho\bar f_1(r,x)\dif r\right|_\mH\\
&\leq&\left|\frac{1}{T}\int_t^{t+[\frac{T}{\varrho}]\varrho}\bar f_1(r,x)\dif r+\frac{1}{T}\int_{t+[\frac{T}{\varrho}]\varrho}^{t+T}\bar f_1(r,x)\dif r-\frac{1}{\varrho}\int_0^\varrho\bar f_1(r,x)\dif r\right|_\mH\\
&\leq&\left|\frac{1}{T}\int_t^{t+[\frac{T}{\varrho}]\varrho}\bar f_1(r,x)\dif r-\frac{1}{\varrho}\int_0^\varrho\bar f_1(r,x)\dif r\right|_\mH+\left|\frac{1}{T}\int_{t+[\frac{T}{\varrho}]\varrho}^{t+T}\bar f_1(r,x)\dif r\right|_\mH\\
&\leq&\left|\frac{[\frac{T}{\varrho}]}{T}\int_0^\varrho\bar f_1(r,x)\dif r-\frac{1}{\varrho}\int_0^\varrho\bar f_1(r,x)\dif r\right|_\mH+\frac{C\varrho(1+|x|_\mH)}{T}\\
&\leq&\frac{2C\varrho(1+|x|_\mH)}{T},
\de
So, we conclude that
\ce
\sup\limits_{t\geq 0}\left|\frac{1}{T}\int_t^{t+T}\bar f_1(r,x)\dif r-\frac{1}{\varrho}\int_0^\varrho\bar f_1(r,x)\dif r\right|_\mH\leq \frac{C(1+|x|_\mH)}{T},
\de
and further noticing $\bar f_{1,p}(x)=\frac{1}{\varrho}\int_0^\varrho \bar f_1(t,x)\dif t$
\be
\sup\limits_{t\geq 0}\left|\frac{1}{T}\int_t^{t+T}\bar f_1(r,x)\dif r-\bar f_{1,p}(x)\right|_\mH\leq \frac{C(1+|x|_\mH)}{T}.
\label{barf1barf1pcond}
\ee

\bl
Assume that $(\mathbf{H}_{f_{1}, g_{1}})$, $(\mathbf{H}^1_{f_{2}, g_{2}})$, $(\mathbf{H}^2_{f_{2}, g_{2}})$, $(\mathbf{H}^p_{f_{1}})$, $(\mathbf{H}^a_{g_{1}})$ and $(\mathbf{H}^p_{f_{2}, g_{2}})$ hold. Then Eq.(\ref{averequap}) has a unique solution $(\bar{X}, \bar{K})$ such that 
\be
\mE\left[\sup _{t \in[0, T]}\left|\bar X_t\right|_\mH^2+\int_0^T\left\|\bar X_t\right\|_\mV^2\dif t\right]\leq C\left(1+|x_0|_\mH^2\right).
\label{barxesti}
\ee
Moreover, set the stopping time $\bar\tau_1:=\inf\{t>0, |\bar X_t|_\mH>R\}$, and for any $\d>0$ small enough,
\be
\mE\int_0^{T\wedge\bar\tau_1}|\bar X_t-\bar X_{t(\d)}|_\mH^2\dif t\leq C\d^{1/2}\left(1+|x_0|_\mH^2\right).
\label{barxincr}
\ee
\el

Since the proof of the above lemma is similar to that for Lemma \ref{checkxwellpose}, we omit it.

Next, we investigate the difference between $\bar X^\e$ and $\bar X$. In order to do this, we define the stopping time
\ce
\bar\tau_2:=\inf\left\{t>0, |\bar X^\e_t|_\mH+|\bar X_t|_\mH+\int_0^t\(\|\bar X^\e_{r}\|_\mV^2+\|\bar{X}_{r}\|^2_\mV\)\dif r>R\right\}.
\de
Then $\bar\tau_2\leq\tau_2\wedge\bar\tau_1$. By the similar deduction to that in the proof of Proposition \ref{barxecheckxespr}, we conclude the following result.

\bp\label{barxebarxestipro}
Suppose that $(\mathbf{H}_{f_{1}, g_{1}})$, $(\mathbf{H}^1_{f_{2}, g_{2}})$, $(\mathbf{H}^2_{f_{2}, g_{2}})$, $(\mathbf{H}^p_{f_{1}})$, $(\mathbf{H}^a_{g_{1}})$ and $(\mathbf{H}^p_{f_{2}, g_{2}})$ hold. Then for any $T>0$
\be
\mE\sup\limits_{t\in[0,T\wedge\bar\tau_2]}|\bar X^\e_t-\bar X_t|^2_\mH\leq C_R\left(\d+\d^{1/2}+\d^{1/4}+\frac{\e}{\d}+\kappa_2(\frac{\d}{\e})\right)(1+|x_0|^2_\mH).
\label{barxebarxesti}
\ee
\ep

At present, we are ready to prove Theorem \ref{peristroapth}.

{\bf Proof of Theorem \ref{peristroapth}.} Following the line in the proof of Theorem \ref{asymstroapth}, we derive Theorem \ref{peristroapth}.

\section{Proof of Theorem \ref{aperistroapth}}\label{aperistroapthproo}

In this section, we prove Theorem \ref{aperistroapth}. First, we establish that the averaged coefficient $\bar f_1(t,x)$ is uniformly almost periodic. Then $\mE\sup\limits_{t\in[0,T]}|\bar X^\e_t-\hat X_t|^2_\mH$ is described. Finally we give the proof of Theorem \ref{aperistroapth}.

First of all, we prepare the following result.

\bl\label{muxunifalmoperi}
Assume that $(\mathbf{H}^1_{f_{2}, g_{2}})$, $(\mathbf{H}^2_{f_{2}, g_{2}})$, $({\bf H}^{ex}_{f_2, g_2})$ and $(\mathbf{H}^{ap}_{f_{2}, g_{2}})$ hold. Then $\mR\times\mH\ni (t,x)\to \mu_t^x\in\sP(\mH)$ from Lemma \ref{mux}  is uniformly almost periodic.
\el
\begin{proof}
First of all, we prove that for any $x\in\mH$, $\mR\ni t\rightarrow \mu_t^x\in\sP(\mH)$ is almost periodic. Indeed, by the definition of $A$, we know that $A$ is continuous. $(\mathbf{H}^{ap}_{f_{2}, g_{2}})$ assures that $f_2(t, x, y)$ is continuous. Then $(H1)$ in \cite[Theorem 4.3]{ChLi0} is satisfied. By $(\mathbf{H}^1_{f_{2}, g_{2}})$, $(\mathbf{H}^2_{f_{2}, g_{2}})$, (\ref{discond1}) and (\ref{discond2}), $(H2')$, $(H3)$ and $(HL)$ in \cite[Theorem 4.3]{ChLi0} hold. Based on the definition of $A$ and (\ref{f2g2linegrow}), $(H4)$ in \cite[Theorem 4.3]{ChLi0} is justified. And $({\bf H}^{ex}_{f_2, g_2})$ implies that $(H5)$ in \cite[Theorem 4.3]{ChLi0} is verified. By \cite[Theorem 4.3]{ChLi0}, we conclude that for any $x\in\mH$, $\mR\ni t\rightarrow \mu_t^x\in\sP(\mH)$ is almost periodic. 

Next, we prove that $\mR\times\mH\ni (t,x)\to \mu_t^x\in\sP(\mH)$ is uniformly almost periodic. Note that for any $x_1, x_2\in\mH$ and $t\in\mR$
\ce
d_{BL}(\mu_t^{x_1}, \mu_t^{x_2})&=&\sup\limits_{\|\phi\|_{BL} \leq 1}\left|\int_\mH \phi \dif (\mu_t^{x_1}-\mu_t^{x_2})\right|\\
&=&\sup\limits_{\|\phi\|_{BL} \leq 1}\left|\mE\phi(Y_t^{x_1})-\mE\phi(Y_t^{x_2})\right|\\
&\leq&(\mE|Y_t^{x_1}-Y_t^{x_2}|^2_\mH)^{1/2}\\
&\leq&C|x_1-x_2|_\mH,
\de
where we use (\ref{yxdiffesti}) and $C>0$ is independent of $t$. Thus, it holds that
\ce
\sup\limits_{t\in\mR}d_{BL}(\mu_t^{x_1}, \mu_t^{x_2})\leq C|x_1-x_2|_\mH,
\de
and this implies that $\mu_t^x$ is continuous in $x$ uniformly in $t$. So, by \cite[Theorem 2.10]{aF}, $\mR\times\mH\ni (t,x)\to \mu_t^x\in\sP(\mH)$ is uniformly almost periodic.
\end{proof}

\br
We note that the similar result is established in \cite[Lemma 4.3]{CSX}, where Theorem \ref{uapchartheo} is employed. In contrast, our proof achieves the same conclusion with greater concision.
\er

Next, keeping the above lemma in mind, we prove that $\mR_+\times\mH\ni(t,x)\to \bar{f}_1(t,x)=\int_\mH f_1(t,x,y)\mu^x_t(\dif y)\in\mH$ is uniformly almost periodic.

\bl\label{barf1unifalmoperi}
Suppose that $(\mathbf{H}_{f_{1}, g_{1}})$, $(\mathbf{H}^1_{f_{2}, g_{2}})$, $(\mathbf{H}^2_{f_{2}, g_{2}})$, $(\mathbf{H}^{ap}_{f_{1}})$, $({\bf H}^{ex}_{f_2, g_2})$ and $(\mathbf{H}^{ap}_{f_{2}, g_{2}})$ hold. Then $\mR_+\times\mH\ni(t,x)\to \bar{f}_1(t,x)=\int_\mH f_1(t,x,y)\mu^x_t(\dif y)\in\mH$ is uniformly almost periodic.
\el
\begin{proof}
For any $m\in\mN$, set 
\ce
f_{1,m}(t,x,y)=
\left\{\begin{array}{ll}
f_1(t,x,y), \qquad |y|_\mH\leq m,\\
f_1(t,x,y\frac{m}{|y|_\mH}), \quad |y|_\mH>m,
\end{array}
\right.
\de
and $f_{1,m}(t,x,\cdot)$ is Lipschitz continuous and bounded. So, it holds that
\ce
\int_\mH f_1(t,x,y)\mu^x_t(\dif y)=\int_\mH f_{1,m}(t,x,y)\mu^x_t(\dif y)+\int_{|y|_\mH>m} (f_1(t,x,y)-f_{1,m}(t,x,y))\mu^x_t(\dif y).
\de
Besides, (\ref{f1g1linegrow}) and (\ref{muxboun}) imply that
\ce
&&\sup\limits_{t\in\mR_+}\left|\int_{|y|_\mH>m} (f_1(t,x,y)-f_{1,m}(t,x,y))\mu^x_t(\dif y)\right|_\mH\\
&\leq&2L_1\sup\limits_{t\in\mR_+}\int_{|y|_\mH>m} (1+|x|_\mH+|y|_\mH)\mu^x_t(\dif y)\\
&\leq&\frac{2L_1}{m}\sup\limits_{t\in\mR_+}\int_\mH[(1+|x|_\mH)|y|_\mH+|y|^2_\mH]\mu^x_t(\dif y)\\
&\leq&\frac{C}{m}(1+|x|^2_\mH).
\de
Then for any $\g>0, R>0$, there exists a $m_0\in\mN$ such that
\ce
\sup\limits_{t\in\mR_+, |x|_\mH\leq R}\left|\int_{|y|_\mH>m_0} (f_1(t,x,y)-f_{1,m_0}(t,x,y))\mu^x_t(\dif y)\right|_\mH<\g.
\de
Combining the above deduction, we obtain that for any $t\in\mR_+, \tau\in\mR_+$ and $x\in\mH$ with $|x|_\mH\leq R$,
\ce
&&\left|\int_{\mH} f_1(t+\tau,x,y)\mu^x_{t+\tau}(\dif y)-\int_{\mH} f_1(t,x,y)\mu^x_t(\dif y)\right|_\mH\\
&\leq&\left|\int_{\mH} f_1(t+\tau,x,y)\mu^x_{t+\tau}(\dif y)-\int_{\mH} f_1(t,x,y)\mu^x_{t+\tau}(\dif y)\right|_\mH\\
&&+\left|\int_{\mH} f_1(t,x,y)\mu^x_{t+\tau}(\dif y)-\int_{\mH} f_1(t,x,y)\mu^x_t(\dif y)\right|_\mH\\
&\leq&\int_{\mH}|f_1(t+\tau,x,y)-f_1(t,x,y)|_\mH\mu^x_{t+\tau}(\dif y)\\
&&+\left|\int_{\mH} f_{1,m_0}(t,x,y)\mu^x_{t+\tau}(\dif y)-\int_\mH f_{1,m_0}(t,x,y)\mu^x_t(\dif y)\right|_\mH+2\g\\
&=:&I_1+I_2+2\g.
\de

For $I_1$, $(\mathbf{H}^{ap}_{f_{1}})$, (\ref{f1g1linegrow}) and (\ref{muxboun}) imply that
\ce
I_1&=&\int_{|y|_\mH\leq R'} |f_1(t+\tau,x,y)-f_1(t,x,y)|_\mH\mu^x_{t+\tau}(\dif y)\\
&&+\int_{|y|_\mH>R'} |f_1(t+\tau,x,y)-f_1(t,x,y)|_\mH\mu^x_{t+\tau}(\dif y)\\
&\leq&\sup\limits_{t\in\mR_+, |x|_\mH\leq R, |y|_\mH\leq R'}|f_1(t+\tau,x,y)-f_1(t,x,y)|_\mH\\
&&+2L_1\int_{|y|_\mH>R'} (1+|x|_\mH+|y|_\mH)\mu^x_{t+\tau}(\dif y)\\
&\leq&\g+\frac{C}{R'}(1+R^2),
\de
where $R'>0$ is large enough. For $I_2$, from $(\mathbf{H}_{f_{1}, g_{1}})$ and (\ref{f1g1linegrow}), it follows that
\ce
I_2\leq \|f_{1,m_0}(t,x,\cdot)\|_{BL}d_{BL}(\mu^x_{t+\tau},\mu^x_t)\leq C_{m_0, R}d_{BL}(\mu^x_{t+\tau},\mu^x_t),
\de
which together with Lemma \ref{muxunifalmoperi} yields that 
\ce
I_2\leq C_{m_0, R}\g.
\de

Finally, collecting the above deduction and Definition \ref{apuapdefi}, we conclude that $\mR_+\times\mH\ni(t,x)\to \bar{f}_1(t,x)=\int_\mH f_1(t,x,y)\mu^x_t(\dif y)\in\mH$ is uniformly almost periodic.
\end{proof}

Based on the above lemma and Theorem \ref{uapmeantheo}, we have that there exists $\hat f_{1,ap}: \mH\rightarrow\mH$ such that for any $T>0$, $R>0$ and $x\in\mH$ with $|x|_\mH\leq R$
\be
\sup\limits_{t\geq 0}\left|\frac{1}{T}\int_t^{t+T}\bar f_1(s,x)\dif s-\hat f_{1,ap}(x)\right|_\mH\leq \hat\kappa_R(T),
\label{hatkapppro}
\ee
where $\hat\kappa_R: [0,\infty)\to [0,\infty)$ is bounded and satisfies
\ce
\lim\limits_{T\rightarrow\infty}\hat\kappa_R(T)=0.
\de

\bl
Assume that $(\mathbf{H}_{f_{1}, g_{1}})$, $(\mathbf{H}^1_{f_{2}, g_{2}})$, $(\mathbf{H}^2_{f_{2}, g_{2}})$, $(\mathbf{H}^{ap}_{f_{1}})$, $(\mathbf{H}^a_{g_{1}})$, $({\bf H}^{ex}_{f_2, g_2})$ and $(\mathbf{H}^{ap}_{f_{2}, g_{2}})$ hold. Then Eq.(\ref{averequaap}) has a unique solution $(\hat{X}, \hat{K})$ such that 
\be
\mE\left[\sup _{t \in[0, T]}\left|\hat X_t\right|_\mH^2+\int_0^T\left\|\hat X_t\right\|_\mV^2\dif t\right]\leq C\left(1+|x_0|_\mH^2\right).
\label{hatxesti}
\ee
Moreover, set the stopping time $\hat\tau_1:=\inf\{t>0, |\hat X_t|_\mH>R\}$, and for any $\d>0$ small enough,
\be
\mE\int_0^{T\wedge\hat\tau_1}|\hat X_t-\hat X_{t(\d)}|_\mH^2\dif t\leq C\d^{1/2}\left(1+|x_0|_\mH^2\right).
\label{hatxincr}
\ee
\el

The proof of the above lemma follows closely along the lines of that for Lemma \ref{checkxwellpose}, and is therefore omitted.

Next, we observe the difference between $\bar X^\e$ and $\hat X$. To do this, we define the stopping time
\ce
\hat\tau_2:=\inf\left\{t>0, |\bar X^\e_t|_\mH+|\hat X_t|_\mH+\int_0^t\(\|\bar X^\e_{r}\|_\mV^2+\|\hat{X}_{r}\|^2_\mV\)\dif r>R\right\}.
\de
Then $\hat\tau_2\leq\tau_2\wedge\hat\tau_1$.

\bp\label{barxehatxestipro}
Suppose that $(\mathbf{H}_{f_{1}, g_{1}})$, $(\mathbf{H}^1_{f_{2}, g_{2}})$, $(\mathbf{H}^2_{f_{2}, g_{2}})$, $(\mathbf{H}^{ap}_{f_{1}})$, $(\mathbf{H}^a_{g_{1}})$, $({\bf H}^{ex}_{f_2, g_2})$ and $(\mathbf{H}^{ap}_{f_{2}, g_{2}})$ hold. Then for any $T>0$
\be
\mE\sup\limits_{t\in[0,T\wedge\hat\tau_2]}|\bar X^\e_t-\hat X_t|^2_\mH\leq C_R\left(\d+\d^{1/2}+\d^{1/4}+\hat\kappa_R(\frac{\d}{\e})+\kappa_2(\frac{\d}{\e})\right)(1+|x_0|^2_\mH).
\label{barxehatxesti}
\ee
\ep
\begin{proof}
Since most of the proof are the same as that for Proposition \ref{barxecheckxespr}, we will only present the differences, that is, the deduction for $J_{31}$. Now we derive as follows:
\ce
&&C_R\mE\sup\limits_{t\in[0,T\wedge\bar\tau_2 ]}\left|\int_{0}^{[\frac{t}{\d}]\d}\<\bar X^\e_{r(\d)}-\hat X_{r(\d)},\bar f_1(\frac{r}{\e},\hat X_{r(\d)})-\hat{f}_{1,ap}(\hat X_{r(\d)})\>_\mH\dif r\right|\\
&=&C_R\mE\sup\limits_{t\in[0,T\wedge\bar\tau_2]}\left|\sum\limits_{k=0}^{[\frac{t}{\d}]-1}\int_{k\d}^{(k+1)\d}\<\bar X^\e_{k\d}-\hat X_{k\d},\bar f_1(\frac{r}{\e},\hat X_{k\d})-\hat{f}_{1,ap}(\hat X_{k\d})\>_\mH\dif r\right|\\
&\leq&C_R\mE\sup\limits_{t\in[0,T\wedge\hat\tau_2]}\sum\limits_{k=0}^{[\frac{t}{\d}]-1}\left|\bar X^\e_{k\d}-\hat X_{k\d}\right|_\mH\left|\int_{k\d}^{(k+1)\d}(\bar f_1(\frac{r}{\e},\hat X_{k\d})-\hat{f}_{1,ap}(\hat X_{k\d}))\dif r\right|_\mH\\
&\leq&C_R\mE\sum\limits_{k=0}^{[\frac{T\wedge\hat\tau_2}{\d}]-1}\left|\bar X^\e_{k\d}-\hat X_{k\d}\right|_\mH\left|\int_{k\d}^{(k+1)\d}(\bar f_1(\frac{r}{\e},\hat X_{k\d})-\hat{f}_{1,ap}(\hat X_{k\d}))\dif r\right|_\mH \\
&\leq&C_R\frac{T}{\d}\sup\limits_{0\leq k\leq [\frac{T}{\d}]-1}\mE\left|\bar X^\e_{k\d}-\hat X_{k\d}\right|_\mH\left|\int_{k\d}^{(k+1)\d}(\bar f_1(\frac{r}{\e},\hat X_{k\d})-\hat{f}_{1,ap}(\hat X_{k\d}))\dif r\right|_\mH I_{\{(k+1)\d\leq \hat\tau_2\}}.
\de
Then (\ref{hatkapppro}) yields that for $k=0, 1, \cdots, [\frac{T}{\d}]-1$
\ce
&&\left|\int_{k\d}^{(k+1)\d}(\bar f_1(\frac{r}{\e},\hat X_{k\d})-\hat{f}_{1,ap}(\hat X_{k\d}))\dif r\right|_\mH I_{\{(k+1)\d\leq \hat\tau_2\}}\\
&=&\d\left|\frac{\e}{\d}\int_{\frac{k\d}{\e}}^{\frac{(k+1)\d}{\e}}\bar f_1(v,\hat X_{k\d})\dif v-\hat{f}_{1,ap}(\hat X_{k\d})\right|_\mH I_{\{(k+1)\d\leq \hat\tau_2\}}\\
&\leq&\d\hat\kappa_R(\frac{\d}{\e}).
\de
And by (\ref{barxee}), (\ref{hatxesti}) and the H\"older inequality it holds that
\ce
\mE\left|\bar X^\e_{k\d}-\hat X_{k\d}\right|_\mH\leq \sqrt 2\left(\mE\left(|\bar X^\e_{k\d}|^2_\mH+|\hat X_{k\d}|^2_\mH\right)\right)^{1/2}\leq C(1+|x_0|_\mH).
\de
So, we obtain that
\ce
&&C_R\mE\sup\limits_{t\in[0,T\wedge\hat\tau_2]}\sum\limits_{k=0}^{[\frac{t}{\d}]-1}\left|\bar X^\e_{k\d}-\hat X_{k\d}\right|_\mH\left|\int_{k\d}^{(k+1)\d}(\bar f_1(\frac{r}{\e},\hat X_{k\d})-\hat{f}_{1,ap}(\hat X_{k\d}))\dif r\right|_\mH\\
&\leq& C_R\hat\kappa_R(\frac{\d}{\e})(1+|x_0|_\mH).
\de
The proof is complete.
\end{proof}

With these results, we proceed to the proof of Theorem \ref{aperistroapth}.

{\bf Proof of Theorem \ref{aperistroapth}.} Theorem \ref{aperistroapth} follows from the same argument as Theorem \ref{asymstroapth}.

\section{An example}\label{exam}

In this section, we provide an example to illustrate our results.

\bx
Take $\mU_1=\mU_2=\mR$ and consider the following nonautonomous multiscale stochastic Burgers equations with reflection on the spatial domain $[0,1]$:
\be\left\{\begin{array}{ll}
\frac{\p X_t^\e(\xi)}{\p t}=\frac{\p^2 X_t^\e(\xi)}{\p\xi^2}+\frac{1}{2}\frac{\p(X_t^\e)^2(\xi)}{\p \xi}+\phi_1(\frac{t}{\e})f(X_t^\e,Y_t^\e)(\xi)+(1-e^{-\frac{t}{2\e}})X_t^\e(\xi)\frac{\p W^1(t,\xi)}{\p t}\\
\qquad\qquad\qquad +\frac{\p K^\e(t,\xi)}{\p t},\\
\frac{\p Y_t^\e(\xi)}{\p t}=\frac{1}{\e}\left[\frac{\p^2 Y_t^\e(\xi)}{\p\xi^2}+X_t^\e(\xi)+\phi_2(\frac{t}{\e})Y_t^\e(\xi)\right]+\frac{1}{\sqrt \e}Y_t^\e(\xi)\frac{\p W^2(t,\xi)}{\p t},\\
X_t^\e(\xi)\geq 0, \\
X_t^\e(0)=X_t^\e(1)=Y_t^\e(0)=Y_t^\e(1)=0, \quad t\geq 0,\\
X_0^\e(\xi)=x_0(\xi)\geq 0, \quad Y_0^\e(\xi)=y_0(\xi),
\end{array}
\right.
\label{exeq}
\ee
where $\phi_1: \mR_+\to \mR, \phi_2: \mR\to \mR$ satisfy $|\phi_1(t)|\leq l_1, |\phi_2(t)|\leq l_2$ for $l_1, l_2>0$ and $\lambda_1-\frac{3}{2}l_2-1>0$, and $f: \mH\times\mH\to \mH$ is Lipschitz continuous with $|f(0,0)|_\mH\leq\|f\|_{Lip}$. Note that 
\ce
&&f_1(t,x,y)=\phi_1(t)f(x,y), \quad g_1(t,x)=(1-e^{-\frac{t}{2}})x, \quad t\in\mR_+, \quad x,y\in\mH,\\
&&f_2(t,x,y)=x+\phi_2(t)y, \qquad g_2(t,x,y)=y, \quad t\in\mR, \quad x,y\in\mH.
\de
On the one hand, it is easy to justify that for any $t\in\mR_+, x_1, x_2, y_1, y_2\in\mH$,
\ce
&&|f_1(t,x_1,y_1)-f_1(t,x_2,y_2)|_\mH\leq l_1\|f\|_{Lip}(|x_1-x_2|_\mH+|y_1-y_2|_\mH),\\
&&|f_1(t,0,0)|_\mH\leq l_1\|f\|_{Lip},\\
&&|g_1(t,x_1)-g_1(t,x_2)|_\mH\leq |x_1-x_2|_\mH,\\
&&|g_1(t,0)|_\mH=0\leq 1.
\de
$({\bf H}_{f_1, g_1})$ holds. On the other hand, for any $t\in\mR, x_1, x_2, y_1, y_2\in\mH$,
\ce
&&|f_2(t,x_1,y_1)-f_2(t,x_2,y_2)|_\mH\leq |x_1-x_2|_\mH+l_2|y_1-y_2|_\mH, \\
&&|f_2(t,0,0)|_\mH\leq 1,\\
&&|g_2(t,x_1,y_1)-g_2(t,x_2,y_2)|_\mH\leq |y_1-y_2|_\mH, \\
&&|g_2(t,0,0)|_\mH\leq 1,\\
&&\lambda_1-\frac{3}{2}l_2-1>0.
\de
In terms of the above deduction, we conclude that $({\bf H}^1_{f_2, g_2})$ and $({\bf H}^2_{f_2, g_2})$ hold. So, from Theorem \ref{genestroapth}, it follows that
\ce
\lim\limits_{\e\rightarrow0}\mP\left(\sup\limits_{t\in[0,T]}|X^\e_t-\bar X^\e_t|_\mH\geq \eta\right)=0,
\de
where $\bar X^\e$ is the solution of the following SPDE
\ce\left\{\begin{array}{ll}
\dif \bar X_t^\e=[A\bar X_t^\e+B(\bar X_t^\e,\bar X_t^\e)+\phi_1(\frac{t}{\e})\int_\mH f(\bar X_t^\e,y)\mu^{\bar X_t^\e}_\frac{t}{\e}(\dif y)]\dif t+(1-e^{-\frac{t}{2\e}})\bar X_t^\e\dif W^1_t+\dif \bar K^\e_t,\\
\bar X_t^\e(\xi)\geq 0, \quad \xi\in[0,1],\\
\bar X_t^\e(0)=\bar X_t^\e(1)=0,\quad t\geq 0,\\
\bar X_0^\e=x_0,
\end{array}
\right.
\de
and $\{\mu_t^x\}_{t\in\mR}$ is the evolution system of measures for the following SPDE: for any $s\in\mR$
\ce\left\{\begin{array}{ll}
\dif Y_t^{x,y_0}=[AY_t^{x,y_0}+x+\phi_2(t)Y_t^{x,y_0}]\dif t+Y_t^{x,y_0}\dif \tilde W^2_t,\\
Y_t^{x,y_0}(0)=Y_t^{x,y_0}(1)=0, \quad t\geq s,\\
Y_s^{x,y_0}=y_0.
\end{array}
\right.
\de

{\bf The asymptotic case}  $(i)$ Assume that there exists a constant $a_1\in\mR, a_1\neq0$ such that 
$$
\lim\limits_{T\rightarrow\infty}\sup\limits_{t\geq 0}\left|\frac{1}{T}\int_t^{t+T}\phi_1(s)\dif s-a_1\right|=0.
$$
For example, $\phi_1(t)=\frac{t}{1+t}, a_1=1$ and 
\ce
\sup\limits_{t\geq 0}\left|\frac{1}{T}\int_t^{t+T}\frac{s}{1+s}\dif s-1\right|\leq \frac{1}{\sqrt T}.
\de
So, $\check f_1(x,y)=a_1 f(x,y)$ and 
\ce
\kappa_1(T)=\sup\limits_{t\geq 0}\left|\frac{1}{T}\int_t^{t+T}\phi_1(s)\dif s-a_1\right|.
\de
$({\bf H}^a_{f_1})$ is verified. 

$(ii)$ Besides, $\bar g_1(x)=x$ and 
\ce
\sup\limits_{t\geq 0}\frac{1}{T}\int_t^{t+T}|(1-e^{-\frac{s}{2}})x-x|^2_{\mH}\dif s\leq\frac{1}{T}(1+|x|^2_\mH).
\de
That is, $({\bf H}^a_{g_1})$ holds.

$(iii)$ Suppose that there exists a constant $a_2\in\mR, a_2\neq0$ such that $\lim\limits_{t\to\infty}\phi_2(t)=a_2$ and 
$$
\lim\limits_{T\rightarrow\infty}\sup\limits_{t\geq 0}\frac{1}{T}\int_t^{t+T}(\phi_2(s)-a_2)^2\dif s=0.
$$
For example, $\phi_2(t)=1-e^{-t}, a_2=1$ and 
$$
\lim\limits_{T\rightarrow\infty}\sup\limits_{t\geq 0}\frac{1}{T}\int_t^{t+T}e^{-2s}\dif s=0.
$$
Then, $ \bar f_2(x,y)=x+a_2 y, \bar g_2(x,y)=y$, $\kappa_3(t)=|\phi_2(t)-a_2|$ and $({\bf H}^a_{f_2, g_2})$ is satisfied. 

So, Theorem \ref{asymstroapth} yields that
\ce
\lim\limits_{\e\rightarrow 0}\mP\left(\sup\limits_{t\in[0,T]}|X_t^\e-\check X_t|_\mH\geq\eta\right)=0,
\de
where $\check X$ is the solution of the following averaged SPDE
\ce\left\{\begin{array}{ll}
\dif \check X_t=[A\check X_t+B(\check X_t,\check X_t)+a_1\int_\mH f(\check X_t,y)\mu^{\check X_t}(\dif y)]\dif t+\check X_t\dif W^1_t+\dif \check K_t,\\
\check X_t(\xi)\geq 0, \quad \xi\in[0,1],\\
\check X_t(0)=\check X_t(1)=0,\quad t\geq 0,\\
\check X_0=x_0,
\end{array}
\right.
\de
and $\mu^x$ is the unique invariant measure of the following SPDE
\ce\left\{\begin{array}{ll}
\dif \bar{Y}_t=\left[A\bar{Y}_t+x+a_2\bar{Y}_t\right]\dif t+\bar{Y}_t\dif W^2_t,\\
\bar{Y}_t(0)=\bar{Y}_t(1)=0,  \quad t\geq 0,\\
\bar{Y}_0=y_0.
\end{array}
\right.
\de

{\bf The periodic case.} Assume that there is a constant $\varrho>0$ such that $\phi_1(t+\varrho)=\phi_1(t), \phi_2(t+\varrho)=\phi_2(t)$ for any $t\geq 0$. For example, $\phi_1(t)=\phi_2(t)=\sin x, \varrho=2\pi$. Then $f_1(t+\varrho,x,y)=f_1(t,x,y)$, $f_2(t+\varrho,x,y)=f_2(t,x,y)$, $g_2(t+\varrho,x,y)=g_2(t,x,y)$, and $({\bf H}^p_{f_1})$ and $({\bf H}^p_{f_2, g_2})$ are satisfied. Then Theorem \ref{peristroapth} implies that
\ce
\lim\limits_{\e\rightarrow 0}\mP\left(\sup\limits_{t\in[0,T]}|X_t^\e-\bar X_t|_\mH\geq\eta\right)=0,
\de
where $\bar X$ is the solution of the corresponding averaged equation:
\ce\left\{\begin{array}{ll}
\dif \bar X_t=[A\bar X_t+B(\bar X_t,\bar X_t)+\frac{1}{\varrho}\int_0^\varrho\phi_1(r)\left(\int_\mH f(\bar X_t,y)\mu^{\bar X_t}_r(\dif y)\right)\dif r]\dif t+\bar X_t\dif W^1_t+\dif \bar K_t,\\
\bar X_t(\xi)\geq 0, \quad \xi\in[0,1],\\
\bar X_t(0)=\bar X_t(1)=0,\quad t\geq 0,\\
\bar X_0=x_0.
\end{array}
\right.
\de

{\bf The almost periodic case.} Assume that $\phi_1(t), \phi_2(t)$ are almost periodic. For example, $\phi_1(t)=\sum\limits_{n=1}^\infty\frac{1}{n^2}\sin(\sqrt{n}t)$, $\phi_2(t)=\sum\limits_{n=1}^\infty\frac{1}{n^2}\cos(\sqrt{n}t)$. Then $({\bf H}^{ap}_{f_1})$ and $({\bf H}^{ap}_{f_2, g_2})$ are satisfied. To verify $({\bf H}^{ex}_{f_2, g_2})$, we take $\mS=\mV$ and $T_n=-\triangle(I-\frac{\triangle}{n})^{-1}$. So, $\{T_n\}_{n\in\mN}$ are continuous on $\mH$. Note that for any $x\in\mH, y\in\mV, t\in\mR$
\ce
_{\mV^*}\left\langle Ay, T_n y\right\rangle_\mV\leq -\l_1\|y\|^2_n, \quad \<\phi_2(t)y, T_n y\>_\mH=\phi_2(t)\|y\|^2_n\leq l_2\|y\|^2_n.
\de
Thus, it holds that
\ce
&&2_{\mV^*}\left\langle Ay, T_n y\right\rangle_\mV+2\<f_2(t,x,y), T_n y\>_\mH+\|g_2(t,x,y)\|_{\cL_2(\mU_2, \mH_n)}^2\\
&\leq&-2\l_1\|y\|^2_n+3l_2\|y\|^2_n+\|y\|^2_n+C|x|^2_\mH\\
&\leq&-2(\l_1-\frac{3}{2}l_2-1)+C|x|^2_\mH.
\de
That is, $({\bf H}^{ex}_{f_2, g_2})$ is satisfied. Then Theorem \ref{aperistroapth} yields that
\ce
\lim\limits_{\e\rightarrow 0}\mP\left(\sup\limits_{t\in[0,T]}|X_t^\e-\bar X_t|_\mH\geq\eta\right)=0,
\de
where $\hat X$ is the solution of the corresponding averaged equation:
\ce\left\{\begin{array}{ll}
\dif \hat X_t=[A\hat X_t+B(\hat X_t,\hat X_t)+\hat f_{1,ap}(\hat X_t)]\dif t+\hat X_t\dif W^1_t+\dif \hat K_t,\\
\hat X_t(\xi)\geq 0, \quad \xi\in[0,1],\\
\hat X_t(0)=\hat X_t(1)=0,\quad t\geq 0,\\
\hat X_0=x_0,
\end{array}
\right.
\de
and for any $x\in\mH$
$$
\hat f_{1,ap}(x)=\lim\limits_{T\rightarrow\infty}\frac{1}{T}\int_0^{T}\int_\mH \phi_1(s)f(x,y)\mu^x_s(\dif y)\dif s.
$$

Finally, we mention that if system (\ref{exeq}) has no the fast part, it falls into the framework in \cite{Q2}. Moreover, our averaging principle in the asymptotic case is the same to Theorem 3.7 there.
\ex

\end{document}